\documentclass[12pt,a4paper]{article}
\usepackage[T2A]{fontenc}
\usepackage[cp1251]{inputenc}
\usepackage{amsmath,amssymb,amsthm,amsfonts,amscd}
\begin{document}

\title{Geometric approach to stable
homotopy groups of spheres. I. The Hopf invariant }
\author{P.M. Akhmet'ev}
\date{dedicated to the memory of Professor M.M. Postnikov}

\sloppy \theoremstyle{plain}
\newtheorem{theorem}{Theorem}
\newtheorem*{main*}{Main Theorem}
\newtheorem*{theorem*}{Theorem}
\newtheorem{lemma}[theorem]{Lemma}
\newtheorem{proposition}[theorem]{Proposition}
\newtheorem{corollary}[theorem]{Corollary}
\newtheorem{conjecture}[theorem]{Conjecture}
\newtheorem{problem}[theorem]{Problem}

\theoremstyle{definition}
\newtheorem{definition}[theorem]{Definition}
\newtheorem{remark}[theorem]{Remark}
\newtheorem*{remark*}{Remark}
\newtheorem*{example*}{Example}
\newtheorem{example}[theorem]{Example}

\def\E{{\bf H}}
\def\II{{\dot{\bf I}}}
\def\Z{{\Bbb Z}}
\def\R{{\Bbb R}}
\def\RP{{\Bbb R}\!{\rm P}}
\def\N{{\bf N}}
\def\L{{\bf L}}
\def\C{{\Bbb C}}
\def\A{{\bf A}}
\def\D{{\bf D}}
\def\Q{{\bf Q}}
\def\i{{\bf i}}
\def\j{{\bf j}}
\def\k{{\bf k}}
\def\H{{\bf E}}
\def\F{{\bf F}}
\def\J{{\bf J}}
\def\G{{\bf G}}
\def\I{{\bf I}}
\def\e{{\bf e}}
\def\b{{\bf b}}
\def\a{a}
\def\bb{{\dot{b}}}
\def\f{{\bf f}}
\def\d{{\bf d}}
\def\H{{\bf H}}
\def\fr{{\operatorname{fr}}}
\def\st{{\operatorname{st}}}
\def\mod{{\operatorname{mod}\,}}
\def\cyl{{\operatorname{cyl}}}
\def\dist{{\operatorname{dist}}}
\def\sf{{\operatorname{sf}}}
\def\dim{{\operatorname{dim}}}
\def\dist{\operatorname{dist}}

\maketitle

\begin{abstract}
A geometric approach to the stable homotopy groups of spheres is
developed in this paper, based on the Pontryagin-Thom
construction. The task of this approach is to obtain an
alternative proof of the Hill-Hopkins-Ravenel theorem [H-H-R] on
Kervaire invariants in all dimensions, except, possibly, a finite
number of dimensions. In the framework of this approach, the Adams
theorem on the Hopf invariant is studied, for all dimensions with
the exception of $15$, $31$, $63$, $127$.  The new approach is based on
the methods of geometric topology.
\end{abstract}

\section*{Introduction}
Let $\pi_{n+m}(S^m)$ be the homotopy groups of spheres. Under the
condition $m \ge n + 2$ this group is independent of $m$ and is
denoted by $\Pi_n$. It is called the stable homotopy group of
spheres in dimension $n$. The problem of calculating the stable
homotopy groups of spheres is one of the main unsolved problems of
topology. A development of the Pontryagin-Thom construction 
leads to various applications having important practical
significance: an approach by V.I.Arnol'd to bifurcations of
critical points in  multiparameter families of functions [V]
chapter 3,  section 2.2 and Theorem 1, section 2.4, Lemma 4,
and gives many unsolved
geometrical problems [E2].

For the calculation of elements of the stable homotopy groups of
spheres, one frequently studies algebraic invariants which are
defined for all dimensions at once (or for some infinite sequence
of dimensions). Nevertheless, as a rule these invariants turn out
 to be trivial, and are nonzero only in exceptional cases,
see [M1]. As Prof. Peter Landweber noted: "`This a very interesting
"`philosophyпї"'. Are there examples to illustrate this, apart from
the Hopf invariant and the Kervaire invariant? There might be one
in N. Minami paper [M2]."'

 A basic invariant is the Hopf invariant, which is defined as
follows in the framework of stable homotopy theory.  The Hopf
invariant (also called the stable Hopf invariant), or the
Steenrod-Hopf invariant is a homomorphism
$$ h: \Pi_{2k-1} \to \Z/2, $$
for details see [W],[M-T]. The stable Hopf invariant is studied in
this paper.

The following theorem was proved by J.F.Adams in [A].

\begin{theorem*}
The stable Hopf invariant $h: \Pi_n \to \Z/2$, $n \equiv 1 \pmod{2}$ is a trivial homomorphism
if and only if $n \ne 1, 3, 7$.
\end{theorem*}

\begin{remark*}
The case n = 15 was proved by Toda (cf. [M-T] Ch.
18).
\end{remark*}

Later Adams and Atiyah offered an alternative approach to the
study of the Hopf invariant, based on results of $K$-theory and
the Bott periodicity theorem, cf. [A-A]. This approach was also
extended in subsequent works. A simple proof of the theorem of
Adams, close to the proof of Adams and Atiyah, was given by V.M.
Buchstaber in [B], section 2.

 The
definition of the stable Hopf invariant is reformulated in the
language of the cobordism groups of immersions of manifolds [E1,
K2, K-S1, K-S2, La]. Using the Pontryagin-Thom
theorem in the form of Wells on the representation of the
stable homotopy groups of infinite dimensional real projective
space (which by the Kahn-Priddy theorem surject onto the
2-components of the stable homotopy groups of spheres), we
classify the cobordism of immersions of (in general nonorientable)
manifolds in codimension 1. The Hopf invariant is expressed as a
characteristic number of  the manifold of double
points of self-intersection of an immersion of a manifold
representing the given element of the stable homotopy group.  This is explicitly formulated in
[E1], Lemma 3.1. This lemma is reformulated in the standard way by
means of the Pontryagin-Thom construction for immersions.

The Theorem of Adams admits a simple geometric proof for
dimensions $n \ne 2^{\ell}-1$. In the case $n \not\equiv 3
\pmod{4}$ a proof, using the elements of the theory of immersions,
was given by A. Sz$\ddot{\rm{u}}$cs in [Sz]. The next case in
complexity arises for $n \ne 2^{\ell}-1$. The proof of AdamsпїЅ
Theorem under this assumption was given by Adem [Adem] using
algebraic methods. In this paper the Adem relations on the
multiplicative generators of the Steenrod algebra were used. The
theorem of Adem was reproved using geometric methods in a joint
paper of the author and A. Sz$\ddot{\rm{u}}$cs [A-Sz].

We assume below that $n=2^{\ell}-1$.
Define a positive integer
 $\sigma=\sigma(\ell)$ by the formula:
\begin{eqnarray}\label{sigma}
\sigma = \left[\frac{\ell}{2}\right]-1.
\end{eqnarray}
In particular, for $\ell=8$, $\sigma=3$. Denote $n_{s} = 2^{s}-1$.
Assume that $s$ is a positive integer, then $n_{s}$ is a positive
integer.

The following is the main result of Part $I$.

\subsubsection*{Main Theorem}
Assume that $\ell \ge 8$, therefore $n \ge 255$. Let $g:
M^{\frac{3n+n_{\sigma}}{4}} \looparrowright \R^n$ be an arbitrary
smooth immersion of a closed manifold $M$,
$\dim(M)=\frac{3n+n_{\sigma}}{4}$, where the normal bundle
$\nu(g)$ to the immersion $g$ is isomorphic to the Whitney sum of
$(\frac{n-n_{\sigma}}{4})$ copies of a line bundle $\kappa$ over
$M$,
 $\nu(g)=(\frac{n-n_{\sigma}}{4})$. (In particular, $w_1(M)=0$, where $w_1$ is the
first Stiefel-Whitney class, because the codimension of the
immersion $g$ is even and
$w_1(M)=(\frac{n-n_{\sigma}}{4})w_1(\kappa)=0$), in general $M$ is
nonconnected.)
 Then  the equation $\langle
w_1(\kappa)^{ \dim(M)};[M] \rangle=0$ is valid.
\[  \]

The Main Theorem is deduced from Theorem $\ref{th9}$. Theorem
$\ref{th9}$ is deduced from Propositions  $\ref{prop22}$,
$\ref{prop23}$; these propositions follow from Lemmas
$\ref{osnlemma1}$ and $\ref{osnlemma1}$.

 The
proof of the Main Theorem is based on the principle of geometric
control due to M.Hirsh, see Proposition $\ref{prop24}$. This proposition 
permits one
to find within a cobordism class of immersions an immersion with
additional properties of self-intersection manifold. In this case
we say that the immersion admits a cyclic or quaternionic
structure (see Propositions $\ref{prop22}$,
$\ref{prop23}$).

We can deduce the following from the Main Theorem by standard
arguments.

\subsubsection*{Main Corollary}
Let $g: M^{n-1} \looparrowright \R^n$ be an arbitrary smooth
immersion of the closed manifold $M^{n-1}$, which in general is
not assumed to be orientable. Then under the assumption
$n=2^{\ell}-1$, $n \ge 255$  (i.e., for $\ell \ge 8$), the equality
$\langle w_1(M)^{n-1};[M] \rangle =0$ is valid.
\[  \]
\begin{remark*}
The equivalence of the preceding assertion and the theorem of
Adams (under the restriction $\ell \ge 4$) is proved in [E1,La].
\end{remark*}

We mention that in topology there are theorems which are close to
the formulation of the theorem of Adams. As a rule, these theorems
are corollaries of AdamsпїЅ theorem. Sometimes these theorems can be
given alternative proofs by simpler methods. As S.P. Novikov
remarks in his survey [N], a theorem of this type is the
Bott-Milnor Theorem that the tangent $n$-plane bundle to the
standard sphere $S^n$ is trivial if and only if $n = 1, 3$ or
$7$.  This theorem was first proved in the paper [B-M]. An elegant
modification of the known proof was recently given in [F]. One
should mention the Baum-Browder Theorem [B-B] about non-immersion
of the standard real projective space $\RP^{2^{\ell-1}+1}$ into
$\R^{2^{\ell}-1}$ for $\ell \ge 4$. It would be interesting to
discover an elementary geometrical proof of this theorem and to
prove the Main Corollary for $\ell \ge 4$ as a generalization of
Baum-Browder Theorem.

We turn our attention to the structure of the paper. In section 1
we recall the main definitions and constructions of the theory of
immersions. The results of this section are formally new, but are
easily obtained by known methods. In section 2 the Main Theorem is
reformulated using the notation of section 1 (Theorem
$\ref{th9}$), which represents a basic step in its proof. The
proof of Main Theorem is based on Lemmas $\ref{osnlemma1}$
and $\ref{osnlemma2}$. 

This paper was written for repeated discussions in the seminar on
algebraic topology under the direction of Professor M.M.
Postnikov.

The author is grateful for the discussions to Prof.
V.M. Buchstaber, Prof. S.A.Bogatyi, Prof. A.V. Chernavsky, Prof.
V. V. Chernov, Prof. P. J. Eccles,
  Prof. P. Landweber, Prof. A.S. Mishenko, Prof. Th. Yu. Popelelsky,  
 Prof. O.Saeki, Prof. E.V. Scepin, Prof. A.B. Skopenkov, Prof. Yu. P. Soloviev,
 Prof. V.A. Vassiliev, N. Brodsky, S.A.Melikhov, R.R.Sadykov and M.B. Skopenkov.
Let me mention that Prof. Peter Landweber and Prof. A.S. Mishenko has devoted time and
care to part 1 of this paper on the Hopf invariant and part 2 of
this paper on the Kervaire invariant. Prof. Peter Landweber pointed out (see [L]) that
the approach in an earlier version of this paper is not valid.
The mistake occurs in section 3 of paper [A1]. I am pleased to acknowledge
his effort.

\section{Preliminary information}
We recall the definition of the cobordism groups of framed immersions in
Euclidean space, which is a special case of a more general construction presented
in the book [K1] on page 55 and in section 10. The connection
with the Pontryagin-Thom construction is explained in [A-E].

Let $f:  M^{n-1}  \looparrowright  \R^{n}$  be a smooth immersion,
where the $(n-1)$-dimensional manifold $M^{n-1}$ is closed, but in
general is nonorientable and nonconnected. We introduce a relation
of cobordism on the space of such immersions. We say that two
immersions $f_0$, $f_1$ are connected by a cobordism, $f_0  \sim
f_1$, if there exists an immersion $ \Phi  :(  W^{n},  \partial  W
=  M^{n-1}_0  \cup M^{n-1}_1) \looparrowright  (\R^{n}  \times
[0;1];\R^n \times \{0,1\})$ satisfying the boundary conditions
$f_i = \Phi \vert_{ M^{n-1}_i} : M^{n-1}_i \looparrowright \R^{n}
\times \{i\} $, $i=0,1$ and, moreover, it is required that the
immersion $\Phi$ is orthogonal to $\R \times \{0,1\}$.

The set of cobordism classes of immersions forms an Abelian group
with respect to the operation of disjoint union of immersions. For
example, the trivial element of this group is represented by an
empty immersion, the element that is inverse to a given element
represented by an immersion $f_0$ is represented by the
composition $S \circ f_0$, where $S$ is a mirror symmetry of the
space $\R^n$.

This group is denoted by $Imm(n-1,1)$. Because $\RP^{\infty}
= MO(1)$, by the Wales theorem [Wa] (see [E1],[Sz2] for references) this
group maps onto the stable homotopy group of
spheres $\lim_{k \to \infty} \pi_{n+k}(S^k)$.

The immersion $ f $ defines an isomorphism of the normal bundle to
the manifold $ M ^{n-1} $ and the orientation  line bundle $
\kappa $, i.e. an isomorphism $ D (f): T (M ^{n-1}) \oplus \kappa
\cong n \varepsilon $, where $ \varepsilon $ is the trivial line
bundle over $ M ^{n-1} $. In similar constructions in surgery
theory of smooth immersions  one requires  a stable isomorphism of
the normal bundle of a manifold $ M ^{n-1} $ and the orientation
line bundle $ \kappa $, i.e. an isomorphism $ T (M ^{n-1}) \oplus
\kappa \oplus N \varepsilon \cong (n + N) \varepsilon $, for $ N>0
$. Using Hirsch's Theorem [Hi], it is easy to verify that if two
immersions of $ f_1, f_2 $ determine isomorphisms $ D (f_1) $, $ D
(f_2) $, that belong to the same class of stable isomorphisms then
the immersions $ f_1 $, $ f_2 $ are regularly cobordant and even
regularly concordant (but, generally speaking, may not be
regularly homotopic).

We also require groups $Imm^{sf}(n-k,k)$. An element of this group is represented
by a triple $(f, \kappa, \Xi)$, where $f: M^{n-k}
\looparrowright \R^{n}$ is an immersion of a
closed manifold, $\kappa: E(\kappa) \to M^{n-k}$ is a line bundle (to shorten notation, we
shall denote the line (one-dimensional) bundle and its characteristic class in $H^1(M^{n-k};\Z/2)$
by the same symbol), and $\Xi$ is a skew-framing of the normal
bundle of the immersion by means of the line bundle $\kappa$, i.e., an isomorphism
of the normal bundle of the immersion $f$ and the bundle $k \kappa$. In the
case of odd $k$, the line bundle $\kappa$ turns out to be orientable over $M^{n-k}$ and
necessarily $\kappa=w_1(M^{n-k})$.

 Two elements of the cobordism group, represented
by triples $(f_1,\kappa_1,\Xi_1)$, $(f_2,\kappa_2,\Xi_2)$ are
equal if the immersions $f_1$, $f_2$ are cobordant (this
definition is analogous to the previous one for representatives of
the group $Imm^{sf}(n - 1, 1)$), where in addition it is required
that the immersion of the cobordism be skew-framed, and that the
skew-framing of the cobordism be compatible with the given
skew-framings on the components of the boundary. We remark that
for $k = 1$ the new definition of the group $Imm^{sf}(n-k,k)$
coincides with the original definition.

We define a homomorphism $$J^{sf}: Imm^{sf}(n-1,1) \to
Imm^{sf}(n-k,k),$$ which is called the homomorphism of transition
to codimension $k$. Consider a manifold $M'^{n-1}$ and an
immersion $f': M'^{n-1} \looparrowright \R^n$ representing an
element of the first group, and consider a classifying map
$\kappa': M'^{n-1} \to \RP^a$ to a real projective space of large
dimension ($a = n-1$ suffices), representing the first
Stiefel-Whitney class $w_1(M'^{n-1})$. Consider the standard
subspace $\RP^{a-k+1} \subset \RP^{a}$ of codimension $(k - 1)$.
Assume that the mapping $\kappa'$ is transverse along the chosen
subspace and define the submanifold $M^{n-k} \subset M'^{n-1}$ as
the complete inverse image of this subspace for our mapping,
$M^{n-k} = \kappa'^{-1}(\RP^{a-k+1})$. Define an immersion $f:
M^{n-k} \looparrowright \R^n$ as the restriction of the immersion
$f'$ to the given submanifold. Notice that the immersion $f:
M^{n-k} \looparrowright \R^n$ admits a natural skew-framing. In
fact, the normal bundle to the submanifold $M^{n-k} \subset
M'^{n-1}$ is naturally isomorphic to the bundle $(k-1) \kappa$,
where $\kappa =\kappa' \vert_{M}$ (here and below, when a manifold
is used as a subscript,  the superscript indicating the dimension
of the manifold is omitted). The isomorphism  $\Xi$ is defined by
the standard skew-framing of the normal bundle to the submanifold
$\RP^{a-k+1}$ of the manifold $\RP^a$, which is transported to the
submanifold $M^{n-k} \subset M'^{n-1}$, since it is assumed that
$\kappa'$ is transverse regular along $\RP^{a-k+1}$. A further
direct summand in the skew-framing of the normal bundle of the
immersion $f$ corresponds to the normal line bundle of the
immersion $f'$. This bundle serves as orientation bundle for
$M'^{n-1}$, hence its restriction to $M^{n-k}$ coincides with
$\kappa$. The homomorphism $J^{sf}$ carries the element
represented by the immersion $f'$ to the element represented by
the triple $(f,\kappa,\Xi)$. Elementary geometrical
considerations, using only the concept of transversality imply
that the homomorphism $ J ^{sf} $ is correctly defined.

We now define the manifold of double points of self-intersection of an
immersion $f: M^{n-k} \looparrowright \R^n$ in general position, and a canonical 2-sheeted
covering over this manifold. Under the assumption that the immersion $f$
is in general position, the subset in $\R^n$ of points of self-intersection of the
immersion $f$ is denoted by $\Delta = \Delta(f)$, $\dim(\Delta)=n-2k$. This subset is defined
by the formula
\begin{eqnarray}\label{Delta}
\Delta=\{ x \in \R^n : \exists x_1, x_2 \in M^{n-k}, x_1 \ne
x_2, f(x_1)=f(x_2)=x\},
\end{eqnarray}  We
define $\bar
\Delta \subset M^{n-k}$ by the formula $\bar \Delta = f^{-1}(\Delta)$.

We recall the standard definition of the manifold of points of
self-intersection and the parameterizing immersion, see e.g. [Ada]
for details.

\subsubsection*{Definition of self-intersection manifold}
The set $N$ is defined by the formula
\begin{eqnarray}\label{N}
N=\{[(x_1,x_2)] \in (M^{n-k} \times M^{n-k})/T' : x_1 \ne x_2,
f(x_1)=f(x_2)\}
\end{eqnarray}
 ($T'$ is the involution
permuting the coordinate factors), and its canonical covering is defined by the
formula
\begin{eqnarray}\label{barN}
\bar N=\{ (x_1,x_2) \in M^{n-k} \times M^{n-k} : x_1 \ne x_2,
f(x_1)=f(x_2)\}.
\end{eqnarray}
Under the assumption that the immersion $ f $ is generic,  $ N $
is a smooth manifold of dimension $ \dim (N) = n-2k $. This
manifold is denoted by $ N^{n-2k} $ and is called the self-intersection manifold of $f$, the projection of the
covering is denoted by $p: \bar N^{n-2k} \to N^{n-2k}$ and is called the canonical 2-sheeted covering.
\[  \]

The immersion $\bar g: \bar N^{n-2k} \looparrowright M^{n-k}$,
parameterizing $\bar \Delta$, is defined by the formula $\bar g =
p \vert_{\bar N}$. The immersion  
\begin{eqnarray}\label{g}
g: N^{n-2k} \looparrowright \R^n, 
\end{eqnarray}
which is a parametrization of $\Delta$, is defined by the formula
$g([x_1,x_2])=f(x_1)$. Notice that the parameterizing immersion $g:
N^{n-2k} \looparrowright \R^n$ of $\Delta$ in general, is not an
immersion in general position. There is a two-sheeted covering $p:
\bar N^{n-2k} \to N^{n-2k}$, for which $g \circ p = f \circ \bar
g$. This covering is called the canonical covering over the
manifold of points of self-intersection.

\begin{definition}\label{def1}
Let $(f, \kappa, \Xi)$ represent an element in the group
$Imm^{sf}(n-k,k)$. Let us define the homomorphism:
$$ h_k : Imm^{sf}(n-k,k) \to \Z/2,$$
called the stable Hopf invariant by the following formula:
$$h_k([f,\kappa,\Xi]) =
\langle \kappa^{n-k};[M^{n-k}] \rangle. $$
\end{definition}

The definitions of the stable Hopf invariant (in the sense of
Definition 1) for distinct values of $k$ are compatible with one
another, and coincide with the definition used in the
introduction. We formulate this as a separate assertion.

\begin{proposition}\label{prop2}
The homomorphism $J^{sf}: Imm^{sf}(n-1,1) \to Imm^{sf}(n-k,k)$
preserves the Hopf invariant, i.e. the invariant $h_1:
Imm^{sf}(n-1,1) \to \Z/2$ and the invariant  $h_k: Imm^{sf}(n-k,k)
\to \Z/2$ are related by the formula:
\begin{eqnarray}\label{1}
 h_1 = h_k \circ J^{sf}.
\end{eqnarray}
\end{proposition}

\subsubsection*{Proof of Proposition $\ref{prop2}$}

Let $f: M^{n-k} \looparrowright \R^n$ be an immersion with a
skew-framing $\Xi$ of its normal bundle and with characteristic
class $\kappa \in H^1(M^{n-k};\Z/2)$, representing an element of
$Imm^{sf}(n-k,k)$, which satisfies $J^{sf}([f'])=[f,\kappa,\Xi]$
for an element $[f'] \in Imm^{sf}(n-1,1)$, where $f$ and $f'$ are
related as in the definition of $J^{sf}$. By definition,
$h_k([f,\kappa,\Xi])=\langle \kappa^{n-k};[M^{n-k}] \rangle$.

On the other hand, $M^{n-k} \subset M'^{n-1}$ is a cycle dual in
the sense of Poincar\'{e} to the cohomology class $ \kappa '^{k-1}
\in H ^{k-1} (M '^{n-1}; \Z / 2) $. The formula $(\ref{1})$ is
valid, since $\langle \kappa'^{n-1};[M'^{n-1}]\rangle =\langle
\kappa^{n-k};[M^{n-k}]\rangle $. In the case $k=n$ this formula is also satisfied. 
Proposition $\ref{prop2}$ is proved.

Let us present an alternative proof.  The image of the
characteristic mapping $\kappa': M'^{n-1} \to \RP^{\infty}$,
without loss of generality
lies in the skeleton of dimension $n-1$ $(n-k)$ of the
classifying space. 
This allows us to write $\kappa': M'^{n-1} \to \RP^{n-1}$.

Moreover, let us assume that the mapping 
 $\kappa'$ is transversal along  $\RP^{n-k} \subset \RP^{n-1}$, i.e.
$M^{n-k} \subset M'^{n-1}$ is determined by the formula  $M^{n-k} = \kappa'^{-1}(\RP^{n-k} \subset \RP^{n-1})$, $\kappa = \kappa' \vert_{M^{n-k}}$. 
The marked point  $pt \in \RP^{n-k} \subset \RP^{n-1}$ is a regular value for $\kappa$.

The characteristic number $\langle
\kappa'^{n-1};[M'^{n-1}]\rangle$ ($\langle
\kappa^{n-k};[M^{n-k}]\rangle$)  coincides with the degree
$\deg(\kappa')$ ($\deg(\kappa)$) of the classifying map $\kappa':
M'^{n-1} \to \RP^{n-1}$ ($\kappa: M^{n-k} \to \RP^{n-k}$),
 which is considered modulo 2 and is determined
as the parity of the number of preimages of regular values of this
map. The value $\deg(\kappa')$ ($\deg(\kappa)$) does not depend on
the choice of the mapping  $\kappa'$ ($\kappa$) to the skeleton,
as above.  The degrees
$\deg(\kappa')$ and $\deg(\kappa)$ coincide, since the regular
value can be chosen common to these mappings.
Proposition $\ref{prop2}$ is alternatively 
proved.
\[  \]

\[  \]

Let us formulate another (equivalent) definition of the stable
Hopf invariant (assuming $ n-2k> 0 $).

\begin{definition}\label{proba}
Let $(f, \kappa, \Xi)$ represent an element in the group
$Imm^{sf}(n-k,k)$, $n-2k>0$.
 Let $N^{n-2k}$ be the manifold of the double points of the immersion $f: M^{n-k} \looparrowright \R^n$,
 $\bar N$ be the canonical 2-sheeted cover over
$N$, $\kappa_{\bar N} \in H^1(\bar N;\Z/2)$ be induced from
$\kappa \in H^1(M^{n-k};\Z/2)$ by the immersion $\bar g: \bar
N^{n-2k} \looparrowright M^{n-k}$.

Let us define the homomorphism $h_k^{sf} : Imm^{sf}(n-k,k) \to \Z/2$ by
the formula:
$$h_k^{sf}([f,\kappa,\Xi]) = \langle \kappa^{n-2k}_{\bar
N};[\bar N^{n-2k}] \rangle. $$
\end{definition}
\[  \]

The following proposition establishes the equivalence of
Definitions $\ref{def1}$ and $\ref{proba}$.

\begin{proposition}\label{prop4}
Let us assume that the conditions of Definition $\ref{proba}$ are
satisfied. Then we have:
\begin{eqnarray}\label{1.400}
\langle \kappa^{n-2k}_{\bar N};[\bar N^{n-2k}] \rangle =\langle
\kappa^{n-k};[M^{n-k}]\rangle.
\end{eqnarray}
\end{proposition}

\subsubsection*{Proof of Proposition $ \ref{prop4} $}

Let $f: M^{n-k} \looparrowright \R^n$ be an immersion  with a skew
framing $ \Xi $  and with the characteristic class  $\kappa \in
H^1(M^{n-k};\Z/2)$; then the triple $(f,\Xi,\kappa)$ represents an
element in the group $Imm^{sf}(n-k,k)$. Let $ N ^{n-2k} $ be the
manifold of self-intersection points of the immersion $ f $, $g:
N^{n-2k} \looparrowright \R^n$ be the parameterizing immersion,
and $\bar N^{n-2k} \to N^{n-2k}$ be the canonical double covering.
Consider the image of the fundamental class $\bar g_{\ast}([\bar
N^{n-2k}]) \in H_{n-2k}(M^{n-k};\Z/2)$ by the immersion $\bar g:
\bar N^{n-2k} \looparrowright M^{n-k}$ and denote by $m \in
H^{k}(M^{n-k};\Z/2)$ the cohomology class Poincar\'{e} dual to the
homology class $\bar g_{\ast}([\bar N^{n-2k}])$. Consider also the
cohomology Euler class of the normal bundle immersion $ f $, which
is denoted by $e \in H^k(M^{n-k};\Z/2)$.

By the Herbert Theorem for the immersion $f: M^{n-k} \looparrowright
\R^n$ with self-intersection manifold $N^{n-2k}$ (see [E-G],
Theorem 1.1 the case $r=1$ coefficients is $\Z/2$; see also this
theorem in the original papers [He], [L-S]) the following formula
is valid:
\begin{eqnarray}\label{1.500}
e+m=0.
\end{eqnarray}

Since the Euler class $ e $ of the normal bundle $ k \kappa $ of
the immersion $ f $ is equal to $ \kappa^k $ (line bundles and
their corresponding  characteristic cohomology classes are denoted
by the same symbols), then the cycle $\bar g_{\ast}([\bar N]) \in
H_{n-2k}(M^{n-k};\Z/2)$ is Poincar\'{e} dual to the cohomology
class $\kappa^{k} \in H^{k}(M^{n-k};\Z/2)$. Therefore, the formula
$(\ref{1.400})$ and Proposition $\ref{prop4}$ are proved.
\[  \]

\[  \]
It is more convenient to reformulate Proposition $\ref{prop4}$ (in
a more general form) by means of the language of commutative
diagrams. The desired reformulation is given in Lemma
$\ref{lemma7}$ below. We turn to the relevant definitions.

Let $g: N^{n-2} \looparrowright \R^n$ be the immersion  of the
double self-intersection points of the immersion $f: M^{n-1}
\looparrowright \R^n$ of codimension $ 1 $. We denote by $\nu_N :
E(\nu_N) \to N^{n-2}$ the normal 2-dimensional bundle of the
immersion $ g $. (Note that the disk bundle associated with the
vector bundle $ \nu_N $ is diffeomorphic to a regular closed
tubular neighborhood of the immersion $ g $.)

In comparison with an arbitrary vector bundle, this bundle carries
an additional structure, namely its structure group as an
$O(2)$-bundle admits a reduction to a discrete dihedral group
which we denote by $\D$. This group has order 8, and is defined as
the group of orthogonal transformations of the plane which carry
the standard pair of coordinate axes into themselves (with
possible change of orientation and order).

In the standard presentation of the group $\D$ there are two
generators $a$, $b$ which are connected by the relations $\{a^4 =
b^2 = 1, [a,b] = a^2 \}$. The generator $a$ is represented by the
rotation of the plane through an angle $\frac{\pi}{2}$, and the
generator $b$ is represented by a reflection with respect to the
bisector of the first and second coordinate axes. Notice that the
element $ba$ (the product means the rule of composition $b \circ
a$ of transformations in $O(2))$ is represented by the reflection
with respect to the first coordinate axis.


\subsubsection*{The structure group of the normal bundle of the manifold
of self intersection points for an immersion $f: M^{n-k} \looparrowright
\R^n$, in case $k=1$}

Let us use the transversality condition for the immersion $f:
M^{n-1} \looparrowright \R^n$. Let $N^{n-2}$ be the
self-intersection manifold of the immersion $ f $, $g: N^{n-2}
\looparrowright \R^n$ be the parameterizing immersion. In the
fiber  $ E (\nu_N) _x $ of the normal bundle $ \nu_N $ over the
point $ x \in N $ an unordered pair of axes is fixed. These axes
are formed by the tangents to the curves of intersection of the
fiber $E(\nu_N)_x$ with two sheets of immersed manifolds
intersecting transversely in the neighborhoods of this point. By
construction, the bundle $ \nu_N $ has the structure group $ \D
\subset O (2) $.

Over the space $ K (\D, 1) $ the universal 2-dimensional $ \D
$--bundle is defined. This bundle will be denoted by $ \psi: E
(\psi) \to K (\D, 1) $. We say that the mapping $ \eta: N \to K
(\D, 1) $ is classifying for the bundle $ \nu_N $, if an
isomorphism $\Xi: \eta^{\ast}(\psi) \cong \nu_N$ is well defined,
where $ \eta ^{\ast} (\psi) $ is the inverse image of the bundle $
\psi $ and $ \nu_N $ is the  normal bundle of the immersion $ g $.
Further a bundle itself and its classifying map will be denoted
the same; in the considered case we have $ \eta \cong \nu_N $. The
isomorphism $ \Xi $ will be called a $ \D $--framing of the
immersion $ g $, and the mapping $ \eta $ will be called the
characteristic mapping of the $ \D $--framing $ \Xi $.


Let a triple $(f,\kappa,\Xi)$, where $f: M^{n-k}
\looparrowright \R^n$ is an immersion and $\Xi$ is a skew--framing of $f$ with the characteristic class
$\kappa \in H^1(M;\Z/2)$, represent an element of the group $Imm^{sf}(n-k,k)$.
We need to generalize the previous construction
for $k = 1$ and to describe the structure group of the normal bundle $\nu_N$
to the manifold $N^{n-2k}$ of self--intersection points of a generic immersion of an arbitrary codimension
$k$, $k \le
[\frac{n}{2}]$.

\begin{proposition}\label{prop5}
The normal bundle $\nu_{N}$, $\dim(\nu_{N}) = 2k$, of the
immersion $g$ is a direct sum of $k$ copies of a two-plane bundle
$\eta$ over $N^{n-2k}$, where each two-plane bundle has structure
group $\D$ and is classified by a classifying mapping $\eta:
N^{n-2k} \to K(\D,1)$ $($ an analogous proposition is proved in
{\rm[Sz2]}$)$.
\end{proposition}

\subsubsection*{Proof of Proposition $ \ref{prop5} $}

Let $x \in N^{n-2k}$ be a point in the manifold of double points. Denote by
$\bar x_1, \bar x_2 \in \bar N^{n-2k} \looparrowright M^{n-k}$ the two preimages of this point under the canonical covering map
by the double covering.
The orthogonal complement in the space $T_{g(x)}(\R^n)$ to the
subspace $g_{\ast}(T_x(N^{n-2k}))$  is the fiber of the
normal bundle $E(\nu_N)$ of the immersion $ g $ over a point $ x \in
N^{n-2k} $. This fiber is represented as a direct sum of two linear
spaces, $E(\nu_N)_{x} = \bar E_{x,1} \oplus \bar E_{x,2}$, where
each subspace $\bar E_{x,i} \subset E(\nu_N)_x$ is a fiber of
the normal bundle of the immersion $ f $ at the point $ \bar x_i $.

Each subspace $\bar E_{x,i}$ of the fiber is canonically a direct
sum of $k$ ordered copies of the fiber of a line bundle, since the
normal bundle to the immersion $f$ is equipped with a
skew-framing. We group the fibers $\bar E(\kappa_{x,j,i})$, $j=1,
\dots, k$, $i=1,2$ with a corresponding  index
into a
two-dimensional subfiber of the fiber of the normal bundle
$E(\nu_N)_x$. As a result, we obtain a decomposition of the fiber
$E(\nu_N)_x$ over each point $x \in N^{n-2k}$ into a direct sum of
$k$ copies of a two-dimensional subspace. This construction
depends continuously on the choice of the point $x$, and can be
carried out simultaneously for each point of the base $N^{n-2k}$.
As a result, we obtain the required decomposition of the bundle
$\nu_N$ into a direct sum of a number of canonically isomorphic
two-plane bundles. Each two-dimensional summand is classified by a
structure  map $\eta: N \to K(\D,1)$, which proves Proposition
$\ref{prop5}$.
\[  \]

\begin{definition}\label{def6}
We define the cobordism group of immersions $ Imm ^{\D} (n-2k, 2k)
$, assuming $ n> 2k $. Let $ (g, \eta, \Psi) $ be a triple, which
determines a $ \D $--framed immersion of codimension $ 2k $. Here
$ g: N^{n-2k} \looparrowright \R ^n $ is an immersion and $ \eta:
N ^{n-2k} \to K (\D, 1) $ is the classifying map of the $ \D
$--framing $ \Psi $. The cobordism relation of triples is
standard.
\end{definition}

\begin{lemma}\label{lemma7}
Under the assumption $k_1 < k$, $2k < n$, the following
commutative diagram of groups is well defined:
\begin{eqnarray}\label{1.5}
\begin{array}{ccccc}
Imm^{sf}(n-k_1,k_1) & \stackrel{J^{sf}}{\longrightarrow} & Imm^{sf}(n-k,k) & \stackrel{h_k^{sf}}{\longrightarrow} & \Z/2\\
\downarrow \delta_{k_1}& &\downarrow \delta_{k} && \| \\
Imm^{\D}(n-2k_1,2k_1) & \stackrel{J^{\D}}{\longrightarrow} &
Imm^{\D}(n-2k,2k) & \stackrel{h^{\D}_k}{\longrightarrow} &
\Z/2.
\end{array}
\end{eqnarray}
\end{lemma}

\subsubsection*{Proof of Lemma $ \ref{lemma7} $}

We define the homomorphisms in the diagram $ (\ref{1.5}) $. The
homomorphism
$$Imm^{sf}(n-k_1,k_1)  \stackrel{J_{sf}}{\longrightarrow}
Imm^{sf}(n-k,k)$$
  is defined exactly  was the
homomorphism $J_{sf}$  for the case $ k_1 = 1 $.

Define the homomorphism
\begin{eqnarray}\label{JD}
Imm^{\D}(n-2k_1,2k_1)   \stackrel{J^{\D}}{\longrightarrow}
Imm^{\D}(n-2k,2k)
\end{eqnarray} 

Let us present 3 (equivalent) definitions.
Assume that a triple  $[(g', \eta', \Psi')]$ represents an element in the cobordism group  $Imm^{\D}(n-2k,2k)$, where $g: N'^{n-2k} \looparrowright
\R^n$--is an immersion, which is equipped with the dihedral framing.

--1.  Take the universal $\D$--bundle $\psi$ over $K(\D,1)$ and take the
pull-back  $\eta'^{\ast}(\psi)$ of this bundle by means of the classifying map $\eta'$. 
Take a submanifold $N^{n-2k_1} \subset
N'^{n-2k}$, which represents the Euler class of the bundle
$(k-k_1)\eta'^{\ast}(\psi)$. The triple $(g, \eta, \Psi)$ is
well defined, where $g = g' \vert_{N}$ and $\eta = \eta'
\vert_{N}$. The $\D$-framing $\Psi$ is defined below.

Let us consider  the normal bundle $\nu_{g'}$ of the immersion
$g'$. The restriction of this bundle on the submanifold 
$N^{n-2k} \subset N'^{n-2k_1}$
is decomposed into the Whitney sum of the two
bundles: $\nu_{g'}=\nu_g \vert_N  \oplus \nu_{N \subset N'}$,
where $\nu_{N \subset N'}$ is the normal bundle of the submanifold
$N^{n-2k} \subset N'^{n-2k_1}$. The bundles $\nu_{N \subset N'}$
and $(k-k_1)\eta^{\ast}(\psi)$ are isomorphic and this bundle is
equipped with the standard $\D$--framing.  Therefore the
bundle $\nu_{g}$ is equipped with the dihedral framing $\Psi:
\nu_{g} \cong k \eta^{\ast}(\psi)$. The triple $(g, \eta,
\Psi)$ represent the required element $J^{\D}(g', \eta', \Psi') \in
Imm^{\D}(n-2k,2k)$.

--2. Consider the configuration space  $Sym_2(\RP^s)$ of non-ordered pairs of distinguished points in  $\RP^s$, 
let us denote this configuration space by 
$\Gamma_{\circ}(s)$ (see the formula $(\ref{99})$, in which we assume that $s=n-k$) in the case $s > 2(n-k_1)+1$. The space 
 $\Gamma_{\circ}(s)$ is an open $2s$--dimensional manifold, with the homotopy $n-k$-type of the space $K(\D,1)$.
Without loss of a generality we may assume that the characteristic mapping is the following: 
$\eta': N'^{n-2k_1} \to \Gamma_{\circ}(s)$. Let us consider the submanifold
$\Gamma_{\circ}(s-k+k_1) \subset \Gamma_{\circ}(s)$ of the codimension $2(k-k_1)$, which is induced by the standard inclusion  $\RP^{s-k+k_1} \subset \RP^s$. Assume, without loss of the generality, that the mapping   $\eta'$ is transversal along the submanifold 
$\Gamma_{\circ}(s-k+k_1) \subset \Gamma_{\circ}(s)$. Define a submanifold  $N^{n-2k} \subset N'^{n-2k_1}$ by the formula 
$N^{n-2k} = \eta'^{-1}(\Gamma_{\circ}(s-k+k_1))$. Obviously, the normal bundle of the submanifold
 $\Gamma_{\circ}(s-k+k_1) \subset \Gamma_{\circ}(s)$ is isomorphic to the bundle $(k-k_1)\psi$. Therefore the bundles 
  $\nu_{N \subset N'}$ and
$(k-k_1)\eta^{\ast}(\psi)$ are isomorphic. Let us define a triple  $(g, \eta, \Psi)$, which represents the element
   $J^{\D}(g', \eta', \Psi') \in Imm^{\D}(n-2k,2k)$, analogously to the Definition --1.

--3.  Let us consider the canonical  2-sheeted covering $\bar N'^{n-2k_1} \to N'^{n-2k_1}$, which is induced from the universal 2-sheeted 
covering 
$K(\I_c,1) \to K(\D,1)$ by the mapping $\eta'$ (the subgroup $\I_c \subset \D$ of the index 2 is defined below by the formula  $(\ref{1.6})$).
Let us denote the involution of this 2-sheeted covering by  $T: \bar N'^{n-2k_1} \to \bar N'^{n-2k_1}$. 
Consider the submanifold  $\bar W^{n-k-k_1} \subset \bar N'^{n-2k_1}$ of the codimension  $k-k_1$, 
which is defined as the transversal preimage of the submanifold
$\RP^{s-k+k_1} \subset \RP^s$ by the mapping $\bar N'^{n-2k_1} \to \RP^s$ in the standard skeleton  $\RP^s \subset K(\Z/2,1)$.
Define the mapping $\bar N'^{n-2k_1} \to \RP^s$  by the composition of the canonical 2-sheeted covering mapping  $\bar \eta':  \bar N'^{n-2k_1} \to K(\I_c,1)$ over the mapping  $\eta'$ with the mapping  $K(\I_c,1) \to K(\Z/2,1)$, the last mapping is induced by the epimorphism  $(\ref{1.7})$, which is defined below. 
Let us consider the manifold  $\bar W^{n-k-k_1} \cap T(\bar W^{n-k-k_1}) \subset \bar N'^{n-2k_1}$, assuming that the manifolds   $\bar W^{n-k-k_1}$ and
$T(\bar W^{n-k-k_1})$ are transversally intersected inside $\bar N'^{n-2k_1}$. Define the manifold  $N^{n-2k}$ as the quotient 
$(\bar W^{n-k-k_1} \cap T(\bar W^{n-k-k_1}))/T$, this manifold is equipped by the natural embedding into the manifold  $N'^{n-2k_1}$.
The normal bundle of the submanifold  $N^{n-2k} \subset N'^{n-2k_1}$ is naturally isomorphic to the Whitney sum of  $k-k_1$ copies of 
a 2-dimensional $\D$--bundle with the classifying mapping 
 $\eta' \vert_{N \subset N'}$. The triple  $(g, \eta, \Psi)$, which represents the element 
   $J^{\D}(g', \eta', \Psi') \in Imm^{\D}(n-2k,2k)$ is defined as in the Definition --1.

Let us prove that Definition -1 and Definition -3 are equivalent.
Let us consider the canonical double covering $p': \bar N'^{n-2k_1}
\to N'^{n-2k_1}$ over the self-intersection points manifold of $g'$.
The manifold $\bar N'^{n-2k_1}$ is naturally immersed into $M^{n-k}$:
$\bar N^{n-2k} \looparrowright M^{n-k}$. Let us consider the
submanifold $M^{n-k} \subset M'^{n-k_1}$, dual to
$\kappa'^{k-k_1}$. This manifold is used in the definition of the homomorphism $J^{sf}$. 
Let us consider the submanifold $M^{n-k} \cap
\bar N'^{n-2k_1} \subset \bar N^{n-2k_1}$,  assuming that $M^{n-k}$
intersects the immersion $\bar N'^{n-2k_1} \looparrowright M^{n-k}$
in a  general position. We denote this submanifold  by $\tilde W$. Obviously, $\dim(\tilde W)=n-k+k_1$
and the codimension of the submanifold
 $\tilde W^{n-k-k_1} \subset \bar N^{n-2k_1}$ is equal to $k-k_1$.

Let us denote by  $T: \bar N^{n-2k} \to \bar
N^{n-2k}$ the involution on the covering space of $p'$. Let us consider the manifold
$T(\tilde W^{n-k+k_1})$  and the intersection $T(\tilde
W^{n-k+k_1}) \cap \tilde W^{n-k+k_1}$ inside $\bar N'^{n-2k_1}$.
Assuming that this intersection is generic,  $T(\tilde
W^{n-k+k_1}) \cap \tilde W^{n-k+k_1}$ is a smooth closed
manifold, let us denote this manifold by $\bar N$,  $\dim(\bar
N)=n-2k$. Moreover, the manifold $\bar N^{n-2k}$ is
equivariant with respect to the involution $T = T'\vert_{\bar
N^{n-2k}}$. The factor-space $\bar N^{n-2k}/T$ is well
defined, this is smooth closed manifold. Let us denote this
manifold by $N^{n-2k}$ and the restriction of the canonical
double cover over this manifold by $p: \bar N^{n-2k} \to
N^{n-2k}$.

Note that the manifold $N^{n-2k}$ is a submanifold in
$N'^{n-2k_1}$ and this submanifold coincides with the
self-intersection manifold of the immersion $f = f' \vert_{M^{n-k} \subset M'^{n-k_1}}$. 
The manifold $N^{n-2k}$ was used in Definition -3 of the homomorphism $J^{\D}$.

Let us denote by $N_2^{n-2k} \subset N'^{n-2k_1}$ the manifold, which represents
the Euler class of the bundle $\eta'^{\ast}(\psi)$. This submanifold is well-defined up to regular $\D$--framed cobordism. Let
us prove that the submanifold $N^{n-2k} \subset N'^{n-2k_1}$ also
represents the Euler class of the bundle $\eta'^{\ast}(\psi)$.
Therefore we may put $N_2^{n-2k}=N^{n-2k}$.

Let us denote the bundle $\eta'^{\ast}(\psi)$ by $\xi'$ for
short. Take the bundle  $p'^{\ast}(\xi')$, let us
denote this bundle by $\bar \xi'$, and take the bundle
$T^{\ast}(p'^{\ast}(\xi')$, which is induced from $\bar \xi'$ by the involution. 
Obviously, the bundle
$\bar \xi'$ decomposes into the Whitney sum of the two
$k-k_1$--dimensional bundles:
\begin{eqnarray}\label{W}
\bar \xi' = \bar \xi'_+ \oplus \bar \xi'_-. 
\end{eqnarray}
In this formula the bundle
 $\bar \xi'_+$ is defined as follows. Take the universal $\D$--bundle $\psi$ and consider $\I_c$--bundle  $\bar \psi$, 
which is induced by 2-sheeted covering constructed by the index 2 subgroup  
$\I_c \subset \D$ (see the formula  $(\ref{1.6})$ below).  Because $\I_c \cong \Z/2 \times \Z/2$, the bundle  $\bar \psi$ 
is naturally isomorphic th the Whitney sum
$\bar \psi = \bar \psi_+ \oplus \bar \psi_-$ of the two line  $\Z/2$--bundles. This Whitney sum induces the Whitney sum $(\ref{W})$.

Analogously, we get 
 $T^{\ast}(\bar \xi') = T^{\ast}(\bar \xi'_+) \oplus T^{\ast}(\bar \xi'_-)$. Moreover, the following equality is satisfied: 
$T^{\ast}(\bar \xi'_+)=\bar \xi'_-$, $T^{\ast}(\bar \xi'_-)=\bar
\xi'_+$. The bundle  $\bar \xi'_+$ is isomorphic to the Whitney sum of  $k-k_1$ copies
of the line bundle   $\kappa'$.

The submanifold $\bar N_2^{n-2k} \subset \bar N'^{n-2k_1}$ is
well defined as the restriction of the covering  $p'$ to the
submanifold $N_2^{n-2k} \subset N'^{n-2k_1}$.
This submanifold represents the equivariant Euler class of the bundle
$\bar \xi'$. This submanifold $\bar N_2^{n-2k}$  is well defined as the intersection
of the two $n-k+k_1$-dimensional submanifolds, which will be defined below and which are denoted by   $\bar W_{2,+}^{n-k+k_1}$ and
  $\bar W_{2,-}^{n-k+k_1}$. The submanifold $\bar W_{2,+}^{n-k+k_1} \subset \bar N'^{n-2k_1}$ 
 represents the Euler class
of the bundle $\bar \xi'_+$. The submanifold $\bar
W_{2,-}^{n-k+k_1} \subset  \bar N'^{n-2k_1}$
represents the Euler class of the
bundle $\bar \xi'_-$. Note that the submanifold 
 $\bar
W_{2,-}^{n-k+k_1} \subset \bar N'^{n-2k_1}$
coincides by definition
with  $\tilde W^{n-k+k_1}$. The submanifold $\bar W_{2,-}^{n-k+k_1}
\subset \bar N'^{n-2k_1}$ coincides by
definition with $T(\tilde
W^{n-k+k_1})$.  Therefore $\bar N_2^{n-2k}$  coincides with $\bar N^{n-2k}$ and $N_2^{n-2k}$
coincides with $N^{n-2k}$. This proves the equivalence of Definition --1 and Definition --3 of
the homomorphism $J^{\D}$.

Recall, that the homomorphism
$$Imm^{sf}(n-k,k)
\stackrel{\delta_{k}}{\longrightarrow}  Imm^{\D}(n-2k,2k)$$
transforms the cobordism class of a triple $ (f, \kappa, \Xi) $ to
the cobordism class of the triple $ (g, \eta, \Psi) $, where $g:
N^{n-2k} \looparrowright \R^n$ is the immersion parameterizing the
self-intersection points manifold of the immersion $ f $ (it is
assumed that the immersion $ f $ intersects itself transversally),
$ \Psi $ is the $\D$--framing of the normal bundle of the
immersion $ g $, and $ \eta $ is the classifying map of the $ \D
$--framing $\Psi $.

We turn to the definition of the homomorphism
\begin{eqnarray}\label{1.66}
Imm^{\D}(n-2k,2k) \stackrel{h^{\D}_k}{\longrightarrow}  \Z/2,
\end{eqnarray}
 which will be called
the dihedral Hopf invariant. Define the subgroup
\begin{eqnarray}\label{1.6}
\I_c \subset \D,
\end{eqnarray}
generated by the transformations of the plane that preserve the
subspace spanned by each basis vector. The subgroup $ \I_c $ is an
elementary abelian 2-group of rank $2$. Define the homomorphism
\begin{eqnarray}\label{1.7}
l: \I_c \to \Z/2,
\end{eqnarray}
by sending an element $x$ of $\I_c$ to $0$ if $x$ fixes the first
basis vector, and to 1 if $x$ sends the first basis vector to its
negative.

 The subgroup $(\ref{1.6})$ has index
$2$ and the following 2-sheeted covering:
\begin{eqnarray}\label{1.8}
 K(\I_c,1) \to K(\D,1),
\end{eqnarray}
induced by this subgroup is well defined.

Denote by
\begin{eqnarray}\label{1.9}
\bar N^{n-2k} \to N^{n-2k}
\end{eqnarray}
the 2-sheeted covering induced by the classifying map $\eta:
N^{n-2k} \to K(\D,1)$
 from the covering $(\ref{1.8})$.
 The following characteristic class is well defined:
$$  \bar \eta_{sf} = l \circ \bar \eta_{\I_c} : \bar N^{n-2k} \to K(\Z/2,1), $$
where $\bar \eta_{\I_c}: \bar N^{n-2k} \to  K(\I_c,1)$ is the
double covering over the classifying map $ \eta: N ^{n-2k} \to K
(\D, 1) $ induced by the coverings $(\ref{1.8})$, $(\ref{1.9})$
over  the target and the source  of the map $ \eta $ respectively.

Let us define the homomorphism $Imm^{\D}(n-2k,2k)
\stackrel{h^{\D}_k}{\longrightarrow}  \Z/2$ by the formula:
\begin{eqnarray}\label{1.99}
h^{\D}_k([g,\eta,\Psi]) = \left\langle  (\bar
\eta_{sf})^{n-2k};[\bar N^{n-2k}] \right\rangle.
\end{eqnarray}

The diagram $(\ref{1.5})$ is now well defined.  Commutativity of
the right square of the diagram follows for $ k_1 = 1 $ by
Proposition $\ref{prop4}$, and for an arbitrary $ k_1 $ the proof
is similar. Let us prove the commutativity of the left square of
the diagram.
Lemma $ \ref{lemma7} $ is proved.
\[  \]

We need an equivalent definition of the dihedral Hopf invariant in
the case of $\D$--framed immersions in the codimension $2k$,
$n-4k>0$. Consider the subgroup of the orthogonal group $ O (4) $
that transforms the set of vectors $ (\pm \e_1, \pm \e_2, \pm
\e_3,\pm \e_4) $ of the standard basis into itself, perhaps by
changing the direction of some vectors and, moreover, preserving
the non-ordered pair of 2-dimensional subspaces $Lin(\e_1, \e_2)
$, $Lin(\e_3, \e_4)$ generated by basis vectors $ (\e_1, \e_2) $,
$ (\e_3, \e_4) $. Thus these 2-dimensional subspaces may be
preserved or interchanged. Denote the subgroup of these
transformations by $ \Z / 2 ^{[3]} $. This group has order $2^7$.
Define the chain of subgroups of index 2:
\begin{eqnarray}\label{1.10}
\I_c \times \D \subset \D \times \D \subset \Z/2^{[3]}.
\end{eqnarray}
The subgroup $ \D \times \D \subset \Z/2^{[3]}$ is defined as the subgroup of
transformations leaving invariant each 2-dimensional subspace
$Lin(\e_1, \e_2) $, $Lin(\e_3, \e_4) $
spanned by pairs of vectors $ (\e_1, \e_2) $, $ (\e_3, \e_4) $.
This subgroup is isomorphic to a direct product of two copies of $
\D $, each factor leaving invariant the corresponding
2-dimensional subspace. The subgroup $ \I_c \times \D \subset \D \times \D$
is defined as the subgroup of transformations that leave invariant
each linear subspace $Lin(\e_1) $, $Lin(\e_2)$ generated by vectors $ \e_1 $, $ \e_2 $.

Let the triple $ (g, \eta, \Psi) $ represent an element of $ Imm
^{\D} (n-2k, 2k) $, assuming that $ n-4k> 0 $ and that $ g $ is an
immersion in general position.
  Let $ L ^{n-4k} $ be the manifold of double self-intersection points of the
   immersion $g: N^{n-2k} \looparrowright \R^n$. The following
tower of 2-sheeted coverings
\begin{eqnarray}\label{1.12}
\bar L_{\I_c \times \D}^{n-4k} \to \bar L_{\D \times \D}^{n-4k} \to L^{n-4k}
\end{eqnarray}
is well defined by the following construction. (The covering $\bar L_{\D \times \D}^{n-4k} \to L^{n-4k}$ was considered above as the canonical covering over self-intersection manifold of the immersion $g$ and was denoted by
$\bar L^{n-4k} \to L^{n-4k}$.)
Let us consider the
parameterizing immersion $ h: L ^{n-4k} \looparrowright \R ^n $.
The normal bundle of the immersion $ h $ will be denoted by
$ \nu_L $. This bundle is classified by a  mapping
  $\zeta: L^{n-4k} \to K(\Z/2^{[3]},1)$.

The chain of subgroups $(\ref{1.10})$ induces a tower of 2-sheeted coverings of classifying spaces:
\begin{eqnarray}\label{1.11}
K(\I_c \times \D,1) \subset K(\D \times \D,1) \subset K(\Z/2^{[3]},1)
\end{eqnarray}
over the target of the classifying map $ \zeta $ and the tower of
2-sheeted coverings $(\ref{1.12})$ over the domain of the mapping
$ \zeta $. The covering $\bar L_{\I_c \times \D} \to  L^{n-4k}$, defined by
the formula $(\ref{1.12})$, will be called the canonical 4-sheeted
covering over the manifold of points of self-intersections of the
immersion $ g $.

Define the epimorphism
\begin{eqnarray}\label{1.51}
l^{[3]}: \I_c \times \D \to \I_d,
\end{eqnarray}
by sending an element $x$ of $\I_c \times \D$ to $0$
if $x$ fixes the first basis vector, and to 1 if $x$ sends the
first basis vector to its negative. This map induces the map of
the classifying spaces:
\begin{eqnarray}\label{1.41}
K(\I_c \times \D,1) \to K(\I_d,1).
\end{eqnarray}
 The classifying mapping $\bar
\zeta_{\I_c \times \D}: \bar L_{\I_c
\times \D}^{n-4k} \to K(\I_c \times \D,1)$ is well
defined as a result of the transition to a 4-sheeted covering over
the mapping $ \zeta $ and the classifying map $\bar \zeta_{sf}:
\bar L^{n-4k}_{\I_c \times \D} \to K(\I_d,1)$ is well
defined  as a result of the composition of the classifying map
$\bar \zeta_{\I_c \times \D}$ with the map $
(\ref{1.41}) $.

\begin{proposition}\label{prop8}
Suppose that a $\D$--framed immersion $ (g, \eta, \Psi) $
represents an element of $Imm^{\D}(n-2k,2k)$, the following
formula is satisfied:
\begin{eqnarray}\label{1.60}
\langle (\bar \zeta_{sf})^{n-4k}; [\bar L_{\I_c \times \D}^{n-4k}]
\rangle =\langle (\bar \eta_{sf})^{n-2k};[\bar N^{n-2k}]\rangle.
\end{eqnarray}
\end{proposition}

\subsubsection*{Proof of Proposition $ \ref{prop8} $}

Let $(g: N^{n-2k} \looparrowright \R^n,\Psi,\eta)$ be a $ \D
$--framed immersion, representing an element of
$Imm^{\D}(n-2k,2k)$ in the image of the homomorphism $\delta^k$.
Let $ L ^{n-4k} $ be the double-points manifold of the immersion $
g $, $ h: L ^{n-4k} \looparrowright \R^n $ be the parameterizing
immersion, and $ \bar L^{n-4k} \to L^{n-4k} $ be the canonical
double covering. Consider the image of the fundamental class $
\bar h_{ \ast} ([\bar L^{n-4k}]) \in H_{n-4k} (N^{n-2k}; \Z / 2) $
by means of the immersion $ \bar h: \bar L^{n-4k} \looparrowright
N^{n-2k} $ and let us denote by $ m \in H^{2k} (N^{n-2k}; \Z / 2)
$ the cohomology class that is Poincar\'{e}-dual to the
homology class $\bar h_{\ast}([\bar L^{n-4k}])$. Consider also the
cohomology Euler class of the normal bundle immersion $ g $, which
is denoted by $e \in H^{2k}(N^{n-2k};\Z/2)$.

By the Herbert Theorem (see [E-G], Theorem 1.1, coefficients $\Z/2$)
for immersion $g: N^{n-2k} \looparrowright \R^n$ with the
self-intersection manifold $L^{n-4k}$ the formula $e=m$ given in $
(\ref{1.500}) $ is valid.
Let us consider the classifying map
 $\eta: N^{n-2k} \to K(\D,1)$. Let us consider the $2$--sheeted cover
$K(\I_c,1) \to K(\D,1)$ over the classifying space. Let us induced
$\eta$ the 2-sheeted covering map over the map $\eta$, denoted by
$\bar \eta_{sf}: \bar N^{n-2k}_{sf} \to K(\I_c,1)$. (Note that in
the case $N^{n-2k}$ is a self-intersection points manifold of a
skew-framed immersion, the manifold
 $\bar N^{n-2k}_{sf}$ was considered above and this manifold was called
the canonical covering manifold over
 $N^{n-2k}$, this manifold was denoted by  $\bar N^{n-2k}$.)

Let us denote by $\bar e \in H^{2k}(\bar N^{n-2k}_{sf};\Z/2)$,
$\bar m \in H^{2k}(\bar N^{n-2k}_{sf};\Z/2)$ the images of the
cohomology classes $ e $, $ m $, respectively, under the canonical
double cover $\bar N^{n-2k}_{sf} \to N^{n-2k}$. The Herbert Theorem
implies that:
$$ \bar e= \bar m, $$
in particular, the following formula holds:
\begin{eqnarray}\label{1.61}
\langle (\bar \eta_{sf})^{n-4k} \bar m;[\bar N_{sf}^{n-2k}]
\rangle = \langle (\bar \eta^{sf})^{n-4k} \bar e;[\bar
N_{sf}^{n-2k}] \rangle.
\end{eqnarray}

 Because $\bar \eta_{sf}$ coincides with $\bar e$, the right side of the formula
 is equal to
$\langle (\bar \eta_{sf})^{n-2k};[\bar N^{n-2k}]\rangle $. Because
$\bar m$ is dual to the cohomology class  $\bar h_{\ast}[\bar
L^{n-4k}_{\H_{\bar c}}]$, the left side of the formula
$(\ref{1.61})$ can be rewritten in the form:  $\langle \bar
\eta_{sf})^{n-4k};[\bar L^{n-2k}_{H_{\bar c}}]\rangle$. Because
the classifying mappings  $\bar \zeta_{sf}: \bar L_{\H_{\bar
c}}^{n-4k} \to K(\Z/2,1)$ and $\bar \eta_{sf} \vert_{\bar L_{
\H_{\bar c}}}$ coincide, the left side of the formula
$(\ref{1.61})$ is equal to the characteristic number $\langle
(\bar \zeta_{sf})^{n-4k}; [\bar L_{\H_{\bar c}}^{n-4k}] \rangle$.
Proposition $\ref{prop8}$ is proved.
\[  \]

Let us generalize Proposition $\ref{prop5}$, Definition
$\ref{def6}$, and Lemma $\ref{lemma7}$.

 Let us assume that the triple  $(g,\eta,\Psi)$, where $g: N^{n-2k}
\looparrowright \R^n$ is an immersion, $\Psi$ is a $\D$--framing
of the normal bundle of the immersion $g$ with the characteristic
class  $\eta: N^{n-2k} \to K(\D,1)$, represents an element in
$Imm^{\D}(n-2k ,2k)$. Let $h: L^{n-4k} \looparrowright \R^n$ be an
immersion, that gives a parametrization of the self-intersection
manifold of the immersion $g$.

\begin{proposition}\label{prop5.1}
The normal $4k$--dimensional bundle  $\nu_{L}$ of the immersion
$h$ is isomorphic to the Whitney sum  of $k$ copies of a
4-dimensional bundle $\zeta$ over $L^{n-4k}$, each 4-dimensional
direct summand has the structure group  $\Z/2^{[3]}$ and is
classified by a classifying mapping  $\zeta: L^{n-4k} \to
K(\Z/2^{[3]},1)$.
\end{proposition}

\subsubsection*{Proof of Proposition $\ref{prop5.1}$}
The proof is omitted,  this proof is analogous to the proof of
Proposition $\ref{prop5}$.

\begin{definition}\label{def6.1}
We define the cobordism group of immersions
$Imm^{\Z/2^{[3]}}(n-4k,4k)$,  assuming $ n> 4k $. Let
$(h,\zeta,\Lambda)$  be a triple, which determines a
$\Z/2^{[3]}$--framed immersion of codimension $ 4k $. Here $h:
L^{n-4k} \looparrowright \R^n$ is an immersion and $\zeta:
L^{n-4k} \to K(\Z/2^{[3]},1)$  is the characteristic map of the
$\Z/2^{[3]}$--framing $ \Lambda $. The cobordism relation of
triples is standard.
\end{definition}

\begin{lemma}\label{lemma11}
Under the assumption $k_1 < k$, $4k < n$, the following
commutative diagram of groups is well defined:
\begin{eqnarray}\label{1.5.1}
\begin{array}{ccccc}
Imm^{\D}(n-2k_1,2k_1) & \stackrel{J^{\D}}{\longrightarrow} & Imm^{\D}(n-2k,2k)& \stackrel{h^{\D}_k}
{\longrightarrow} & \Z/2\\
\downarrow \delta^{\D}_{k_1}& &\downarrow \delta^{\D}_{k} && \| \\
Imm^{\Z/2^{[3]}}(n-4k_1,4k_1) &
\stackrel{J^{\Z/2^{[3]}}}{\longrightarrow} &
Imm^{\Z/2^{[3]}}(n-4k,4k) &
\stackrel{h^{\Z/2^{[3]}}_k}{\longrightarrow} & \Z/2.
\end{array}
\end{eqnarray}
\end{lemma}

\subsubsection*{Proof of Lemma $ \ref{lemma11} $}

Let us define the homomorphisms in the diagram   $(\ref{1.5.1})$. The
homomorphism  $h^{\D}_k$  is given by the characteristic number in the right
side of the formula   $(\ref{1.99})$. The homomorphism
$h^{\Z/2^{[3]}}_k$ is defined by means of the characteristic
number in the left side of the formula $(\ref{1.60})$. The
commutativity of the right square of the diagram is proved in
Proposition  $\ref{prop8}$.

We define the further homomorphisms in the diagram $ (\ref{1.5.1}) $. The
homomorphism
$$Imm^{\Z/2^{[3]}}(n-4k_1,4k_1)  \stackrel{J^{\Z/2^{[3]}}}{\longrightarrow}
Imm^{\Z/2^{[3]}}(n-4k,4k)$$
  is defined exactly  was the
homomorphism $J^{\D}$ in the bottom row of the diagram
$(\ref{1.5})$. Namely, let  a triple $(h, \zeta, \Lambda)$ 
represent an element in the cobordism group
$Imm^{\Z/3^{[3]}}(n-4k,4k)$, where $h: L^{n-4k} \looparrowright
\R^n$ is an immersion with $\Z/2^{[3]}$--framing $\Lambda$. Take
the universal $\Z/2^{[3]}$--bundle $\psi_{[3]}$ over
$K(\Z/2^{[3]},1)$  and take the pull-back of this bundle by means
of the classifying map $\zeta$, $\zeta^{\ast}(\psi_{[3]})$. Take a
submanifold $L'^{n-4k} \subset
  L^{n-4k_1}$  which represents the Euler class of
the bundle $(k-k_1)\eta^{\ast}(\psi_{[3]})$. The triple $(h',
\zeta', \Lambda')$  is well defined, where $h' = h \vert_{L'}$ пїЅ
$\zeta' = \zeta
  \vert_{L'}$. Let us define the $\Z/2^{[3]}$-framing $\Lambda'$.

Let us consider  the normal bundle $\nu_{h'}$ of the immersion
$h'$. This bundle decomposes into the Whitney sum of the two
bundles: $\nu_{h'}=\nu_h \vert_L  \oplus \nu_{L' \subset L}$,
where $\nu_{L' \subset L}$ is the normal bundle of the submanifold
$L'^{n-2k_1} \subset L^{n-2k}$. The bundles $\nu_{L' \subset L}$
and $(k-k_1)\zeta^{\ast}(\psi_{[3]})$ are isomorphic and this
bundle is equipped with the standard $\Z/2^{[3]}$--framing. The
bundle $\nu_h \vert_L$ is also equipped with the
$\Z/2^{[3]}$--framing. Therefore the bundle $\nu_{h}\vert_L$ is
equipped with the dihedral framing $\Lambda': \nu_{h'} \cong k_1
\zeta^{\ast}(\psi_{[3]})$. The triple $(h', \zeta', \Lambda')$
represent the element $J^{\Z/2^{[3]}}(h, \zeta, \Lambda) \in
Imm^{\Z/2^{[3]}}(n-4k_1,4k_1)$.

The homomorphism
$$Imm^{\D}(n-2k,2k)
\stackrel{\delta^{\D}_{k}}{\longrightarrow}
Imm^{\Z/2^{[3]}}(n-4k,4k)$$ transforms the cobordism class of a
triple $(g,\eta,\Psi)$ to the cobordism class of the triple
$(h,\zeta,\Lambda)$,  where $h: L^{n-4k} \looparrowright \R^n$ is
the immersion parameterizing the self-intersection points manifold
of the immersion $ g $ (it is assumed that the immersion $ g $
intersects itself transversely), $ \Lambda $ is the
$\Z/2^{[3]}$--framing of the normal bundle of the immersion $ h $,
and $ \zeta $ is the classifying mapping of the $
\Z/2^{[3]}$--framing $\Lambda $.

The commutativity of the left square in the diagram
$(\ref{1.5.1})$ is proved analogously with the commutativity of
the left square in the diagram $(\ref{1.5})$. Lemma
$\ref{lemma11}$ is proved.

\section{Proof of the main theorem}

We reformulate the Main Theorem ($\sigma(n)$ is defined by the formula $(\ref{sigma})$),
taking into account the notation of the previous section.

\begin{theorem}\label{th9}
For $\ell \ge 8$ the homomorphism $h_{2^{\ell-2}-n_{\sigma}}^{sf}:
Imm^{sf}(3 \cdot
2^{\ell-2}+n_{\sigma-4},2^{\ell-2}-1-n_{\sigma-4}) \to \Z/2$,
given by the equivalent Definitions $\ref{def1}$ and $\ref{proba}$
is trivial. (Recall, by the formula  $(\ref{sigma})$, $n_{\sigma-2}=2^{\sigma-2}-1 \ge 1$).
\end{theorem}

\begin{remark}\label{127}
In the cases $\ell=4$ the proof by means of the
considered approach is unknown, because of the dimensional
restriction in Lemma $\ref{osnlemma1}$, in the
case $\ell=5,6,7$ because of the dimensional restrictions in Lemma
$\ref{osnlemma2}$. For a possible approach in the cases $n=31$,
$n=63$ and $n=127$ see a remark in a remark at the beginning of Section 5.
\end{remark}

\subsubsection*{Proof of Main Corollary from Theorem $\ref{th9}$}

Consider the homomorphism $J^{sf}: Imm^{sf}(n-1,1) \to Imm^{sf}(3
\cdot 2^{\ell-2}+n_{\sigma-4},2^{\ell-2}-1-n_{\sigma-4})$.
According to Proposition $\ref{prop2}$,
$h_1^{sf}=h_{2^{\ell-2}-n_{\sigma}-1}^{sf} \circ J^{sf}$.  Let an
element of the group $Imm^{sf}(n-1,1)$ be represented by an
immersion $f: M^{n-1} \looparrowright \R^n$, $w_1(M)=\kappa$. The
value $h_k^{sf}(J^{sf}(f))$, where  $k=2^{\ell-2}-n_{\sigma}-1$,
coincides with the characteristic number $\langle
\kappa^{n-1};[M^{n-1}]\rangle$. Applying Theorem $\ref{th9}$, we
conclude the proof of the Main Theorem.
\[  \]

For the proof of Theorem $\ref{th9}$ we shall need the fundamental
Definitions $\ref{defcycl}$, $\ref{defQ}$,  $\ref{def16}$, whose formulation will
require some preparation.

\subsubsection*{Definition of the subgroups $\I_d \subset \I_a \subset \D$, $\I_{b \times \bb}$}

We denote by $\I_a \subset
 \D$ the cyclic subgroup of order $4$ and index $2$,
containing the nontrivial elements $a,a^2,a^3 \in \D$ (i.e.,
generated by the plane rotation which exchanges the coordinate
axes). We denote by $\I_d \subset \I_a$ the subgroup of index 2
with nontrivial elements $a^2$. We denote by $\I_b \subset \D$ an elementary
subgroup with the generator $b$, by $\II_b \subset \D$ an elementary
subgroup with the generator $ab^2$ 2 
(i.e., $\I_b$ and $\II_b$ are generated by the reflections with respect to the bisectors
of the coordinate axes). The subgroup $\I_b \times \II_b \subset \D$ denote by  $\I_{b \times \bb} \subset \D$.

The following inclusion homomorphisms of subgroups are well
defined: $i_{d,a}: \I_d \subset \I_a$,  $i_{d,b}: \I_d \subset
\I_b$. When the image coincides with the entire group $ \D $ the
corresponding index for the inclusion homomorphism will be
omitted: $i_d: \I_d \subset \D$, $i_a: \I_a \subset \D$, $i_b:
\I_b \subset \D$.
\[  \]

\subsubsection*{Definition of the subgroup $i_{\Q}: \Q \subset \Z/2^{[3]}$}

Let $ \Q $ be the quaternion group of order $8$. This group has presentation $\{ \i,\j,\k \quad \vert \quad
\i\j = \k = -\j\i,
\j\k=\i = -\k\j, \k\i = \j = -\i\k, \i^2 =\j^2 = \k^2 =-1 \}$.
 There is a
standard representation  $ \chi: \Q \to O(4) $. The
representation $ \chi $ (a matrix acts to the left on a vector)
carries the unit quaternions $ \i, \j, \k $ to the matrices
\begin{eqnarray}\label{Q a1}
\left(
\begin{array}{cccc}
0 & -1 & 0 & 0 \\
1 & 0 & 0 & 0 \\
0 & 0 & 0 & -1 \\
0 & 0 & 1 & 0 \\
\end{array}
\right),
\end{eqnarray}
\begin{eqnarray}\label{Q a2}
\left(
\begin{array}{cccc}
0 & 0 & -1 & 0 \\
0 & 0 & 0 & 1 \\
1 & 0 & 0 & 0 \\
0 & -1 & 0 & 0 \\
\end{array}
\right),
\end{eqnarray}
\begin{eqnarray}\label{Q a3}
\left(
\begin{array}{cccc}
0 & 0 & 0 & -1 \\
0 & 0 & -1 & 0 \\
0 & 1 & 0 & 0 \\
1 & 0 & 0 & 0 \\
\end{array}
\right).
\end{eqnarray}
These matrices give the action of $\bf{left}$ multiplication by
$\i$, $\j$, and $\k$ on the standard basis for the quaternions.
This representation $\chi$ defines the subgroup
 $i_{\Q}: \Q \subset
\Z/2^{[3]} \subset O(4)$.

\subsubsection*{Definition of the subgroups $\I_d \subset \I_a \subset \Q$}

Denote by $i_{\I_d,\Q}: \I_d \subset \Q$
 the central subgroup of the quaternion group,
which is also the center of the whole group $\Z/2^{[3]}$.

Denote by $i_{\I_a,\Q}: \I_a \subset \Q$
 the  subgroup of the quaternion group
generated by the quaternion $ \i $.

The following inclusions are well defined:
$i_{\I_d}: \I_d \subset \Z/2^{[3]}$,
$i_{\Q}: \Q \subset \Z/2^{[3]}$.
\[  \]

\begin{definition}\label{defcycl}
We say that a classifying map $\eta: N^{n-2k} \to K(\D,1)$ is
cyclic if it can be factored as a composition of a map $\mu_a:
N^{n-2k} \to K(\I_a,1)$ and the inclusion $i_a: K(\I_a,1) \subset
K(\D,1)$. We say that the mapping 
$\mu_a$ determines a reduction of the classifying mapping  $\eta$.
\end{definition}

\begin{definition}\label{defQ}
We say that a classifying map $\zeta: L^{n-4k} \to K(
\Z/2^{[3]},1)$ is quaternionic if it can be factored as a
composition of a map $\lambda:  L^{n-4k} \to K(\Q,1)$ and the
inclusion $i_{\Q}: K(\Q,1) \subset K(\Z/2^{[3]},1)$. We also say that the mapping 
$\lambda$ determines a reduction of the classifying mapping  $\zeta$.
\end{definition}

We shall later require the construction of the Eilenberg-Mac Lane
spaces $K(\I_a,1)$, $K(\Q,1)$ and a description of the finite
dimensional skeletпїЅ of these spaces, which we now recall.

Consider the infinite dimensional sphere $S^{\infty}$ (a contractible
space), which it is convenient to define as a direct limit of an infinite
sequence of inclusions of standard spheres of odd dimension,
$$S^{\infty}=
\lim_{\longrightarrow} \quad (S^1 \subset S^3 \subset \cdots \subset
S^{2j-1} \subset S^{2j+1} \subset \cdots  \quad) .$$

Here $S^{2n-1}$ is defined by the formula $S^{2j-1}=\{ \- (z_1,
\dots, z_j) \in \C^j,  \quad \vert z_1 \vert ^2+ \dots + \vert
z_j\vert ^2=1 \- \}$. Let $\i(z_1, \dots, z_j)=(\i z_1, \dots, \i
z_j)$.

Then the space $S^{2j-1}/\i$, which is called the
$(2j-1)$-dimensional lens space over $\Z/4$, is the
$(2j-1)$--dimensional skeleton of the space $K(\I_a,1)$. The space
$S^{\infty}/\i$ itself is the Eilenberg-Mac Lane space
$K(\I_a,1)$. The cohomology ring of this space is well-known, see
e.g. [A-M].

Let us define the Eilenberg-Mac Lane space $K(\Q,1)$.
Consider the infinite dimensional sphere $S^{\infty}$, which now it is convenient to define as a direct limit of an infinite
sequence of inclusions of standard spheres of dimensions $4j+3$:
$$S^{\infty}=
\lim_{\longrightarrow} \quad (S^3 \subset S^7 \subset \cdots \subset
S^{4j-1} \subset S^{4j+3} \subset \cdots  \quad).$$

 A coordinate action $\Q
\times (\C^{2})^j \to (\C^{2})^j$ is defined on each direct
summand $\mathbb{H}=\C^2$ in accordance with the formulas $(\ref{Q
a1}), (\ref{Q a2}), (\ref{Q a3})$.
  Thus,
the space $S^{4j-1}/\Q$ is a $ (4j-1) $-dimensional skeleton of
the space $S^{\infty}/\Q$ and this space  is called the $ (4j-1)
$-dimensional lens space over $ \Q $. The space $S^{\infty}/\Q$
itself is the Eilenberg-Mac Lane space $K(\Q,1)$. The cohomology
ring of this space is well known, see [At] section 13.

Denote by
 $\psi_+: E(\psi_+) \to S^{4k-1}/\i$ the restriction of the universal $SO(2)$--bundle  over  $K(\Z/4,1)$ on the standard skeleton
 of the dimension  $4k-1$. Denote by 
$\psi_-: E(\psi_-) \to S^{4k-1}/\i$ the $SO(2)$--bundle, which is the bundle  $\psi_+$ with the opposite orientation of fibers. 
Denote by
 $T_{\Q}: S^{4k-1}/\i \to S^{4k-1}/\i$ the free involution, which is constructed by means of the normal subgroup $\I_a \subset \Q$ 
 of the index 2.  

\begin{lemma}\label{TQ}
The following isomorphism of 
$SO(2)$--bundles is well defined:  
\begin{eqnarray}\label{TQQ}
\psi_- \cong T_{\Q}^{\ast}(\psi_+).
\end{eqnarray}
\end{lemma}

\subsubsection*{Proof of Lemma $\ref{TQ}$}
Consider the trivial 2-bundle 
 $\C \times S^{4k-1} \to S^{4k-1}$. Define the action 
$\I_a \times \C \times S^{4k-1} \to \C \times S^{4k-1}$ by the formula: 
$(\i,x,y) \mapsto (\i x,\i y)$, $\i \in \I_a$, $x \in \C$, $y \in S^{4k-1}$. 
The quotient of the considered action is the total space of the bundle $\psi_+$, which will be denoted by $E(\psi_+)$. 
The following projection  $E(\psi_+) \to S^{4k-1}/\i$ is well-defined, this is the required 
$SO(2)$--bundle  $\psi_+$. Let us define the bundle $\psi_-$ by the analogous construction, using the action  $(\i,x,y) \mapsto (-\i x,\i y)$.
The bundles $\psi_+$ and $\psi_-$ are isomorphic as  $O(2)$--bundles, the isomorphism is given by the formula:
$(x,y)/\i \mapsto (R(x),y)/\i$, where $R(x)$ is the reflection with respect to the line, which is determined by the first base vector of the real fiber. 
 
Let us prove that the bundles  
$T_{\Q}(\psi_+)$ and $\psi_-$ are isomorphic. Define the automorphism of the trivial bundle
$\C \times S^{4k-1}$ (which is not the identity on the base) by the formula:
$(R(x),y) \mapsto (x,T_{\Q}(y))$. Because $T_{\Q}(\i y))=-\i T_{\Q}(y)$ and 
$R(\i (x))=-\i R(x)$, we get
 $(R(\i x),\i y) \mapsto (-R(\i (x)),-T_{\Q}(\i y))$. Therefore, the quotient $\sim \I_a$ is well defined and the required isomorphism
 $(\ref{TQQ})$ is well defined. Lemma  $\ref{TQ}$ is proved.

\subsubsection*{Definition of the characteristic number $h_{\mu_a,k}$}

Let us assume  $ n> 4k $ and let us assume that on the manifold
$N^{n-2k}$ of self-intersection points of  a skew-framed immersion
there is defined a map $\mu_a: N^{n-2k} \to K(\I_a,1)$. Define the
characteristic value $h_{\mu_a,k}$  by the formula:
\begin{eqnarray}\label{mu_a}
 h_{\mu_a,k}=\langle \bar e_g \bar \mu_a^{\ast}x;[\bar
N_a^{n-2k}]\rangle,
\end{eqnarray}
where  $\bar \mu_a: \bar N_a^{n-2k} \to K(\I_d,1)$  is a double
cover over the map $\mu_a: N_a^{n-2k} \to K(\I_a,1)$, induced by
the cover $K(\I_d,1) \to K(\I_a,1)$, $x \in
H^{n-4k}(K(\I_d,1);\Z/2)$ is the generator, $\bar e_g \in H^{k}(\bar N_a^{n-2k};\Z/2)$
is the image of the Euler class  $e_g \in H^k(N^{n-2k};\Z/2)$ of the immersion $g$ by means of the covering
$p_a: \bar N_a^{n-2k} \to N^{n-2k}$,
$\bar e_g = p_a^{\ast}(e_g)$, and $[\bar
N_a^{n-2k}]$ is the fundamental class of the manifold $\bar
N_a^{n-2k}$. (The manifold $\bar N_a^{n-2k}$ coincides with
the canonical 2-sheeted covering $\bar N^{n-2k}$, if and only if the
classifying mapping $\eta: N^{n-2k} \to K(\D,1)$ is cyclic, see
Definition $\ref{defcycl}$.)

Let us assume that
 $k \equiv 0\pmod{2}$. Then the characteristic number $(\ref{mu_a})$ is the reduction modulo 2 of
 the following characteristic number, denoted the same, determined modulo 4:
\begin{eqnarray}\label{2909}
\langle e_g \mu_a^{\ast}x;[N^{n-2k}]\rangle,
\end{eqnarray}
where $x \in H^{n-4k}(K(\I_a,1);\Z/4)$ is the generator,  $e_g \in
H^{k}(N^{n-2k};\Z/4)$ is the Euler class of the co-oriented
immersion $g$ with coefficients modulo 4, $[N^{n-2k}]$ is the
fundamental class of the oriented manifold  $N^{n-2k}$ with
coefficients modulo 4.

\subsubsection*{Definition of the characteristic number $h_{\lambda,k}$}

Let us assume  $ n> 4k $ and let us assume that on the manifold
$L^{n-4k}$ of self-intersection points of  a $\D$-framed immersion
there is defined a map $\lambda: L^{n-4k} \to K(\Q,1)$. Define the
characteristic value $h_{\lambda}^{k}$  by the formula:
\begin{eqnarray} \label{lambda_a}
 h_{\lambda,k}=\langle \bar
\lambda^{\ast}y;[\bar L_{\I_d}]\rangle,
\end{eqnarray}
where  $y \in H^{n-4k}(K(\I_d,1);\Z/2)$ is a generator,
\begin{eqnarray}\label{302}
\bar \lambda_{\I_d}: \bar L_{\I_d}^{n-4k} \to K(\I_d,1)
\end{eqnarray}
 is a 4-sheeted cover over
the map  $\lambda: L^{n-4k} \to K(\Q,1)$, induced by the cover
$K(\I_d,1) \to K(\Q,1)$, and $[\bar L_{\I_d}]$  is  the
fundamental class of the manifold $\bar L_{\I_d}^{n-4k}$. (The
manifold $\bar L_{\I_d}^{n-4k}$ coincides with the canonical
4-sheeted covering $\bar L^{n-4k}_{\I_c \times \D}$,
if and only if the classifying mapping $\zeta: L^{n-4k} \to
K(\Z/2^{[3]},1)$ is quaternionic, see Definition $\ref{defQ}$.)

Let us assume that
$k \equiv 0\pmod{2}$. Then the characteristic number $(\ref{lambda_a})$ is the reduction
modulo 2 of the following characteristic number, denoted the same, determined modulo 4:
\begin{eqnarray}\label{3024}
\langle \bar
\lambda^{\ast}y;[\bar L]\rangle,
\end{eqnarray}
where $y \in
H^{n-4k}(K(\I_a,1);\Z/4)$ is the generator, $[\bar L]$
is the fundamental class of
the oriented manifold  $\bar L^{n-4k}$ (the manifold   $\bar L^{n-4k}$ is the canonical 2-sheeted covering over the manifold $L^{n-4k}$) with coefficients modulo 4.

\begin{lemma}\label{lemma12}
For an arbitrary skew-framed immersion $(f: M^{n-k}
\looparrowright \R^n, \kappa, \Xi)$ with
self-intersection manifold $N^{n-2k}$ for which the classifying
mapping $ \eta $ of the normal bundle is cyclic, the following
equality is satisfied:
$$h_k^{sf}(f,\kappa, \Xi) = h_{\mu_a,k},$$
where the characteristic value on the right side is calculated for
a mapping $\mu_a$, satisfying the condition $\eta = i_a \circ
\mu_a$,  $i_a: K(\I_a,1) \subset K(\D,1)$.
\end{lemma}

\subsubsection*{Proof of Lemma $\ref{lemma12}$}

Consider the double cover $\bar \mu_a: \bar N_a^{n-2k} \to
K(\I_d,1)$ over the mapping $\mu_a: N^{n-2k} \to K(\I_a,1)$,
induced by the double cover $K(\I_d,1) \to K(\I_a,1)$ over the
target space of the map. Since the structure mapping $ \eta$ is
cyclic, the manifold $\bar N^{n-2k}_a$ coincides with the
canonical 2-sheeted cover $\bar N^{n-2k}$ over the
self-intersections manifold $ N^{n-2k} $ of the immersion $ f $, the class
 $\bar e_g \in H^{2k}(N^{n-2k};\Z/2)$ coincides with the class  $\bar \eta_{sf}^{2k}$,
$\bar \eta_{sf} \in H^1(\bar N^{n-2k};\Z/2)$.
The proof of the lemma follows from Lemma $\ref{lemma7}$ since the
mappings $\bar \mu_a$ and $\bar \eta$  coincide and the
characteristic number $h_{\mu_a,k}$ is computed as in the right
side of the equation ($ \ref{1.99} $).
\[  \]

\begin{lemma}\label{lemma13}
For an arbitrary $\D$-framed immersion $(g: N^{n-2k}
\looparrowright \R^n, \eta, \Psi)$ with a self-intersection
manifold $L^{n-4k}$, for which the classifying mapping $ \zeta $
of the normal bundle is quaternionic, the following equality is
satisfied:
$$h_k^{sf}(g,\eta, \Psi) = h_{\lambda,k},$$
where the characteristic value on the right side is calculated by the formula $(\ref{lambda_a})$ for
a mapping $\lambda$, satisfying the condition $\zeta = i_{a} \circ
\lambda$, $i_{a}: K(\Q,1) \subset K(\Z/2^{[3]},1)$.
\end{lemma}

\subsubsection*{Proof of Lemma $\ref{lemma13}$}

Define the 2-sheeted covering
\begin{eqnarray}\label{301}
\bar \lambda_{\I_a}: \bar L_{\I_a}^{n-4k} \to K(\I_a,1)
\end{eqnarray}
 over the mapping $\lambda: L^{n-4k} \to K(\Q,1)$,
induced by the 2-sheeted cover $K(\I_a,1) \to K(\Q,1)$ over the
target space of the map. Let us consider the 4-sheeted covering
$\bar \lambda_{\I_d}: \bar L_{\I_d}^{n-4k} \to K(\I_d,1)$
over the mapping $\lambda: L^{n-4k} \to K(\Q,1)$, defined by the formula $(\ref{302})$,  over the
target space of the map.

Since the structure mapping $ \zeta$ is quaternionic, the manifold
$\bar L^{n-4k}_{\I_d}$ coincides with the canonical 4-sheeted covering
$\bar L^{n-4k}$ over the self-intersections manifold $L^{n-4k}$ of
the immersion $ g $. The proof of the lemma follows from
Proposition $\ref{prop8}$, since the mappings $\bar \lambda$ and
$\bar \zeta$ coincide and the characteristic number
$h_{\lambda}^{k}$ is computed as in the left side of the equation
($\ref{1.60}$).
\[  \]

  



Let us define the subgroup $\H_b \subset \Z/2^{[3]}$ as a product
of the subgroup $\I_a \subset \Q \subset \Z/2^{[3]}$ and an
elementary subgroup, with the only non-trivial element $t$ given
by the transformation transposing each pair of the corresponding
basis vectors $\e_1 ={\bf 1}$ and $\e_3 =\j$ and
 the pair of the basis vectors $ \e_2 = \i $ and $ \e_4 = \k $, preserving their direction.
It is easy to verify that the group $ \H_b $ has the order $8$ and
this group is isomorphic to $ \Z / 4 \times \Z / 2 $.
The groups
$\H_b$ and $\Q$ contains the common index 2 subgroup: $\I_a \subset \H_b$,
$\I_a \subset \Q$.

\begin{definition}\label{def16}
Let $(g,\Psi,\eta)$ be a $\D$-framed immersion, where the
immersion $g: N^{n-2k} \looparrowright \R^n$ is assumed to be in a
general position with  self-intersection manifold denoted by
$L^{n-4k}$. Assume that the manifold $N^{n-2k}$ contains a marked component
$N^{n-2k}_a \subset N^{n-2k}$, with self-intersection manifold $L^{n-4k}_a \subset L^{n-4k}$.

Let the component
$N^{n-2k}_a$ be equipped with a mapping
$\mu_a: N^{n-2k}_a \to K(\I_a,1)$, which is determined a reduction of the restriction
of the classifying mapping 
$\eta$ to the component  $N^{n-2k}_a$ (see property  1 in Definition $(\ref{def14})$).

 Assume that the manifold  $L^{n-4k}_a$ is the disjoint union of the two closed submanifolds:
 $L^{n-4k}_a = L^{n-4k}_{\Q} \cup L^{n-4k}_{\H_b}$.
 Moreover, there exists a pair of
  mappings
  $(\mu_a,
 \lambda)$, where
 $\mu_a: N^{n-2k}_a \to K(\I_a,1)$, $\lambda = \lambda_{\Q} \cup \lambda_{\H_b}: L^{n-4k}_{\Q}
\cup  L^{n-4k}_{\H_b} \to K(\Q,1) \cup
K(\H_b,1)$.
Define the manifold
 $\bar L_{\Q}^{n-4k} \cup \bar L_{\H_b}^{n-4k}$ and its mapping
 \begin{eqnarray}\label{barA}
\bar \lambda = \bar \lambda_{\Q} \cup \bar \lambda_{\H_b}: \bar L_{\Q}^{n-4k} \cup
\bar L_{\H_b}^{n-4k} \to  K(\I_a,1) \cup K(\I_a,1),
\end{eqnarray}
as the 2-sheeted covering mapping over the disjoint union of the
mappings $\lambda_{\Q} : L^{n-4k}_{\Q} \to K(\Q,1)$,
$\lambda_{\H_b} : L^{n-4k}_{\H_b} \to K(\H_b,1)$ which is induced
from 2-sheeted coverings $K(\I_a,1) \to K(\Q,1)$, $K(\I_a,1) \to
K(\H_b,1)$ over the target space of the mapping $\lambda$. We say
that this $\D$-framed immersion $(g, \Psi, \eta)$
  admits a quaternionic structure if the following two conditions are
satisfied:

--1. The manifold $\bar L_{\Q}^{n-4k} \cup \bar L_{\H_b}^{n-4k}$ is diffeomorphic to the
canonical 2-sheeted covering manifold $\bar L^{n-4k}$ over the
self-intersection manifold of the immersion $g$, and the mapping
\begin{eqnarray}\label{barAA}
(Id \cup Id) \circ \bar \lambda : \bar L_{\Q}^{n-4k} \cup
\bar L_{\H_b}^{n-4k} \to  K(\I_a,1) \cup K(\I_a,1) \to K(\I_a,1) \to K(\I_a,1),
\end{eqnarray}
where $Id \cup Id: K(\I_a,1) \cup K(\I_a,1) \to K(\I_a,1)$  is the identity map
over each component, coincides with the restriction
of the mapping $\mu_a$ to the submanifold $\bar
L^{n-4k} \looparrowright N^{n-2k}$.

--2. The following equation is satisfied: 
\begin{eqnarray}\label{uravn}
h_{\mu_a,k}=h_k^{\D}(g,\eta,\Xi),
\end{eqnarray}
where the characteristic number in the left side of the formula is given by $(\ref{mu_a})$ and on the right side is given by the formula $(\ref{1.99})$. 
\end{definition}

\begin{example}\label{def17}
Let us assume that the classifying map $ \zeta $ is quaternionic.
The quaternionic structure is defined by the mappings  $\mu_a:
N^{n-2k} \to K(\I_a,1)$, $\lambda: L^{n-4k} \to K(\Q,1) \subset
K(\Z/2^{[3]},1)$, where $i_{\I_a} \circ \mu_a = \eta$,
$i_{\I_a}: K(\I_a,1) \subset K(\D,1)$, $i_{\Q} \circ \lambda =
\zeta$, $i_{\Q}: K(\Q,1) \subset K(\Z/2^{[3]},1)$.
\end{example}

\begin{lemma}\label{lemma18}
Assume that a pair $(\mu_a, \lambda)$ determines a quaternionic structure
for a
$\D$--framed immersion $(g,\eta,\Psi)$. Then the characteristic number
$h_{\lambda,k}$, determined by the formula
$(\ref{lambda_a})$, coincides with the characteristic number
$h_{\mu_a}^{k}$, given by the formula $(\ref{mu_a})$.
\end{lemma}

\subsubsection*{Proof of Lemma $\ref{lemma18}$}
The proof is analogous to the proof of Proposition $\ref{prop8}$.
\[  \]

\begin{corollary}\label{cor18}
Assume that a pair $(\mu_a, \lambda)$ determines a quaternionic structure
for a
$\D$--framed immersion $(g,\eta,\Psi)$. Then
\begin{eqnarray}\label{hQH}
 h_{\mu_a,k}=\langle \bar e_g \bar \mu_a^{\ast}x;[\bar
N_a^{n-2k}]\rangle = h_{\lambda,k}(L_{\Q}) + h_{\lambda,k}(L_{\H_b}),
\end{eqnarray}
where the terms in the right side are defined by the formula
$(\ref{lambda_a})$ for each corresponding component of the
manifold $L^{n-4k}$.
\end{corollary}

\section{Cyclic structure for formal (equivariant) mappings}


\begin{theorem}\label{prop22}
For  $k \ge 5$, an arbitrary element of the group 
$Imm^{\D}(n-2k,2k)$, which is in the image of the homomorphism  $(\ref{JD})$
$J^{\D}: Imm^{sf}(n-k,k) \to Imm^{\D}(n-2k,2k)$, 
is represented by a $\D$--framed immersion 
$(g,\eta, \Psi)$, where $g: N^{n-2k} \looparrowright \R^n$, 
$N^{n-2k} = N^{n-2k}_{a} \cup N^{n-2k}_{\D}$,
the restriction of the classifying mapping   $\eta \vert_{N^{n-2k}_a}$ to the component   $N^{n-2k}_a \subset N^{n-2k}$ 
is cyclic in the sense of Definition  $\ref{defcycl}$, where the characteristic number $(\ref{1.99})$, which is calculated 
for the mapping 
$\eta \vert_{N^{n-2k}_{\D}}$, is trivial. 
\end{theorem}
\[  \]

We need to reformulate the notion of a cyclic structure and of a quaternionic structure without the assumption
that the corresponding maps $f: M^{n-k} \to
\R^n$ and $g: N^{n-2k} \to \R^n$ are  immersions. We formulate
the necessary definition in minimal generality, under the assumption that
$M^{n-k}=\RP^{n-k}$, $N^{n-2k} = S^{n-2k}/\i$.

Let 
\begin{eqnarray}\label{d}
d: \RP^{n-k} \to \R^n
\end{eqnarray}
  be an arbitrary $PL$-mapping.
Consider the two-point configuration space
\begin{eqnarray}\label{99}
(\RP^{n-k} \times
\RP^{n-k} \setminus \Delta_{\RP^{n-k}})/T',
\end{eqnarray}
which is called the пїЅdeleted squareпїЅ of the space $\RP^{n-k}$. This space is obtained
as the quotient of the direct product without the diagonal by the
involution $T': \RP^{n-k} \times \RP^{n-k} \to
\RP^{n-k} \times \RP^{n-k}$, exchanging the coordinates.
This space is an open manifold. It is convenient to define an analogous
space, which is a manifold with boundary.

Define the space $\bar \Gamma$ as a spherical blow-up of the
space $\RP^{n-k} \times \RP^{n-k} \setminus \Sigma_{diag}$ in the
neighborhood of the diagonal. The spherical blow-up is a manifold
with boundary, which is defined as a result of compactification of
the open manifold $\RP^{n-k} \times \RP^{n-k} \setminus
\Sigma_{diag}$ by the fiberwise glue-in of the fibers of the unit
sphere bundle $ST \Sigma_{diag}$ of the tangent bundle $T
\Sigma_{diag}$ in the neighborhood of zero--sections of the normal
bundle of the diagonal $\Sigma_{diag} \subset \RP^{n-k} \times
\RP^{n-k}$. The following natural inclusions are well defined:
$$\RP^{n-k} \times \RP^{n-k} \setminus \Sigma_{diag} \subset \bar \Gamma,$$
$$ST \Sigma_{diag} \subset \bar \Gamma.$$
On the space  $\bar \Gamma$ the free involution $\bar T': \bar
\Gamma_0 \to \bar \Gamma$, which is  an extension of the
involution $ T' $ is well  defined.

The quotient
$\bar \Gamma / \bar T'$
is denoted by $ \Gamma $, and the corresponding double covering by
$$
p_{\Gamma}:  \bar \Gamma / \bar T' \to \Gamma.
$$
The space $ \Gamma $ is a manifold with boundary and it is
called the resolution space of the configuration space $
(\ref{99}) $. The projection $p_{\partial \Gamma}: \partial
\Gamma \to \RP^{n-k}$ is well defined, and is called a
resolution of the diagonal.


 


\subsubsection*{Formal (equivariant) mapping with holonomic self-intersection}
Denote by 
 $T_{\RP^{n-k}}$, $T_{\R^n}$ the standard involutions on the spaces  $\RP^{n-k} \times \RP^{n-k}$, $\R^n \times \R^n$, which permutes the coordinates.
 Let 
\begin{eqnarray}\label{d2}
d^{(2)}: \RP^{n-k} \times \RP^{n-k} \to \R^n \times \R^n 
\end{eqnarray}
be an arbitrary
 $(T_{\RP^{n-k}}$, $T_{\R^n})$--equivariant mapping, which is transversal along the diagonal of the source space.
Denote  $(d^{(2)})^{-1}(\R^n_{diag})/T_{\RP^{n-k}}$ by $\N=\N(d^{(2)})$, let us call this polyhedron a self-intersection (formal) polyhedron of 
the mapping  $d^{(2)}$. In the case the formal mapping $d^{(2)}$ is the extension of a mapping $(\ref{d})$,
the polyhedron $\N(d^{(2)})$ coincides with the polyhedron, denoted by the formula:
\begin{eqnarray}\label{Nd}
\N(d) = Cl \{ ([x,y]) \in \Gamma_{\circ} : y \ne x, d(y)=d(x) \}.
\end{eqnarray}
Note that 
$int(\Gamma)$ is homeomorphic to   $\Gamma_{\circ}$. Denote $int(\N(d^{(2)}) = \N(d^{(2)}) \setminus 
 (\N(d^{(2)}) \cap \Delta_{\Gamma})$ by $\N(d^{(2)})_{\circ}$.

 There is a
canonical double covering
\begin{eqnarray}\label{pNd}
p_{\N}: \bar \N \to \N,
\end{eqnarray}
ramified over the boundary $\partial \N$ (above this
boundary the cover is a diffeomorphism). The following diagram is
commutative:
$$
\begin{array}{ccc}
i_{\bar \N}: (\bar \N, \partial \N) & \subset & (\bar \Gamma, \partial \Gamma) \\
\\
\downarrow p_{\N}&   &  \downarrow  p_{\Gamma} \\
\\
i_{\N}: (\N, \partial \N) & \subset & (\Gamma, \partial \Gamma).
\end{array}
$$

\subsubsection*{Structural map  $\eta_{\N\circ}: \N(d)_{\circ} \to
K(\D,1)$}
Define the mapping  
\begin{eqnarray}\label{stG}
\eta_{\Gamma}: \Gamma \to K(\D,1),
\end{eqnarray}
which we shall call the structure mapping of the пїЅdeleted squareпїЅ. Note that the
inclusion $\bar \Gamma \subset \RP^{n-k} \times \RP^{n-k}$
induces an isomorphism of fundamental groups, since the
codimension of the diagonal $\Delta_{\RP^{n-k}} \subset \RP^{n-k}
\times \RP^{n-k}$ is equal to $n-k$ and satisfies the inequality $ n-k \ge 3 $.
Therefore, the equality is satisfied:
\begin{eqnarray}\label{pi_1}
\pi_1(\bar \Gamma) = H_1(\bar \Gamma;\Z/2)= \Z/2 \oplus \Z/2.
\end{eqnarray}

Consider the induced automorphism $T'_{ast}: H_1(\bar
\Gamma;\Z/2) \to H_1(\bar \Gamma;\Z/2)$.  Note that this
automorphism is not the identity. Fix an isomorphism of the groups
$H_1(\bar \Gamma;\Z/2)$ and  $\I_c$, which maps the generator of
the first (respectively second) factor  of $ H_1 (\bar \Gamma; \
Z / 2) $, see (\ref{pi_1}), into the generator $ ab \in \I_c
\subset \D $ (respectively, $ ba \in \I_c \subset \D $), which in
the standard representation of the group $ \D $ is defined by the
reflection with respect to the second (respectively, the first)
coordinate axis.

It is easy to verify that the  automorphism of the conjugation
with respect to the subgroup $ \I_c \subset \D $ by means of the
element $ b \in \D \setminus \I_c $ (in this formula the element $
b $ can be chosen arbitrarily), defined by the formula $ x \mapsto
bxb^{-1} $, corresponds to the automorphism $ T'_{\ast} $. The
fundamental group $ \pi_1 (\Gamma) $ is a quadratic extension of
$ \pi_1 (\bar \Gamma) $ by means of  the element $ b $, and this
extension is uniquely defined up to isomorphism  by the
automorphism $ T'_{\ast} $. Therefore
  $ \pi_1 (\Gamma) \simeq \D $, and hence the mapping
$ \eta_{\Gamma}: \Gamma \to K (\D, 1) $ is well defined.

It is easy to verify that the mapping $\eta_{\Gamma} \vert_
{\partial \Gamma}$ takes values in the subspace $K(\I_{b \times \bb},1)
\subset K(\D,1)$.
 The mapping $ \eta_{\Gamma} $, which is defined by the formula $(\ref{stG})$,
induces the mapping
\begin{eqnarray}\label{stN}
\eta_{\N\circ}: (\N_{\circ},
U(\partial \N)_{\circ}) \to (K(\D,1),K(\I_{b \times \bb},1)), 
\end{eqnarray} 
which we call the structure mapping.

 Also, it
is easy to verify that the homotopy class of the composition
$U(\partial \N)_{\circ} \stackrel{\eta_{\N\circ}}{\longrightarrow} K(\I_{b \times \bb},1)
\stackrel{p_b}{\longrightarrow} K(\I_{d},1)$, where the map 
$K(\I_{b \times \bb},1)
\stackrel{p_b}{\longrightarrow} K(\I_{d},1)$ 
is induced by the homomorphism
$\I_{b \times \bb} \to \I_d$ with the kernel $\I_{b}$,
$\partial \N(d) \stackrel{\eta}{\longrightarrow} K(\I_b,1)
\stackrel{p_b}{\longrightarrow} K(\I_d,1)$
is extended to a map on $\partial \N$ and this extension coincides to the
map $\kappa \circ res_{d}: \partial
\N(d) \to \RP^{n-k} \to K(\I_d,1)$, which is the composition of the resolution map
$res_d: \partial \N(d)
\to  \RP^{n-k}$ and the embedding of the skeleton $\RP^{n-k} \subset K(\I_d,1)$ in the
classifying space.

Let us assume that in the polyhedron 
$\N_{\circ}$ there is a marked component 
\begin{eqnarray}\label{Na}
\N_a \subset \N_{\circ} 
\end{eqnarray}
 and 
assume that the following mapping
\begin{eqnarray}\label{cycl}
\mu_{a}: \N_{a} \to K(\I_a,1)
\end{eqnarray}
determines the reduction of the restriction of the structure mapping $(\ref{stN})$ to the marked component
in the formula $(\ref{Na})$.

The following characteristic number 
\begin{eqnarray}\label{h mu_a}
\langle \mu_{a}^{\ast}(t);[\N_{a}]
\rangle,
\end{eqnarray}
is well defined, where   $t \in
H^{n-2k}(K(\I_a,1);\Z/2)$ is the generic cohomology class
$[\N_{a}]$ is the fundamental class of the polyhedron  $\N_{a}$ (assuming that the equivariant $PL$-mapping $d^{(2)}$
is transversal along $\Delta_{\R^n}$, this polyhedron is a $PL$-manifold).
\[  \]

\begin{definition}$\label{cyclicd}$
{\bf{Cyclic structure of a formal mapping $d^{(2)}$}}

Assume that a component $(\ref{Na})$  of the polyhedron $\N_{\circ}$ is marked and a mapping $(\ref{cycl})$ is well defined.
The mapping
 $(\ref{cycl})$ is called the cyclic structure of the equivariant mapping  with holonomic self-intersection  along
 the polyhedron $\N_{a}$,
  if the characteristic number $(\ref{h mu_a})$ satisfies the following equation:
  \begin{eqnarray}\label{h mu_a1}
\langle \mu_{a}^{\ast}(t);[\N_{a}]
\rangle =1.
\end{eqnarray}
\end{definition} 
\[  \]

 We need a criterion to verify that the mapping $\mu_{a}$
   satisfies the equation ($\ref{h mu_a1}$).

 
 
Let us consider the canonical 2-sheeted covering over the polyhedron $\N$, which is, probably, branched over the boundary:
\begin{eqnarray}\label{p_a}
p: \bar \N \to \N.
\end{eqnarray}
The total space of this covering is a closed polyhedron $\bar \N$ of the dimension
 $n-2k$. This polyhedron is decomposed into the union of the two subpolyhedra:
$\bar \N = \bar \N_{a} \cup \bar \N_{\b}$, this decomposition corresponds with the decomposition
in the formula $(\ref{Na})$.
The polyhedron  $\bar \N_{\a}$ is a closed $PL$--manifold of the dimension  $n-2k$, the polyhedron $\bar \N_{\b}$ 
is a compactification of an open  $PL$--manifold of the dimension $n-2k$ by means the boundary 
$\partial \N_{\b}$. The restriction of the involution $T_p$  $(\ref{p_a})$ on the boundary is free and a closed polyhedron 
$\bar \N_{\b} \cup_{T_p} T(\bar \N_{\b})$ is well-defined, denote this polyhedron by  $\N_{\b}/_{\sim}$. Denote $\bar \N_{\a} \cup  \N_{\b}/_{\sim}$ by
$\bar \N/_{\sim}$.

Let us consider the mapping
$p_{\I_c,\I_d} \circ \bar \eta: \bar \N/_{\sim} \to K(\I_c,1) \to
K(\I_d,1)$, where the mapping $p_{\I_c,\I_d}: K(\I_c,1) \to K(\I_d,1)$ is induced by the epimorphism 
$\I_c \to \I_d$, with the kernel, generated by the element   $ab \in \I_c$. Denote this mapping by 
$\bar \eta: \bar \N/_{\sim} \to K(\I_d,1)$. Denote the restriction of $\bar \eta$ on $\bar \N_{\b}/_{\sim}  \subset  \N/_{\sim}$ 
by
$\bar \eta_{\b}: \bar \N_{\b}/_{\sim} \to K(\I_d,1)$.

\begin{lemma}\label{lemma20}
Assume that the restriction of the structure mapping $\eta_{\circ}$ on the component $\N_{\a}$ as in Definition  $\ref{cyclicd}$
and a mapping $\mu_{a}$ determines a reduction of the structure mapping to the mapping into $K(\I_a,1) \subset K(\D,1)$.
The condition $(\ref{h mu_a})$ is a corollary  of the following two conditions: the homology class
\begin{eqnarray}\label{Nb}
\bar \eta_{\b\ast}([\bar \N_{\b}/_{\sim}]) \in  H_{n-2k}(K(\I_d,1);\Z/2) 
\end{eqnarray}
is trivial.
\end{lemma}

\subsubsection*{Proof of Lemma $\ref{lemma20}$}
Let us consider a sketch of the proof. The homology class $ \bar \eta_{\ast}([\bar \N/_{\sim}]) \in  H_{n-2k}(K(\I_d,1);\Z/2)$,
is the generator, because the fundamental class
 of the subpolyhedron $\bar \N/_{\sim} = Cl(\bar \N_{\circ}) \subset \RP^{n-k}$
is dual to the (normal) characteristic class of the dimension  $2k$, which is the generator  in $H_{n-2k}(\RP^{n-k})$, because $n=2^{\ell}-1$. 
Therefore the characteristic number 
 $(\ref{h mu_a1})$ is equal to 1 iff the class $(\ref{Nb})$ is trivial.
Lemma 
$\ref{lemma20}$ is proved. 
\[  \]

Theorem $ \ref{prop22} $ is
based on the application of the following principle of density of
the subspace of immersions in the space of continuous maps
equipped with the compact-open topology, see [Hi, Theorem 5.10].

\begin{proposition}\label{prop24}
Let $ f_0: M \looparrowright R $ (we will use the case $R=\R^n$) be a smooth immersion between
manifolds, where the manifold $M$ is compact, the manifold $ R $ is equipped with the
metric $ \dist $ and $ \dim (M) <\dim (R) $. Let $ g: M \to R $ be a
continuous mapping homotopic to the immersion $ f_0 $. Then $
\forall \varepsilon> 0 $ there exists an immersion $ f: M
\looparrowright R $, regularly homotopic to the immersion $ f_0 $,
for which $\quad \dist (g; f)_{C^0} < \varepsilon$ in the space of
maps with the induced metric.
\end{proposition}

\begin{proposition}\label{prop24bis}
Let
 $(M^{(2)},\partial M^{(2)})$ be a smooth manifold with boundary, assume that 
a free involution  
$T_{M^{(2)}}: (M^{(2)},\partial M^{(2)}) \to (M^{(2)},\partial M^{(2)})$ is well defined. 
Let 
 $R^{(2)}$ be a smooth manifold with the metric, denoted by $\dist$, equipped the standard involution  
$T_{R^{(2)}}: R^{(2)} \to R^{(2)}$, which is free outside the diagonal  $\Delta_R \subset R^{(2)}$ and $\dim(R^{(2)}) = 2n$, $\dim(\Delta_R)=n$, $\dim(M^{(2)}) < \dim(R^{(2)})$. Assume that there exists a
$(T_{M},T_R)$--equivariant immersion  $F_0^{(2)}: M^{(2)} \looparrowright R^{(2)}$, and the image of the boundary has no intersection with
the fixed point manifold: 
\begin{eqnarray}\label{uslF0}
Im(F_0^{(2)}(\partial M^{(2)})) \subset R^{(2)} \setminus \Delta_R. 
\end{eqnarray}
Let $G^{(2)}: M^{(2)} \to R^{(2)}$ be a continuous $(T_{M^{(2)}},T_{R^{(2)}})$--equivariant mapping, for which the following condition
is satisfied 
\begin{eqnarray}\label{uslG} 
Im(G^{(2)}(\partial M^{(2)})) \subset R^{(2)} \setminus \Delta_R. 
\end{eqnarray}
Moreover, assume that the equivariant mapping $G^{(2)}$ is equivariant homotopic to the immersion 
 $F_0^{(2)}$ in the space of mappings with the condition above. Then 
$\forall
\varepsilon > 0$ there exists a $(T_{M^{(2)}},T_{R^{(2)}})$--equivariant immersion 
$F_1^{(2)} : M^{(2)} \looparrowright R^{(2)} $, which is regular equivariant homotopic to the equivariant immersion  $F_0^{(2)}$, 
for which the analogous condition $(\ref{uslF0})$ is satisfied, and additionally, the condition  $\dist(F_1^{(2)};G^{(2)})_{C^0} < \varepsilon $ in the space of equivariant maps with induced metric.
\end{proposition}

\subsubsection*{A sketch of a proof of Proposition  $\ref{prop24bis}$}

Consider an equivariant triangulation of the manifold  $(M^{(2)},\partial M^{(2)})$ of the caliber much less then $\varepsilon$.
The proof is possible by the induction over the skeletons by analogous arguments as the Hirsch Theorem $\ref{prop24}$.
\[   \]

 We shall use Proposition 
 $\ref{prop24bis}$ as follows. Let $(f,\kappa,\Psi)$ be a skew-framed immersion
$f: M^{n-k} \looparrowright \R^n$. Consider an open manifold  $M^{n-k} \times M^{n-k} \setminus \Delta_M$, which is equipped with the standard  involution $T_M$, this involution is free outside the diagonal. Denote by 
$\bar{M}^{(2)}$ the spherical blow-up of the manifold $M^{n-k} \times M^{n-k} \setminus \Delta_M$, equipped with the free involution  $T^{(2)}_M$.
Denote by $M^{(2)}$ the quotient $\bar{M}^{(2)}/T^{(2)}_M$. The boundary $\partial M^{(2)}$ coincides with the projectivization   $TP(M^{n-k})$ of the tangent bundle $T(M^{n-k})$. Denote the manifold $M^{(2)} \setminus \partial M^{(2)}$ by $M^{(2)}_{\circ}$. 
Denote by $(\R^n)^{(2)}$ the manifold  $\R^n \times \R^n$, equipped with the standard involution  $T^{(2)}_{\R^n}$. 
The following mapping of the classifying spaces is well defined:
\begin{eqnarray}\label{Gauss}
\bar{f}^{(2)}: \bar{M}^{(2)} \to \R^n \times \R^n, 
\end{eqnarray}
This mapping is 
$(T^{(2)}_M,T^{(2)}_R)$--equivariant immersion. The equivariant immersion 
$(\ref{Gauss})$ satisfies the condition, which is analogous to the condition for the mapping $F_0$ from Proposition  $\ref{prop24bis}$. 

Let us calculate the normal bundle $\nu_{\bar{f}^{(2)}}$ of the immersion  $\bar{f}^{(2)}$, using the skew-framing $(\kappa,\Psi)$ of the immersion 
$f$. Evidently, the framing
$\Psi$ induces the isomorphism 
\begin{eqnarray}\label{kappa+kappa}
\bar{\Psi}^{(2)}: \nu_{\bar{f}^{(2)}} = k(\kappa_1 \oplus \kappa_2), 
\end{eqnarray}
where
$\kappa_i$ is the line bundle, which is induced by the immersion of the $i$-th factor, $i=1,2$. 
The involution $T^{(2)}_M$ is covered by an involution of the bundle $\nu_{\bar{f}^{(2)}}$ (this involution is not the identity over the bases of the bundles), which permutes the corresponding line factors in the right side of the formula  $(\ref{kappa+kappa})$.  
Therefore the vector bundle over  $M^{(2)}$ is well-defined, denote this bundle by $\nu_{f^{(2)}}$, the isomorphism 
$(\ref{kappa+kappa})$ induces the isomorphism 
\begin{eqnarray}\label{D}
\Psi^{(2)}: \nu_{f}^{(2)} = k(\eta), 
\end{eqnarray}
where $\eta$ is a 2-dimensional $\D$--bundle. 
 
Let us assume that the mapping   $(\ref{Gauss})$ is transversal along the diagonal  $\Delta_{\R^n} \subset \R^n \times \R^n$ 
and let us consider the inverse image  $(\bar{f}^{(2)})^{-1}(\Delta_{\R^n})$ of this diagonal, and denote this inverse image by 
$\bar N \subset \bar{M}^{(2)}$. Evidently, that $\bar N$ is a closed $n-2k$-dimensional manifold, which coincides with the manifold $(\ref{barN})$. The manifold  $\bar N$ is equipped with the standard free involution $T^{(2)}_M \vert_{\bar N}: \bar N \to \bar N$, denote the quotient 
with respect to this involution by
$N^{n-2k}$. The manifold $N^{n-2k}$ is closed and this manifold coincides with the manifold $(\ref{N})$. 
This new definition of the manifold is more general, because this definition is possible without an assumption that the equivariant immersion
 $\bar{f}^{(2)}$ is holonomic. We shall use this definition for $(T^{(2)}_M,T^{(2)}_R)$--equivariant immersions, which satisfy the condition $(\ref{uslF0})$ from the equivariant regular homotopy class of the equivariant immersion $\bar{f}^{(2)}$.
The immersion 
 $\bar N^{n-2k} \looparrowright M^{n-k}$ and the restriction $f \vert_{\bar N^{n-2k}}: \bar N^{n-2k} \looparrowright \R^n$ is well defined.
The normal bundle of the immersion  $f \vert_{\bar N^{n-2k}}$ is isomorphic to the restriction of the bundle  $\nu_{\bar{f}^{(2)}}$ 
over the submanifold 
$\bar N \subset \bar{M}^{(2)}$. 

Let us consider the immersion
 $g: N^{n-2k} \looparrowright \R^n$, which is defined by the formula  $(\ref{g})$. The normal bundle 
 $\nu_g$ of this immersion  $g$ is naturally isomorphic to the restriction of the bundle  $\nu_{f}^{(2)}$ over the submanifold   
$N^{n-2k} \subset M^{(2)}$. The dihedral framing, which is determined by the formula   $(\ref{D})$, coincides with the dihedral framing  $\Psi$, which is determined in Proposition  $\ref{prop5}$.

 Assume, that a $(T^{(2)}_{M},T^{(2)}_{\R^n})$--equivariant immersion  
$G^{(2)}: M^{(2)} \looparrowright \R^n \times \R^n$ is well defined,  this immersion satisfies the condition $(\ref{uslG})$,  and is regular homotopic to the immersion 
$(\ref{Gauss})$ in the space of equivariant immersions with the prescribed conditions, 	in particular, regular homotopies keep the isomorphism $(\ref{D})$. 
 Assume that the regular homotopy keeps the condition    $(\ref{uslG})$ and is transversal to 
$\Delta_{\R^n} \subset \R^n \times \R^n$.  Let us consider the manifold  $N^{n-2k}(G)$. This manifold is close, additionally
the immersion 
 $g(G): N^{n-2k}(G) \looparrowright \Delta_{\R^n} = \R^n$ is well defined as the restriction of the immersion  $G$.
The normal bundle of the considered immersion, which is denoted by   $\nu_{g(G)}$, is isomorphic to the restriction of the normal bundle 
of the immersion $G/T_M$. The normal bundle  $\nu_{g(G)}$ is equipped with a dihedral framing $(\ref{D})$.
Therefore a 
 $\D$--framed immersion  $(g(G),\eta(G),\Psi(G))$ is regular cobordant to a  $\D$--framed immersion 
$(g,\eta,\Psi)$ of the self-intersection manifold of the given skew-framed immersion   $(f,\kappa,\Xi)$.
The $\D$--framed immersion $(g,\eta,\Psi)$ satisfies additional properties, which are defined using the mapping 
$(\ref{Gauss})$.



\subsubsection*{Proof of Theorem  $\ref{prop22}$ from Lemma  $\ref{osnlemma1}$}

Consider the skew-framed immersion 
 $(f,\kappa,\Xi)$, which is determined in the statement of  Proposition  $\ref{prop22}$.
Consider the equivariant mapping  $d^{(2)}: \RP^{n-k}\times\RP^{n-k} \to \R^n\times\R^n$, which is constructed in Lemma  $\ref{osnlemma1}$.  
Consider the mapping $\kappa: M^{n-k} \to \RP^N$, $N >> n-k$. Without loss of a generality, we may assume that  $Im(\kappa) \subset \RP^{n-k} \subset \RP^{\infty}$. Consider the mapping  $\kappa \times \kappa:
 M^{n-k} \times M^{n-k} \to \RP^{n-k} \times \RP^{n-k}$. Denote by  $(T_{M},T_{\RP^{n-k}})$ the standard involutions in the image and in the preimage of the mapping $\kappa \times \kappa$. Denote the  $(T_{M},T_{\RP^{n-k}})$--equivariant mapping  $\kappa \times \kappa$ by $\kappa^{(2)}$. 
Let us consider the composition  
$d^{(2)} \circ \kappa^{(2)}:  M^{n-k} \times M^{n-k} \to \R^n \times \R^n$,
where the equivariant mapping   $d^{(2)}$ is constructed in Lemma $\ref{osnlemma1}$.
 
 Let us consider the restriction of the mapping
 $d^{(2)} \circ \kappa^{(2)}$ on the diagonal $\Delta_{M} \to \Delta_{\R^n}$.
 By Proposition 
 $\ref{prop24}$ the mapping $d^{(2)} \circ \kappa^{(2)}$ is arbitrary closed to an equivariant immersion from the equivariant regular homotopy class
 of the equivariant immersion  $f^{(2)}$. Therefore there exists an equivariant mapping 
$G^{(2)}: (M^{n-k})^{(2)} \to (\R^n)^{(2)}$, such that the restriction of this mapping on the interior of  $M^{n-k})^{(2)}_{\circ} \subset M^{n-k})^{(2)}$ is arbitrary closed to the restriction 
 $d^{(2)} \circ \kappa^{(2)}$ on $M^{n-k})^{(2)}_{\circ} = M^{n-k} \times M^{n-k} \setminus \Delta_{M^{n-k}}$, and which satisfies the condition 
 $(\ref{uslG})$. By Proposition  $\ref{prop24bis}$ the equivariant immersion with the prescribed properties 
exists.

 Consider the marked closed component 
 $\N_{a}$ of the polyhedron of the formal self-intersection of the formal mapping  $d^{(2)}$.
The closed component $N'_a \subset  (M^{n-k})^{(2)}$, which is mapped to the component  $\N_{a}$ by the mapping $\kappa^{(2)}$ is well-defined.
Therefore, if 
 $\varepsilon >0$ is sufficiently small,  the immersion  $G^{(2)}$ contains the component  $N_a^{n-2k}(G^{(2)})$, which is projected on
 $\N_{a}$, and the degree of this projection is equal to  $\deg^2(\kappa)=h_k \pmod{2}$, where the characteristic number  $h_k$ is given by 
Definition  $\ref{def1}$ for $(f,\kappa,\Xi)$.  

The following statement is evident.
\subsubsection*{Statement}
Assume that an arbitrary element
$x \in Imm^{sf}(n-k,k)$ is given, and a  $\D$--framed immersion 
$(g_0,\eta_0,\Psi_0)$ represents an element $y \in Imm^{\D}(n-2k,2k)$. 
Then there exists a skew-framed immersion 
$(f,\kappa,\Xi)$, which represents the element $x$, for which the self-intersection manifold contains a closed component, which represents $\D$-framed immersion  $(g_0,\eta_0,\Psi_0)$ up to regular homotopy.
\[  \]

Using the statement, we may assume that the self-intersection manifold
 $N_a^{n-2k}$ of the immersion $f$, which is equipped by a characteristic class  $\eta_{ N_a^{n-2k}}$, contains a closed component 
$N^{n-2k}(G^{(2)})_a \subset N^{n-2k}(G^{(2)})$. Denote by $(g,\eta,\Psi)$ a $\D$-framed immersion, which is the parametrization
of the self-intersection manifold of $f$. In particular, the restriction of $\eta$ of the component  
    $N^{n-2k}(G^{(2)})_a \subset N^{n-2k}(G^{(2)})$ was investigated above. 
This restriction admits a reduction mapping
$\eta_a: N^{n-2k}_a(G^{(2)}) \to K(\I_a,1)$, which is induced from the cyclic mapping of the cyclic structure of the mapping  $d^{(2)}$. 
The manifold $N^{n-2k}(G^{(2)})$ is equipped with a $\D$--framed immersion  $(g(G^{(2)}),\eta(G^{(2)}),\Psi(G^{(2)}))$, 
this $\D$-framed immersion  represents the prescribed element $[(g,\eta,\Psi)]$ in $Imm^{\D}(n-2k,2k)$.

In particular, the characteristic number 
 $(\ref{1.99})$, which is calculated for $(N_a^{n-2k},\eta_a)$, coincides with  
$h_k(f,\kappa,\Xi)$ by the property $(\ref{h mu_a1})$. Therefore, because  
$h_k(f,\kappa,\Xi) = h^{\D}(g,\eta,\Psi)$, the last component (probably, non-connected) $(N^{n-2k} \setminus N_a^{n-2k},\eta)$ 
has the trivial Hopf invariant. Theorem 
$\ref{prop22}$ is proved.

\section{Quaternionic structure of $PL$--mappings with singularities}

Let
$c: S^{n-2k}/\i \to \R^n$ be a $PL$--mapping in a general position.

Consider the configuration space
\begin{eqnarray}\label{99.1}
((S^{n-2k}/\i \times S^{n-2k}/\i) \setminus \Delta_{S^{n-2k}})/T',
\end{eqnarray}
which is called the пїЅdeleted squareпїЅ of the lens space $S^{n-2k}/\i$. This space is obtained
as the quotient of the direct product without the diagonal by the
involution $T': S^{n-2k}/\i \times S^{n-2k}/\i \to
S^{n-2k}/\i \times S^{n-2k}/\i$, exchanging the coordinates.
This space is an open manifold. It is convenient to define an analogous
space, which is a manifold with boundary.

Define the space $\bar \Gamma_1$ as a spherical blow-up of the
space $(S^{n-2k}/\i \times S^{n-2k}/\i) \setminus
\Delta_{S^{n-2k}}$ in the neighborhood of the diagonal. The
spherical blow-up is a manifold with boundary, which is defined as
a result of compactification of the open manifold $(S^{n-2k}/\i
\times S^{n-2k}/\i) \setminus \Sigma_{diag}$ by the fiberwise
glue-in of the fibers of the spherization  $ST \Sigma_{diag}$ of the
tangent bundle $T \Sigma_{diag}$ in the neighborhood of
zero--sections of the normal bundle of the diagonal $\Sigma_{diag}
\subset S^{n-2k}/\i \times S^{n-2k}/\i$. The following natural
inclusions are well defined:
$$S^{n-2k} \times S^{n-2k} \setminus \Sigma_{diag} \subset \bar \Gamma_1,$$
$$ST \Sigma_{diag} \subset \bar \Gamma_1.$$
On the space  $\bar \Gamma_1$ the free involution
\begin{eqnarray}\label{barT1}
\bar T': \bar \Gamma_1 \to \bar \Gamma_1
\end{eqnarray}
which is  an extension of an involution $ T' $ is well defined.

The quotient $\bar \Gamma_1 / \bar T'$
is denoted by $ \Gamma_1 $, and the corresponding double covering by
$$
p_{\Gamma_1}:  \bar \Gamma_1 / \bar T' \to \Gamma_1.
$$
The space $ \Gamma_1 $ is a manifold with boundary and it is called the resolution space of the
configuration space $ (\ref{99.1}) $. The projection $p_{\partial \Gamma_1}: \partial \Gamma_1 \to S^{n-2k}/\i$
is well defined, this map is called
a resolution of the diagonal.

For an arbitrary mapping $ c: S^{n-2k}/\i \to \R^n $ the polyhedron $\L (c) $ of self-intersection points of the mapping $ c $ is defined by the formula:
\begin{eqnarray}\label{Nc}
\L(c) = Cl \{ ([x,y]) \in int(\Gamma_1) : y \ne x, c(y)=c(x) \}.
\end{eqnarray}

By Porteous' Theorem [Por] under the assumption that the map $ c $ is
smooth and  generic, the polyhedron $ \L (c) $ is a manifold with
boundary of dimension $ n-4k $. This polyhedron is denoted by
$\L^{n-4k}(c)$ and called the polyhedron of self-intersection 
of the map $ c $. This formula $(\ref{Nc})$ defines an embedding
of polyhedra into manifold:
$$i_{\L(c)}: (\L^{n-4k}(c), \partial \L^{n-4k}(c)) \subset (\Gamma_1,\partial \Gamma_1). $$

The boundary $\partial \L^{n-4k}(c)$ of the manifold $\L^{n-4k}(c)$
is called the resolution manifold of critical
  points of the map $ c $. The map $p_{\partial \Gamma_1} \circ  i_{\partial \L(c)} \vert_{\partial \L(c)}:
 \partial \L^{n-4k}(c) \subset \partial \bar \Gamma_1  \to S^{n-2k}$ is called the
  resolution map of singularities of the map $ c $, we denote this mapping by
$res_{c}: \partial \L(c) \to S^{n-2k}/\i$.

 Consider the
canonical double covering
\begin{eqnarray}\label{pNc}
p_{\L(c)}: \bar \L(c)^{n-4k} \to \L(c)^{n-4k},
\end{eqnarray}
ramified over the boundary $\partial \L(c)^{n-4k}$ (over this
boundary the cover is a diffeomorphism). The next diagram is
commutative:
$$
\begin{array}{ccc}
i_{\bar \L(c)}: (\bar \L^{n-4k}(c), \partial \L^{n-4k}(c)) & \subset & (\bar \Gamma_1, \partial \Gamma_1) \\
\\
\downarrow p_{\L(c)}&   &  \downarrow  p_{\Gamma_1} \\
\\
i_{\L(c)}: (\L^{n-4k}(c), \partial \L^{n-4k}(c)) & \subset & (\Gamma_1, \partial \Gamma_1).
\end{array}
$$

\subsubsection*{Definition of the subgroup $\H \subset \Z/2^{[3]}$}
Consider the space $ \R^4 $ with the basis $(\e_1,\e_2,\e_3,\e_4)$.
The basis vectors are conveniently identified with the basic unit
quaternions $({\bf 1}, \i,\j,\k)$, which sometimes will be
used to simplify the formulas for some
transformations. Define the subgroup
\begin{eqnarray}\label{H}
\H \subset \Z/2^{[3]}
\end{eqnarray}
as a subgroup of transformations, of the following two types:

-- in each plane $Lin(\e_1 ={\bf 1},\e_2=\i)$,
$Lin(\e_3=\j,\e_4=\k)$ may be (mutually independent)
transformations by the multiplication with the quaternion $ \i $. The subgroup of
all such transformations is denoted by $ \H_c $, this subgroup is
isomorphic to $ \Z / 4 \times \Z / 4 $.

- the transformation exchanging each pair of the corresponding
basis vectors $\e_1 ={\bf 1}$ and $\e_3 =\j$ and
 the pair of the basis vectors $ \e_2 = \i $ and $ \e_4 = \k $, preserving their direction. Denote this transformation
 by
 \begin{eqnarray}\label{t}
t \in \H \setminus \H_c.
 \end{eqnarray}

It is easy to verify that the group itself is $ \H $, has order $ 32 $, and is a subgroup
$ \H \subset \Z / 2^{[3]} $
of
index 4.
\[  \]

\subsubsection*{Definition of the subgroup $ \H_{b \times \bb} \subset \H $ and the monomorphism
$i_{\I_a,\H}: \I_a \subset \H$}
Define the inclusion   $i_{\I_a,\H}: \I_a \subset \H$,
which translates the generator of the group $ \I_a $ into the operator of
multiplication by the quaternion $ \i $, acting simultaneously in each
plane $Lin(\e_1 ={\bf 1},\e_2=\i)$, $Lin(\e_3=\j,\e_4=\k)$.
Define the subgroup $ \H_{b \times \bb} \subset \H $ as the product  of the subgroup $i_{\I_a,\H}: \I_a \subset \H$
 and the subgroup generated by the generator
$ t \in \H $. It is easy to verify that the group  $ \H_{b \times \bb}
$ is of the order $8$ and this group is isomorphic to $ \Z / 4 \times
\Z / 2 $. The subgroup $ \I_a \subset \H_{b \times \bb} $ is of the  index 2.
The inclusion homomorphism $i_{\I_a,\H_{b \times \bb}}: \I_a \subset
\H_{b \times \bb}$ of the subgroup and the projection $p_{\H_{b \times \bb},\I_a}: \H_{b \times \bb} \to
\I_a$ are defined,
such that the composition  $\I_a \stackrel{i_{\I_a,\H_{b \times \bb}}}{\longrightarrow} \H_{b \times \bb}
\stackrel{p_{\H_a,\I_a}}{\longrightarrow} \I_a$
 is the
identity.
\[ \]

\subsubsection*{Structure map  $\zeta_{\circ}: \L_{\circ}(c) \to
K(\H,1)$}

Define the map $\zeta_{\Gamma_1}: \Gamma_1 \to K(\H,1)$, which we
call the structure mapping of the пїЅdeleted squareпїЅ. Note that
the inclusion $\bar \Gamma_1 \subset S^{n-2k}/\i \times
S^{n-2k}/\i$ induces an isomorphism of fundamental groups, since
the codimension of the diagonal $\Delta_{S^{n-2k}/\i} \subset
S^{n-2k}/\i \times S^{n-2k}/\i$  satisfies the inequality $ n-2k
\ge 3 $. Therefore, the following equality is satisfied:
\begin{eqnarray}\label{pi_1H}
\pi_1(\bar \Gamma_1) = H_1(\bar \Gamma_1;\Z/4)= \Z/4 \times \Z/4.
\end{eqnarray}

Consider the induced automorphism $\bar T'_{\ast}: H_1(\bar
\Gamma_1;\Z/2) \to H_1(\bar \Gamma_1;\Z/2)$, induced by the
involution $(\ref{barT1})$. Note that this automorphism is not the
identity and permutes the factors. Fix an isomorphism of the
groups $H_1(\bar \Gamma_1;\Z/2)$ and $\H_c$, which maps the
generator of the first (respectively second) factor of $ H_1 (\bar
\Gamma_1; \Z / 2) $, see (\ref{pi_1H}), into the generator,
defined by the multiplication by the quaternion $\i$ in the plane
$Lin({\bf 1},\i)$ (respectively, in the plane $Lin(\j, \k)$) and
is the identity on the complement.

It is easy to verify that the  automorphism of the conjugation
with respect to the subgroup $\H_c \subset \H$ by means of the
element $t \in \H \setminus \H_c$,  (in this formula the element $
t $ is given by the equation $(\ref{t})$), defined by the formula $
x \mapsto txt^{-1} $, corresponds to conjugation by means of 
the automorphism $\bar T'_{\ast} $. The fundamental group $ \pi_1
(\Gamma_1) $ is a quadratic extension of $ \pi_1 (\bar \Gamma_1) $
by means of  the element $ t $, and this extension is uniquely
defined up to isomorphism by the automorphism $ T'_{\ast} $.
Therefore
  $ \pi_1 (\Gamma_1) \simeq \H $, and hence the mapping
$ \zeta_{\Gamma_1}: \Gamma_1 \to K (\H, 1) $ is well defined.

It is easy to verify that the mapping $\zeta_{\Gamma_1} \vert_
{\partial \Gamma_1}$ takes values in the subspace $K(\H_{b \times \bb},1)
\subset K(\H,1)$.
 The mapping $ \zeta_{\Gamma_1} $
induces the map $ \zeta_{\circ}: (\L_{\circ} (c),
U(\partial \L(c))_{\circ}) \to (K (\H, 1), K (\H_{b \times \bb}, 1)) $, which we call the structure map.
(In the considered case the notion of the structure mapping is
analogous to the notion of the classifying mapping for
$\Z/2^{[3]}$--framed immersion.) 
Also, it is easy to verify that
the homotopy class of the composition $U(\partial \L(c)_{\circ})
\stackrel{\zeta_{\circ}}{\longrightarrow} K(\H_{b \times \bb},1)
\stackrel{p_{\H_{b \times \bb},\I_a}}{\longrightarrow} K(\I_a,1)$ is extended to $\partial \L(c)$
 and coincides with the
characteristic map $\eta \circ res_{c}: \partial \L(c) \to
S^{n-2k}/\i \to K(\I_a,1)$, which is the composition of the
resolution map $res_c: \partial \L(c) \to S^{n-2k}/\i$ and
the embedding of the skeleton $S^{n-2k}/\i \subset K(\I_a,1)$ in
the classifying space.

\begin{definition}$\label{quaterniond}$
{\bf{Quaternionic structure for  a mapping $c: S^{n-2k}/\i \to
\R^n$ with singularities}}

Let  $c: S^{n-2k}/\i \to \R^n$ be a map in general position,
having critical points, where $k \equiv 0 \pmod{2}$.
 Let  $\L(c)_{\circ}$ be the polyhedron of double self-intersection points of the map $c$
with the boundary $\partial \L(c)$.

Let us assume that the polyhedron  $\L(c)_{\circ}$ is the disjoint union of the two components
\begin{eqnarray}\label{Lc}
\L(c)_{\circ} = \L_{\Q} \cup \L_{\H_{b \times \bb}\circ},
\end{eqnarray}
where the polyhedron  $\L_{\Q}$ is closed, and the polyhedron
$\L_{\H_{b \times \bb}\circ}$, generally speaking, contains a regular neighborhood of the boundary
$U(\partial \L_{\H_{b \times \bb}})_{\circ} = U(\partial \L(c))_{\circ}$.

We
say that this map $c$ admits a (relative) quaternionic structure, if the
structure map $\zeta_{\circ} : (\L(c)_{\circ},U(\partial \L(c))_{\circ}) \to (K(\H,1),K(\H_{b \times \bb},1))$ is given by the composition:
\begin{eqnarray}
\lambda_{\L(c)}: \L_{\Q} \cup \L_{\H_{b \times \bb}\circ} \stackrel{\lambda_{\L_{\Q}} \cup \lambda_{H_{b \times \bb}\circ}}{\longrightarrow} 
K(\Q,1) \cup
K(\H_{b \times \bb},1)
\end{eqnarray}
$$
 \stackrel{i_{\Q,\H} \cup i_{\H_{b \times \bb},\H}}{\longrightarrow} K(\H,1) \cup K(\H,1)
\stackrel{Id \cup Id}{\longrightarrow} K(\H,1).
$$
\end{definition}
\[  \]

\begin{proposition}\label{prop23}
If $n=4k+n_{\sigma}$,  $n\ge 255$,  an arbitrary element of the
group $Imm^{\D}(n-2k,2k)$, in the image of the homomorphism
$\delta_k: Imm^{sf}(n-k,k) \to Imm^{\D}(n-2k,2k)$,
 is
represented by a $\D$--framed immersion $(g,\eta, \Psi)$,
admitting a quaternionic structure in the sense of
{\rm{Definition}} $\ref{def16}$.
\end{proposition}
\[  \]


 

\subsubsection*{Proof of Proposition  $\ref{prop23}$ from Proposition  $\ref{prop22}$ and Lemma 
$\ref{osnlemma2}$ (Section 5)}

Consider a skew-framed immersion 
 $(f,\kappa,\Xi)$, such that $\delta([(f,\kappa,\Xi)]) = [(g_1: N^{n-2k} \looparrowright
\R^n,\eta,\Psi)]$ represents the prescribed element in the group 
$Imm^{\D}(n-2k,2k)$. By Proposition  $\ref{prop22}$, without loss of the generality, we may assume that the immersion $f: M^{n-k}
\looparrowright \R^n$ admits a closed marked component   $g: N_a^{n-2k} \looparrowright \R^n$, equipped with the reduction mapping
 $\mu_a:
N_a^{n-2k} \to K(\I_a,1)$.

Let us consider the mapping 
 $c:
S^{n-2k}/\i \to \R^n$, which is constructed in Lemma $\ref{osnlemma2}$ and consider the immersion 
 $g_a: N^{n-2k}_a \looparrowright \R^n$, which is defined as a
$C^0$--small approximation of the composition $c \circ
\mu_a: N_a^{n-2k} \to S^{n-2k}/\i \to \R^n$ in the prescribed regular homotopy class of the restriction of the immestion  $g$. 

Let us consider the formal (holonomic) mapping of $2$-configuration spaces. Analogously to Theorem
 $\ref{prop22}$ a quaternionic structure 
of the $\D$--framed immersion  $[(g,\eta,\Psi)]$ is well defined, this structure is induced from the quaternionic structure 
of the mapping $c$. Proposition 
$\ref{prop23}$ is proved.

\section{Cyclic and quaternionic
structure for formal (equivariant) mappings}

\begin{lemma}\label{osnlemma1}

Assuming the dimension restriction
\begin{eqnarray}\label{dimdimdim}
k \ge 17, \quad n - k \equiv 0\pmod{4}, 
\end{eqnarray}
 there exists a formal (equivariant) mapping $d^{(2)}$,
which admits a cyclic structure in the sense of Definition $\ref{cyclicd}$.
\end{lemma}

\begin{lemma}\label{osnlemma2}
Assuming the dimension restriction
$n=4k+(2^{\sigma}-1)$, $n=2^{\ell}-1$, $\ell \ge 8$, $\sigma =
\left[ \frac{\ell}{2}\right]-1$, there exists a generic $PL$-mapping  $d_1: S^{n-2k}/\i \to \R^n$ with a singularity, 
which admits a relative quaternionic structure in the sense of Definition $\ref{quaterniond}$.
\end{lemma}

In this section we prove Lemmas 
$\ref{osnlemma1}$-$\ref{osnlemma2}$  from a unified point of view.
 The possibility of
such an approach in the case of cyclic structure was discovered by
Prof. A.V.Chernavsky at the end of the last century, and by Dr.
S.A.Melikhov (2005) in the case of quaternionic structure.
Preliminary results for cyclic and $\H_{b \times \bb}$--structure
in the case of weaker restrictions on the codimension of the
immersion, are given in the papers  \cite{A1}.

Note that the cases $n=15$, $n=31$, $n=63$, and $n=127$
were not considered. The construction of quaternionic structure in Lemma
$\ref{osnlemma2}$ does not require the Massey  embedding
$S^3/\Q \subset \R^4$ \cite{Ma}, see also \cite{Me}.  Such an embedding was known earlier
to W.Hantzsche \cite{He}. By means of such
an approach, it might be possible to weaken the dimensional
restrictions in Lemma $\ref{osnlemma2}$.
For example, the Massey embedding allows to generalize Lemma
$\ref{osnlemma2}$ for maps in the range $\frac{4}{5}$ (for maps $M^{m} \to \R^n$, $\frac{m}{n} \le \frac{4}{5}$).
This means
that one may consider an extra two quadratic extensions of the
quaternionic group as the structure group of framing of
immersions.

Additional arguments, in particular, might yield a proof of
the last cases in the Adams Theorem on Hopf invariants,
and clarify the remaining case in dimension $126$ not covered by the
Hill-Hopkins-Ravenel Theorem on Kervaire invariants.

\subsection{Auxiliary mappings}
We start by
construction of auxiliary mappings. In Lemma  $\ref{osnlemma1}$ 
there are axillary mappings $\hat c$, $c$ for the mapping
 $d$; in Lemma $\ref{osnlemma2}$
there are axillary mappings $c_1$, $\tilde c_1$  for the mapping
 $d_1$.

The transformation in Lemma
$\ref{osnlemma1}$ to the required formal (equivariant) mapping $d^{(2)}$ from the mapping $c$
 is given by an approximation, which is constructed in Lemma  $\ref{Ycirc}$. 

To proof the mentioned lemmas and propositions we introduce on the
singular set of auxiliary mappings the coordinate system called
\emph{angle-momentum}. By means of this coordinate system in
Lemmas $\ref{lemma280}$,$\ref{lemma291}$.
The configuration space in Lemma $\ref{lemma280}$  is defined as finite-dimensional  resolution
spaces for the singularity of the mapping $c$. In
Lemma $\ref{lemma291}$ the resolution
spaces is much simpler, because the mapping under investigation is
close to stable. 

\subsubsection*{Construction of an axillary mapping  $c: \RP^{n-k} \to
\R^n$, $\hat c: S^{n-k}/\i \to \R^n$ in Lemma $\ref{osnlemma1}$}


The mapping $ p': S^{n-k} \to J $ is well defined as the join of
$r$ 
copies of the standard 4-sheeted coverings $S^1 \to S^1 / \i$. The
standard action $\I_a \times S^{n-k} \to S^{n-k}$ commutes with
the mapping $p'$. Thus, the map $ \hat p: S^{n-k} /\i \to J$ is
well defined and the map $ p: \RP^{n-k} \to J$ is well defined
as the composition $\hat p \circ \pi: \RP^{n-k} \to J$ of the standard double covering $ \pi:
\RP^{n-k} \to S^{n-k} / \i$ with the map $ \hat p$.

The required mapping $c$ is defined by the formula 
\begin{eqnarray}\label{c}
i_{J} \circ p : \RP^{n-k} \to J \subset \R^n. 
\end{eqnarray}
The required mapping $\hat c$ is defined by the formula 
\begin{eqnarray}\label{hatc}
i_{J}
\circ \hat p : S^{n-k}/\i \to \R^n.
\end{eqnarray}

\subsubsection*{Construction of axillary mappings $c_1: S^{n-2k+2^{\sigma-1}}/\i \to \R^n$,
$\tilde c_1: S^{n-2k}/\i \to \R^n$ }

Let a positive integer parameter $k$ and a positive integer $n$
are given as in Lemma $\ref{osnlemma2}$.
Let us denote by $J_1$ a $(n-2k+2^{\sigma-1})$--dimensional
polyhedron (the equation  $n-2k+2^{\sigma-1} =
\frac{n-1}{2} + 2^{\sigma}$ is satisfied), this polyhedron is defined as the
join of
\begin{eqnarray}\label{r_1}
\frac{n+1}{2^{\sigma+1}}+1=r_1
\end{eqnarray}
copies of the standard quaternionic lens space
$S^{2^{\sigma}-1}/\Q$. Below we shall used the following notation
$n_{\sigma} = 2^{\sigma}-1$, as in  [A1] and $m_{\sigma} =
2^{\sigma}-2$, as in  [A2]). By the Hirsch Theorem an embedding
 $i_{\Q}: S^{n_{\sigma}}/\Q \subset \R^{n_{\sigma}-3}$ is well defined.

Assuming $n=4k+2^{\sigma}-1$, $\ell \ge 7$ the embedding $J_1
\subset \R^n$, as the join of $r_1$ copies of the embedding
$i_{\Q}$, is well defined; let us denote this embedding by
$i_{J_1}: J_1 \subset \R^n$ (comp. with the mapping in [Lemma 35, A2].

The mapping $p'_1: S^{n-2k+n_{\sigma-1}-1} \to J_1$ is well
defined as the join of $r_1$ copies of the standard coverings
$S^{n_{\sigma}} \to S^{n_{\sigma}}/\Q$. The action $\Q \times
S^{n-2k+n_{\sigma-1}-1} \to S^{n-2k+n_{\sigma-1}-1}$ is
well defined as the standard diagonal action on the join of $r_1$ copies of the standard action on $S^3$, 
this action commutes with the mapping $p'_1$.

 Thus, the map $\hat p_1: S^{n-2k+n_{\sigma-1}-1}/\Q \to J_1$ is well defined and the map
 \begin{eqnarray}\label{p_1}
p_1 \cong \hat p_1 \circ \pi_1: S^{n-2k+n_{\sigma-1}-1}/\i \to
J_1,
\end{eqnarray}
  as the composition of the standard double
covering $\pi_1:
S^{n-2k+n_{\sigma-1}-1}/\i \to S^{n-2k+n_{\sigma-1}-1}/\Q$ with the map $ \hat
p_1$.

Define the required mapping
 $c_1$ as the composition
$i_{J_1} \circ p_1: S^{n-2k+n_{\sigma-1}-1}/\i \to
S^{n-2k+n_{\sigma-1}-1}/\Q \to J_1 \subset \R^n$. Consider the submanifold  $i: S^{n-2k}/\i \subset
S^{n-2k+n_{\sigma-1}-1}/\i$, this submanifold is in  general position with respect to strata
of the manifold $S^{n-2k+n_{\sigma-1}-1}/\i$, the strata are determined by the join structure.
Define the mapping
\begin{eqnarray}\label{tildep_1}
\tilde p_1 \cong \hat p_1 \circ \pi_1 \circ i: S^{n-2k}/\i \subset
S^{n-2k+n_{\sigma-1}-1}/\i \to J_1.
\end{eqnarray}
Define the required mapping
 $\tilde c_1$ as the composition 
 \begin{eqnarray}\label{tildec_1}
\tilde c_1: S^{n-2k}/\i \subset S^{n-2k+n_{\sigma-1}-1}/\i
\stackrel{c_1}{\longrightarrow} \R^n.
\end{eqnarray}

\subsubsection*{Subspaces and factorspaces of the 2-configuration space for
$ \RP^{n-k} $, related with the axillary mappings $c$, $\hat c$ in Lemma $\ref{osnlemma1}$}

The space $ \Gamma $,
the subspace $
\Gamma_{ \circ} \subset \Gamma $, its double coverings $
\bar \Gamma $, $\bar \Gamma_{\circ}$ were defined above. The structural mapping $\eta_{\Gamma_{ \circ}}: \Gamma_{ \circ} \to K (\D, 1)$
also were defined.

Denote by
\begin{eqnarray}\label{Sigmacirc}
\Sigma_{\circ} \subset \Gamma_{\circ}
\end{eqnarray} 
the polyhedron of double-points singularities of the map $ p:
\RP^{n-k} \to J $, this polyhedron is defined by the formula
 $ \{[(x, y)] \in \Gamma_{ \circ}, p(x) = p(y),
x \ne y \} $. This polyhedron is equipped with a structural mapping
$\eta_{\Sigma_{\circ}}: \Sigma_{\circ} \to K(\D,1),$
which is induced by the restriction of the structural mapping $
\eta_{\Gamma_{ \circ}} $ on the subspace $ \Sigma_{\circ} $.

The standard  free involution $ \i: \RP^{n-k} \to \RP^{n-k} $ is well
defined. This involution permutes points in each fiber of the
standard double covering $ \RP^{n-k} \to S^{n-k} / \i $. The space
$ \bar \Gamma_{ \circ} $ admits an involution (with fixed points)
\begin{eqnarray}\label{Ticirc}
 T_{\i\circ}
: \bar \Gamma_{ \circ} \to \bar \Gamma_{ \circ} ,
\end{eqnarray}
which is defined as the restriction of an involution 
$\i \times \i: \RP^{n-k} \times \RP^{n-k} $, constructed by the involution $
\i $ on each factor, on the subspace $ \bar \Gamma_{ \circ}
\subset \RP^{n-k} \times \RP^{n-k} $. On the quotient $\bar \Gamma_{\circ}/T=\Gamma_{\circ}$
of $ \Gamma_{\circ} $ by the another involution $ T $, which permutes the coordinates, the factorinvolution $
T_{\i \circ}: \Gamma_{ \circ} \to \Gamma_{ \circ} $ is
well defined.

Let us denote by $ \Sigma_{antidiag} \subset \Gamma_{\circ} $ a
subspace, called the antidiagonal, which is formed by all
antipodal pairs $ \{[(x, y)] \in \Gamma_{\circ}: x,y \in
\RP^{n-k}, x \ne y, \i (x) = y \} $. It is easy to verify that the antidiagonal $
\Sigma_{antidiag} \subset \Gamma_{\circ} $ is the set of fixed
points for the involution $ T_{\i\circ} $.


The subpolyhedron $ \Sigma_{\circ} \subset \Gamma_{\circ} $ of
multiple-points of the map $ p $ is represented by a union $
\Sigma_{ \circ} = \Sigma_{antidiag} \cup K_{\circ} $, where
$K_{ \circ} $ is an open  subpolyhedron contains all points of $
\Sigma_{\circ} $ outside the antidiagonal. The subpolyhedron $
K_{\circ} \subset \Gamma_{K_{ \circ}} $ is invariant under the
involution $ T_{\i\circ}$.

Define the restriction of the involution
 $T_{\i\circ} \vert_{K_{\circ}}$ by
$T_{K_{\circ}}$. 
The considered restriction is a free involution. Denote the factorspace 
 $K_{\circ}/ T_{K_{\circ}}$ by  $\hat K_{\circ}$.
 The restriction of the structure mapping
$\eta_{\Gamma_{\circ}}: \Gamma_{\circ}
\to K(\D,1)$ on  $K_{\circ}$  denote by 
 $\eta_{K_{\circ}}$.

Denote the closure of $ Cl (K_{
\circ}) $ of the polyhedron  $ K_{ \circ} $ (respectively, the
closure of the polyhedron $ Cl (\hat K_{ \circ}) $ polyhedron $
\hat K_ ( \circ) $) 
by $
K $ (respectively, by $ \hat K $).
Denote by
$ Q_{diag}$ the space
$\partial \Gamma_{diag} \cap K$. Obviously, $ Q_{diag} \subset K $.
We shall call this subspace  the component of the boundary of the polyhedron $ K $.
Similarly, we denote by $ \hat Q_{diag} $ the component of the boundary of the polyhedron
$ \hat K $.


Note that the mapping  $\eta_{K_{\circ}}$ is not expendable to
boundary component $ Q_{diag} $. The mapping $ \kappa_{diag}:
Q_{diag} \to K (\I_d, 1) $ is well defined. Let us denote by $U(Q_{diag})_{\circ} \subset K_{\circ} $ a small regular deleted
neighborhood of  $Q_{diag} $. The projection $ proj_{diag}: U(Q_{diag})_{\circ} \to Q_{diag} $ of the regular deleted
neighborhood to  $Q_{diag} $. The restriction of the structural
mapping $ \eta_{K_{\circ}} $ to the neighborhood  $U(Q_{diag})_{\circ}$ is represented by a composition of the map $ \eta_{U(Q_{diag})_{\circ}}: U(Q_{diag})_{\circ} \to K(\I_b, 1) $ and the maps
$ i_{\I_b, \D}: K(\I_b, 1) \to K (\D, 1) $.
Homotopy classes of
maps $ \eta_{K} \vert_{Q_{diag}} $ and $ \eta_{U(Q_{diag})_{\circ}}$ satisfy the equation:
$$\eta_{diag} \circ proj_{diag} = p_{\I_b,\I_d} \circ \eta_{UQ_{diag\circ}}.$$

Let us investigate the polyhedron of singularities of an
axillary mapping $\hat c$. define the following commutative diagram of subgroups:
\begin{eqnarray}\label{HH}
\begin{array}{ccccccc}
&&&& \I_{b \times \bb} &&\\
&& \nearrow && \cap &&\\
\I_d& \subset &\I_a& \subset & \D & \subset & \E . \\
&& \searrow && \cap &&\\
&&&& \I_c &&\\
\end{array}
\end{eqnarray}

In this diagram, the inclusion
$\D \subset \E$ is a central quadratic extension of $\D$ by the element $\i$ (of the order 4),
for which  
$\i^2$ coincides with the generator $-1$ of the subgroup $\I_d \subset \D$.
The abelian groups $\E_a,
\E_{b \times \bb}, \E_c, \E_d$ are the subgroups in  $\E$, this groups are the quadratic extensions of the corresponding subgroups $\I_a, \I_b \times \II_b, \I_c, \I_d$  by means of the element $\i$. 
Note that the groups $\E_{b \times \bb}$ and $\H_{b \times \bb}$ (see above  [formula (84), A2]) are isomorphic.

The difference between the considered  groups $\E_{b \times \bb}$ and $\H_{b \times \bb}$ are the following: the representation of
$\E_{b \times \bb} \to \Z^{[3]}$ (see below [Example 16, A1]) and 
$\H_{b \times \bb} \to \Z/2^{[3]}$ (see [Diagram (85), A2]) are different. 
The kernel of the epimomorphism 
\begin{eqnarray}\label{epiE}
\E_{b \times \bb} \to \Z/2^{[3]} \to \Z/2, 
\end{eqnarray}
where $\Z/2^{[3]} \to \Z/2$ corresponds to the 
subgroup [(19),A2] of the index 2, contains an element $\i \in \E_d \subset \E_{b \times \bb}$ of the order 4
(comp. with Diagram $(\ref{140})$ below, in which $\E_d = \E_c \cap \E_{b \times \bb}$).
The kernel of the homomorphism $\H_{b \times \bb} \to \Z/2^{[3]} \to \Z/2$ coincides
with the subgroup $\I_{b \times \bb} \subset \E_{b \times \bb}$, which is an elementary 2-group.

The induced automorphism
$\chi^{[3]}: \Z/2^{[3]} \to \Z/2^{[3]}$ of the group $\E_{b \times \bb}$,
re-denoted by 
\begin{eqnarray}\label{hatchiE}
\hat \chi^{[2]}: \E_{b \times \bb} \to \E_{b \times \bb}, 
\end{eqnarray}
is defined by the formula $\hat \chi^{[2]}(\i) = \i$, where 
$\i \in \E_d$--is the generator.

The following natural mapping
$\eta_{\hat K_{\circ}}: \hat K_{\circ} \to K(\E,1)$, which corresponds to
the mapping of canonical 2-sheeted covering, is well-defined:

\begin{eqnarray}\label{140}
\begin{array}{ccccccc}
\bar K_{\circ} & \stackrel{\bar r}{\longrightarrow} & \tilde{K}_{\circ}& &K(\I_c,1)& \longrightarrow & K(\E_c,1)\\
\downarrow & & \downarrow & \longrightarrow & \downarrow & & \downarrow\\
K_{\circ} & \stackrel{r}{\longrightarrow} & \hat K_{\circ}& &K(\D,1)& \longrightarrow & K(\E,1).\\
\end{array}
\end{eqnarray}

Horizontal maps between the spaces of the diagrams we re-denote
for brevity by $ \bar \eta, \check \eta, \eta, \hat \eta $,
respectively.

\subsubsection*{Subspaces and factorspaces of the 2-configuration space for $S^{n-2k}/\i$,
related with the axillary mapping $c_1$}

The space $ \Gamma_1 $, its double covering $ \bar \Gamma_1 $, and
the structural map $ \eta_{\Gamma_1}: \Gamma_1 \to
K (\E, 1) $ was defined in [A1, Section 4, (62) and below]. The space $ \Gamma_1 $ is a manifold with boundary. Denote the interior of this manifold
by $ \Gamma_{1 \circ} $. The restriction of the structural map $ \eta_{\Gamma_1} $ to $ \Gamma_{1 \circ} $
will be denoted by $ \eta_{\Gamma_{1 \circ}}: \Gamma_{1 \circ} \to
K (\E,1) $.

Denote by $ \Sigma_{1 \circ} \subset \Gamma_{1 \ circ} $ the
polyhedron of double-points singularities of the map $ p: S^{n-2k}
\to J_1 $, this polyhedron is obtained by the blowing up of the
polyhedron $ \{[(x, y)] \in \Gamma_{1 \circ}, p(x) = p(y), x \ne y
\} $. This polyhedron is equipped with a structural mapping
\begin{eqnarray}\label{strukt2}
\zeta_{\Sigma_{1\circ}}: \Sigma_{1\circ} \to K(\E,1),
\end{eqnarray}
which
is induced by the restriction of the structural mapping $
\zeta_{\Gamma_{1 \circ}} $ on the subspace $ \Sigma_{1 \circ} $.
 
 The subpolyhedron $ \Sigma_{1 \circ} \subset \Gamma_{1 \circ} $ of
multiple-points of the map $ p_1 $ is represented by a union $
\Sigma_{1 \circ} = \Sigma_{antidiag} \cup K_{1 \circ} $, where
$K_{1 \circ} $ is an open  subpolyhedron, this subpolyhedron
contains all points of $ \Sigma_{1 \circ} $ outside the
antidiagonal. Let us denote the restriction of the structural
mapping $ \zeta_{\Gamma_{1 \circ}}: \Gamma_{1 \circ} \to K (\E, 1)
$ on $ \Gamma_{K_{1 \circ}} $ and on $ K_{1 \circ} $ by $
\zeta_{\Gamma_{1 \circ}} $ and by $ \zeta_{K_{1 \circ}} $
respectively.

 Denote the closure of $ Cl (K_{1
\circ}) $ of the polyhedron  $ K_{1 \circ} $ in $ \Gamma_{1} $ (respectively, the
closure of the polyhedron $ Cl (\hat K_{1 \circ}) $ polyhedron $
\hat K_{1 \circ} $ in $\hat \Gamma_{1} $) by $
K_{1} $ (respectively, by $ \hat K_{1} $).
  Denote by
$ Q_{antidiag}$ the space $\Sigma_{antidiag} \cap K_{1} $, denote by $ Q_{diag}$ the space
$\partial \Gamma_{diag} \cap K_{1} $. Obviously, $ Q_{diag} \subset K_{1} $,
$ Q_{antidiag} \subset K_{1} $. We shall call these subspaces  the components of the boundary of the polyhedron $K_{1}$.

Note that the structural mapping of $ \zeta_{K_{1 \circ}} $ is extended from $ K_{1 \circ} $ to the component
$ Q_{antidiag} $ of the boundary. Denote this extension by $ \zeta_{Q_{antidiag}}: Q_{antidiag} \to K(\E, 1) $.
The mapping $ \zeta_{Q_{antidiag}} $ is the composition
$ \zeta_{antidiag}: Q_{antidiag} \to K (\Q, 1) $ and the
inclusion $ i_{\Q, \E}: K (\Q, 1) \subset K (\E, 1) $.

Note that the mapping $ \zeta_{K_{1}} $ is not expendable to
boundary component $ Q_{diag} $. The mapping $ \zeta_{diag}:
Q_{diag} \to K (\I_a, 1) $ is well defined. Let us denote by $
U(Q_{diag})_{ \circ} \subset K_{1 \circ} $ a small regular deleted
neighborhood of  $Q_{diag} $. The projection $ proj_{diag}: U
(Q_{diag})_{ \circ} \to Q_{diag} $ of the regular deleted
neighborhood to $Q_{diag} $  to the central manifold is well defined.

The restriction of
the structural mapping $ \zeta_{K_{1 \circ}} $ to the neighborhood
$ U(Q_{diag})_{ \circ} $ is represented by a composition of the map $
\zeta_{UQ_{diag \circ}}: UQ_{diag \circ} \to K(\E_{b \times \bb}, 1) $ and
the maps $ i_{\E_{b \times \bb}, \E}: K(\E_{b \times \bb}, 1) \to K (\E, 1) $.

Homotopy
classes of maps $\zeta_{diag} $ and $\zeta_{UQ_{diag\circ}}$ are related by the equation: 
$$ \zeta_{diag} \circ proj_{diag} = p_{\E_{b \times \bb}, \I_a} \circ \zeta_{UQ_{diag\circ}}.$$
\[  \]

\section{Resolution spaces for singularities}

\subsection*{Resolution spaces for polyhedra  $K_{\circ}$  and $\hat K_{\circ}$}

We construct a space $RK_{\circ}$, which we call the
resolution space of the polyhedron  $K_{\circ}$. In [A2] the group $(\I_b \times \II_b)_{\chi^{[2]}} \Z$, equipped with the homomorphism
$\Phi^{[2]}: (\I_b \times \II_b)_{\chi^{[2]}} \Z \to \D$,
and the subgroup $\I_b \times \II_b \subset (\I_b \times \II_b)_{\chi^{(2)}} \Z$ are well defined. 

Consider the following diagrams: 

\begin{eqnarray}\label{16.2}
\begin{array}{ccc}
 RK_{\circ} & \stackrel {pr}{\longrightarrow} & K_{\circ}  \\
&&\\
  \phi \downarrow & & \\
&&\\
K((\I_b \times \II_b)_{\chi^{(2)}} \Z,1), &  & \\
\end{array}
\end{eqnarray}

\begin{eqnarray}\label{118.2}
\begin{array}{cccc}
RQ_{diag \circ} & \qquad \stackrel{pr}{\longrightarrow} \qquad & UQ_{diag \circ}\\
\phi \searrow & & \swarrow \eta_{diag \circ } \\
& K(\I_b \times \II_b,1), & 
\end{array}
\end{eqnarray}
where  $RQ_{diag \circ} = 
(pr)^{-1}(
UQ_{diag \circ})$.

\begin{lemma}\label{lemma28}
There exists the space
$RK_{\circ}$, which is included into the commutative diagram  $(\ref{16.2})$. 
The following diagram
$(\ref{118.2})$ determines the boundary conditions.
\end{lemma}

\subsubsection*{Resolution spaces for polyhedra   $\Sigma$ and $\hat K$}

Define a space 
$R\Sigma_{\circ}$, which is called the resolution space for the polyhedron  $\Sigma_{\circ}$, 
which is given by the formula  $(\ref{Sigmacirc})$.

The space
$R\Sigma_{\circ}$ contains two components, which is denoted by  $R\Sigma_{a}$, $RK_{b \times \bb\circ}$:
 \begin{eqnarray}\label{RK0}
R\Sigma_{a} \cup RK_{b \times \bb\circ} = R\Sigma_{\circ}.
\end{eqnarray}

The space
$R\Sigma_{a}$ is a closed polyhedron, for which the structured mapping
\begin{eqnarray}\label{phia}
\phi_a: R\Sigma_a \to K(\I_a,1)
\end{eqnarray}
is well-defined. The mapping
 $(\ref{phia})$ is included into the following commutative diagram:

\begin{eqnarray}\label{16.20}
\begin{array}{ccc}
\Sigma_{\circ} & \stackrel {pr}{\longleftarrow}& R\Sigma_{a}\\
&&\\
\downarrow \eta_{\circ} & & \downarrow \phi_a \\
&&\\
K(\D,1) & \supset &  K(\I_a,1).\\  
\end{array}
\end{eqnarray}

The space
$RK_{b \times \bb\circ}$  is a 2-sheeted covering space of the covering
 $Rr_{b \times \bb} : RK_{b\times \bb\circ} \to 
R\hat K_{b\times \bb\circ}$.

\begin{eqnarray}\label{16.23}
\begin{array}{ccccccc}
K_{\circ} & \stackrel {pr}{\longleftarrow}& RK_{b\times \bb\circ}&  \stackrel {Rr_{b\times \bb}}{\longrightarrow} &  R\hat K_{b\times \bb\circ}& \stackrel {p \hat r}{\longrightarrow} & \hat K_{\circ}  \\
&&&&&&\\
& & \downarrow \hat \phi_{b\times \bb} & & \downarrow \phi_{b\times \bb} &  & \\
&&&&&&\\
&& K(\I_{b \times \bb} \int_{\chi^{[2]}} \Z,1) & \subset &  K(\E_{b \times \bb} \int_{\hat \chi^{[2]}} \Z,1). &&  
\end{array}
\end{eqnarray}
The group
  $\E_{b \times \bb} \int_{\hat \chi^{[2]}}$, which is used in Diagram $(\ref{16.23})$
 is defined analogously to the group  $(\H_{b \times \bb}) \int_{\chi^{[3]}} \Z$,
 [Formula (68), A2], using the automorphism (involution)  $(\ref{hatchiE})$.

Denote
$(p \hat r)^{-1}(\hat
UQ_{diag\circ})$ by $R\hat Q_{diag\circ}$. The following inclusion 
$R\hat Q_{diag\circ} \subset R\hat K_{b \times \bb\circ} $ is well-defined.

Let us denote by $RQ_{diag\circ}$ the boundary of the corresponding 2-sheeted covering space over
 $R\hat Q_{diag\circ}$. The following diagram is well-defined.


\begin{eqnarray}\label{118.20}
\begin{array}{ccc}
R\hat Q_{diag\circ} & \qquad \stackrel{p\hat r}{\longrightarrow} \qquad & U\hat Q_{diag\circ} \\
&&\\
\hat \phi_{b \times \bb} \downarrow & &\hat \eta_{diag\circ} \downarrow \\
&&\\
K(\E_{b \times \bb} \int_{\chi^{[2]}} \Z,1) & \supset & K(\E_{b \times \bb},1). \\
\end{array}
\end{eqnarray}

To prove the main result of the section we will use the following lemma.

\begin{lemma}\label{lemma280}
There exists a space $R\Sigma_{\circ}$, which
is satisfies the equation $(\ref{RK0})$.

The component  $R\Sigma_a$ is equipped by the mapping
$(\ref{phia})$, which is included into the commutative diagram $(\ref{16.20})$.

The component
$RK_{b\times \bb\circ}$ is the total space of a regular 2-sheeted covering over the space 
 $R\hat K_{b\times \bb\circ}$ such that the commutative diagram  $(\ref{16.23})$
 is well-defined. Moreover, the commutative diagram
$(\ref{118.20})$, which determines boundary conditions, is well-defined. 
\end{lemma}

\subsubsection*{Resolution space for the polyhedron $\Sigma_1$}

We shall define a space
 $R\Sigma_{1\circ}$, which we call resolution
space of the polyhedron $\Sigma_1$.
The space $R\Sigma_{1\circ}$ contains two components, which is denoted by  $R\Sigma_{\Q}$, $RK_{\H_{b \times \bb}\circ}$, as follows: 
\begin{eqnarray}\label{RK1}
R\Sigma_{\Q} \cup RK_{\E_{b \times \bb}\circ} = R\Sigma_{1\circ}.
\end{eqnarray}

Let us consider the following diagrams:

\begin{eqnarray}\label{16.21}
\begin{array}{ccc}
R\Sigma_{\Q} \cup RK_{\E_{b \times \bb}\circ} & \stackrel {pr_1}{\longrightarrow} & \Sigma_1  \\
&&\\
  \phi_1 \downarrow & & \\
&&\\
K(\Q,1) \cup K(\E_{b \times \bb},1), &  & \\
\end{array}
\end{eqnarray}


\begin{eqnarray}\label{118.21}
\begin{array}{ccc}
RQ_{diag} & \qquad \stackrel{pr_1}{\longrightarrow} \qquad & Q_{diag} \\
\phi_1 \searrow & & \swarrow \zeta_{diag} \\
& K(\E_{b \times \bb},1), &
\end{array}
\end{eqnarray}
in which $RQ_{diag} = (pr_1)^{-1}(
Q_{diag})$.

The following lemma is analogous to Lemma 
$\ref{lemma280}$

\begin{lemma}\label{lemma291}
There exists a space  $RK_1$, which is satisfies the equation $(\ref{RK1})$,
an which is included in the commutative
diagram $ (\ref{16.21}) $. Moreover, the commutative
diagrams $(\ref{118.21})$ determines 
boundary conditions.
\end{lemma}

\subsubsection*{The complex stratification of polyhedra $ J $, $ \Sigma $, 
$ \Sigma_{ \circ} $ by means of the coordinate
system  angle - momentum}

Let us order lens spaces, which form the join, by the integers
from 1 up to $ r $ and let us denote by $ J (k_1, \dots, k_s)
\subset J $ the subjoin,  formed by a selected set of circles
(one-dimensional lens spaces) $ S^1 / \i $ with indexes $ 1 \le
k_1 <\dots <k_s \le r $, $ 0 \ge s \ge r $. The stratification
above is induced from the standard stratification of the open
faces of the standard $ r $-dimensional simplex $ \delta^{r} $
under the natural projection $ J \to \delta^{r} $. The
preimages of vertexes of a simplex are the lens spaces $ J (j)
\subset J $, $ J (j) \approx S^1 / \i $, $ 1 \le j \le r $,
generating the join.

  Define the space $ J^{[s]} $ as a subspace of $ J $, obtained
by the union of all subspaces $ J (k_1, \dots, k_s) \subset J
$.

Thus, the following stratification
\begin{eqnarray}\label{stratJ}
J^{(r)} \subset \dots \subset J^{(1)} \subset J^{(0)},
\end{eqnarray}
of the space $
J $ is well-defined. For the considered stratum a number $ r-s $ of missed
coordinates to the full set of coordinates is called the deep of
the stratum. 

Let us introduce the following denotation:
\begin{eqnarray}\label{strat[J]}
J^{[i]} = J^{(i)} \setminus J^{(i+1)}.
\end{eqnarray}
 Denote the maximum open cell of the space $\hat p^{-1} (J
(k_1, \dots, k_s)) $ by $\hat U (k_1, \dots, k_s) \subset S^{n-k}/\i$.
 This open cell is called an elementary stratum of the
depth $ (r-s) $. A point at an elementary stratum $ U (k_1,
\dots, k_s) \subset S^{n-k} / \i $ is defined by a set of
coordinates $ (\check x_{k_1}, \dots, \check x_{k_s}, \lambda) $, where
$ \check x_{k_i} \in S^1 $ is a coordinate on the 1-sphere (circle),
covering lens space with the number $ k_i $, $ \lambda = (l_{k_1},\dots, l_{k_s}) $ is a barycentric coordinate
on the corresponding $ (s-1) $-dimensional simplex of the join.
Thus if the two sets of coordinates are identified  under the
transformation of the cyclic $ \I_a $-covering by means of the
generator, which is common to the entire set of coordinates, then
these sets define the same point on $ S^{n-k} / \i $. Points on
elementary stratum $ \hat U (k_1, \dots, k_s) $ belong in the
union of simplexes with vertexes belong to the lens spaces of the join
with corresponding coordinates. Each elementary strata $ \hat U
(k_1, \dots, k_s) $ is a base space of the double covering $ U
(k_1, \dots, k_s) \to \hat U (k_1, \dots, k_s) $, which is induced
from the double covering $ \RP^{n-k} \to S^{n-k} / \i $ by the
inclusion $ \hat U (k_1, \dots, k_s) \subset S^{n-k} / \i $.

The polyhedron  $\Sigma_{\circ}$ is split into the union of
open subsets (elementary strata),
these elementary strata are defined as the connected components of 
the inverse images of  elementary strata  $(\ref{strat[J]})$.
Denote these elementary strata by
\begin{eqnarray}\label{compstrat}
K^{[r-s]}(k_1, \dots, k_s), \qquad 1 \le s
\le r. 
\end{eqnarray}

Anti-diagonal strata belong to 
 $\Sigma_{\circ}$.  Each elementary stratum, except anti-diagonal strata, is a double covering over the corresponding stratum
 of the polyhedron 
$\hat K_{\circ}$.  Denote these elementary strata by
\begin{eqnarray}\label{hatcompstrat}
\hat K^{[r-s]}(k_1, \dots, k_s), \qquad 1 \le s \le r.
\end{eqnarray}

Let us describe an elementary stratum $ K^{[r-s]} (k_1, \dots, k_s) $ by
means of the coordinate system. To simplify the notation let us
consider the case $ s = r $. Suppose that for a pair of points $
(x_1 $, $ x_2) $, defining a point on $ K^{[0]} (1, \dots, r) $,
 the following pair of points $ (\check x_1, \check x_2) $ on the covering space
  $ S^{n-k} $ is fixed, and the pair $ (\check x_1, \check x_2) $ is mapped to the pair  $ (x_1 $,
  $ x_2) $ by means of the projection of $ S^{n-k} \to \RP^{n-k} $.
  Accordingly to  the construction above, we denote by
  $ (\check x_{1, i}, \check x_{2, i}) $, $ i = 1, \dots, r $ a set of spherical coordinates of
  each point. Each such coordinate with the number $i$ defines a point on
  1-dimensional sphere
  (circle)
$ S^1_i $ with the same number $ i $,
  which covers the corresponding circle
$ J (i) \subset J $ of the join. Note that the pair of
coordinates with the common number determines the pair of points
in a common layer of the standard cyclic $ \I_a $-covering $ S^1
\to S^1 / \i $.

The collection of coordinates $ (\check x_{1, i}, \check x_{2, i})
$ are considered up to independent changes to the antipodal. In
addition, the points in the pair $ (x_1, x_2) $ does not admit a
natural order and the lift of the point in $ K $ to a pair of
points $ (\bar x_1, \bar x_2) $ on the sphere $ S^{n-k} $, is well
determined up to $8$ different possibilities. (The order of the
group $ \D $ is equal to $8$.)

An analogous construction holds for points on deeper elementary
strata $ K^{[r-s]} (k_1, \dots, k_s) $, $ 1 \le s \le r $.


\subsubsection*{The coordinate description of elementary strata of the polyhedra
$ K_{ \circ} \subset \Sigma_{\circ}$}

Let $x \in  K^{[r-s]} (k_1, \dots, k_s)$ be a point  on an
elementary stratum. Consider the sets of spherical coordinates
$\check x_{1,i}$ пїЅ $\check x_{2,i}$, $k_1 \le i \le k_s$ of the
point $ x $. For each $ i $ the following cases: a pair of $ i
$-th coordinates coincides; antipodal, the second coordinate is
obtained from first by the transformation by means of the
generator (or by the minus generator) of the cyclic cover.
Associate to an ordered pair of coordinates $ \check x_{1, k_i} $
and $ \check x_{2, k_i} $, $ 1 \le i   \le s $ the residue $ v_{k_i} = \check x_{1, k_i}(\check x_{2, k_i})^{-1} $
of a value $ +1 $, $ -1 $, $ + \i $ or $ - \i $, respectively.
It is easy to check that the collection of residues $\{v_{k_i}\}$ is changed by the following  transformation.
When  the collection of coordinates of a point is changed to the
antipodal collection, say, the collection of coordinates of the
point $ x_2 $ is changed to the antipodal collection, the set of
values of residues
 of the new pair $(\bar x_1, \bar x_2) $ on the spherical
covering is obtained from the initial set of residues by changing
of the signs:
$$\{(\check x_{1, k_i},\check x_{2, k_i})\} \mapsto \{(-\check x_{1, k_i},\check x_{2, k_i})\}, \quad \{v_{k_i}\} \mapsto \{-v_{k_i}\},$$ 
$$\{(\check x_{1, k_i},\check x_{2, k_i})\} \mapsto \{(\check x_{1, k_i},-\check x_{2, k_i})\}, \quad \{v_{k_i}\} \mapsto \{-v_{k_i}\}.$$ 
The residues of the renumbered pair of points change
by the inversion:
$$\{(\check x_{1, k_i},\check x_{2, k_i})\} \mapsto \{(\check x_{2, k_i},\check x_{1, k_i})\}, \quad \{v_{k_i}\} \mapsto \{\bar {v}_{k_i}\},$$ 
where  $v \mapsto \bar {v}$ means the complex conjugation.
  Obviously, the set of residues does not change,
if we choose another point on the same elementary stratum of the
space $ K_{\circ} $. 

Elementary strata of the space $ K (k_1, \dots, k_s) $, in
accordance with sets of residues, are divided into 3 types: $
\I_a, \I_{b \times \bb}, \I_d $. If among the set of residues are  only
residues $\{+\i,-\i\}$ (respectively, only residues $\{+1,-1\}$),
we shall speak  about the elementary stratum of the type $\I_a $
(respectively of the type $ \I_{b \times \bb} $). If among the residues are
residues from the both set $\{+\i,-\i\}$ and  $\{+1,-1\}$, we
shall speak about elementary stratum of the type $ \I_d $. It is
easy to verify that the restriction of the structure mapping $
\eta: K_{0 \circ} \to K (\D, 1) $ on an elementary stratum of
the type $ \I_a, \I_{b \times \bb}, \I_d $ is represented by the composition of
a map in the space $ K (\I_a, 1) $ (respectively in the space $ K
(\I_{b \times \bb}, 1) $ or $ K (\I_d, 1) $) with the map $i_a: K (\I_a, 1)
\to K (\D, 1) $ (respectively, with the map $ i_{b \times \bb}: K (\I_{b \times \bb}, 1)
\to K (\D, 1) $ or $ i_d: K (\I_d, 1) \to K (\D, 1) $). For
the first two types of strata the reduction of the structural
mapping (up to homotopy) is not well defined, but is defined only
up to a composition with the conjugation in the subgroups $ \I_a
$, $ \I_{b \times \bb} $.

The polyhedron $\Sigma_{\circ}$ contains the polyhedron  $K_{\circ}$ and $\Sigma_{\circ} \setminus K_{\circ}$ consists of antidiagonal elementary strata. For an arbitrary elementary antidiagonal stratum
$ K (k_1, \dots, k_s) $ the residue of the each angle coordinate is equal to  $+\i$.
A antidiagonal stratum is an elementary stratum of the type $\I_a$.
The polyhedron  $\Sigma$ is derived from  
 $\Sigma_{\circ}$ by the joining of all diagonal strata
 (on each diagonal strata the residue of an arbitrary angle coordinate is equal $+1$), which is in the boundary of the polyhedron. It is easy to verify that  $\Sigma \setminus \Sigma_{\circ}$ contains all elementary diagonal strata of the deep greater, or equal, then $1$.

Define the following open subpolyhedra
\begin{eqnarray}\label{Ia}
K_{a \circ} \subset K_{\circ} \subset \Sigma_{\circ},
\end{eqnarray}
\begin{eqnarray}\label{Ib}
K_{b \times \bb \circ} \subset K_{\circ} \subset \Sigma_{\circ},
\end{eqnarray}
\begin{eqnarray}\label{Id}
K_{d \circ} \subset K_{\circ} \subset \Sigma_{\circ}
\end{eqnarray}
as the unions of all elementary strata of the corresponding type.

The following polyhedron
\begin{eqnarray}\label{hatIb}
\hat K_{b \times \bb \circ} \subset \hat K_{\circ}
\end{eqnarray}
is defined as the base of 2-sheeted covering over the polyhedron $(\ref{Ib})$. 
The description of 
$(\ref{Ib})$ by means of the coordinates is obvious and is omitted.

\subsubsection*{Description of the structural map $\eta_{\circ}: \Sigma_{\circ} \to K(\D,1)$,
by means of the
coordinate system}

Let $x = [(x_1, x_2)]$ be a marked a point on $ K_{\circ} $,
on a maximal elementary stratum. Consider closed path $\lambda:
S^1 \to K_{\circ}$, with the initial and ending points in this
marked point, intersecting the singular strata of the depth 1 in a
general position in a finite set of points. Let $(\check x_1,
\check x_2)$ be the two spherical preimages of the point $ x $.
Define another pair $(\check x'_1,\check x'_2)$ of spherical
preimages of $ x $, which will be called  coordinates, obtained in
result of the natural transformation of the coordinates $(\check
x_1, \check x_2)$ along the path $ \lambda $.

At regular points of the path $ \lambda $ the family of pairs of
spherical preimages in the one-parameter family is changing
continuously, that uniquely identifies the  inverse images of the
end point of the path by the initial data. When crossing the path
with the strata of depth 1, the corresponding pair of spherical
coordinates with the number $ l $ is discontinuous. Since all the
other coordinates remain regular, the extension of regular
coordinates along the path at a critical moment time is uniquely
determined. For a given point $ x $ on elementary stratum of the
depth $0$ of the spaces $ K_{\circ} $ the choice of at least one
pair of spherical coordinates  is uniquely determines the choice
of spherical coordinates with the rest numbers. Consequently, the
continuation of the spherical coordinates along a path is uniquely
defined in a neighborhood of a singular point of the path.

The transformation of the ordered pair $(\check x_1, \check x_2)$
to the ordered pair $(\check x'_1,\check x'_2)$ defines an element
the group $ \D $. This element does not depend on the choice of
the path $ l $ in the class of equivalent paths, modulo homotopy
relation in the group $\pi_1(\Sigma_{\circ},x)$. Thus, the
homomorphism $\pi_1(\Sigma_{\circ},x) \to \D$ is well defined and the
induced map
\begin{eqnarray}\label{eta}
\eta_{\circ}: \Sigma_{\circ} \to K(\D,1)
\end{eqnarray}
coincides with structural mapping, which was determined earlier.
 It is easy to verify that the restriction
of the structural mapping $ \eta_{\circ} $ on the
connected components of a single elementary stratum $K_{
\circ}(1, \dots, r)$  is
homotopic to a map with the image in the subspeces $K(\I_a,1)$,
$K(\I_{b \times \bb},1)$, $K(\I_d,1)$,
which corresponds to the type and subtype elementary stratum.

\subsubsection*{Coordinate description of the canonical covering over an elementary stratum}

Consider an elementary stratum $ K^{[r-s]}(k_1, \dots, k_s) \subset
K^{(r-s)}_{\circ}$  of the depth  $ (r-s) $.
Denote by
\begin{eqnarray}\label{pi}
\pi: K^{[r-s]}(k_1, \dots, k_s) \to K(\Z/2,1)
\end{eqnarray}
 the classifying map, that is responsible for the permutation of a
pair of points around a closed path on this elementary stratum.
This mapping is called the {\it{classified}} mapping for the corresponding  2-sheeted covering.

The mapping $ \pi $  coincides with the
composition
$$ K^{[r-s]}(k_1, \dots, k_s) \stackrel{\eta}{\longrightarrow} K(\D,1) \stackrel{p}{\longrightarrow} K(\Z/2,1), $$
 where $K(\D,1)
\stackrel{p}{\longrightarrow} K(\Z/2,1)$ be the map of the
classifying spaces, which is induced by the epimorphism $ \D \to
\Z / 2 $ with kernel $ \I_c \subset \D $ .
The canonical 2-sheeted covering, which is associated with the mapping
$\pi$ let us denote  by
\begin{eqnarray}\label{bar}
\bar K^{[r-s]}(k_1, \dots, k_s) \to K^{[r-s]}(k_1, \dots, k_s).
\end{eqnarray}

With the mapping 
 $(\ref{pi})$ the following equivariant mapping is associated:
\begin{eqnarray}\label{eqvpi}
\bar{\pi}: \bar K^{[r-s]}(k_1, \dots, k_s) \to S^{\infty},
\end{eqnarray}
where the involution in the image is the standard antipodal involution.
This mapping is a 2-sheeted covering over the mapping  $(\ref{pi})$.

For an elementary strata of the type
 $\I_{b \times \bb}$ with the mapping  $(\ref{eqvpi})$ the following equivariant mapping is associated:
 \begin{eqnarray}\label{hateqvpi}
\tilde{\pi}: \tilde K^{[r-s]}(k_1, \dots, k_s) \to S^{\infty},
\end{eqnarray}
where the mapping
$\tilde K^{[r-s]}(k_1, \dots, k_s) \subset \tilde K(\E_{b \times \bb},1)$,
$(\ref{hateqvpi})$ is a 2-sheeted covering over the mapping  $(\ref{eqvpi})$.

\begin{lemma}\label{lemma32a}

The restriction of the map  $(\ref{eqvpi})$ to the canonical 2-sheeted covering over an elementary strata 
of an arbitrary type is homotopic to the following composition
\begin{eqnarray}\label{pi1a}
\bar{\pi}: \bar K^{[r-s]}(k_1, \dots, k_s) \to S^1 \subset S^{\infty}.
\end{eqnarray}
The restriction of the equivariant map  $(\ref{hateqvpi})$ to the canonical 2-sheeted covering over an elementary strata 
of the type $\E_{b \times \bb}$ is homotopic to the following composition
\begin{eqnarray}\label{hatpi1a}
\tilde{\pi}: \tilde K^{[r-s]}(k_1, \dots, k_s) \to S^1 \subset S^{\infty},
\end{eqnarray}
where $S^1 \subset S^{\infty}$ is the equivariant embedding of the standard 1-dimensional skeleton
of the classifying space.
\end{lemma}

\subsubsection*{Proof of Lemma $\ref{lemma32a}$}
Let us prove the lemma by means of explicit formulas for the mappings
$(\ref{pi1a})$ $(\ref{hatpi1a})$. An arbitrary point
$[(x_1,x_2)] \in \hat K^{[r-s,i]}(k_1, \dots, k_s)$, or $[(x_1,x_2)] \in K^{[r-s,i]}(k_1, \dots, k_s)$ 
is determined by the equivalence class of the collection of angle coordinates and the momentum coordinate. 
The structure mapping  $\eta_{\circ}$, $\hat \eta_{b \times \bb \circ}$ is determined by a transformation of angle coordinates. 
Let us define the mappings 
$(\ref{pi1a})$, $(\ref{hatpi1a})$ by the corresponding transformation of the  \emph{marked} pair of the angle coordinates. 
Below the prescribed pair of the angle coordinates for an elementary stratum of each arbitrary type is defined.   

Assume that a point
 $[(\hat{x}_1,\hat{x}_2)] \in \hat K^{[r-s]}(k_1, \dots, k_s)$ is belong to the stratum of the type $\E_{b \times \bb}$.
 Because the residue of the prescribed pair of the angle coordinates is well-defined, a non-ordered pair of the angle coordinates with the residue
 $-1$ it is convenient to denote by 
 $[(\check{x}_{1,-}, \check{x}_{2,-})]$, a pair of the angle coordinates with the residue
 $+1$ denote by  $[(\check{x}_{1,+}, \check{x}_{2,+})]$. 

The each coordinate 
 $\check{x}_{1,-}$,  $\check{x}_{2,-}$, $\check{x}_{1,+}$, $\check{x}_{2,+}$ determines the corresponding point on  $S^1$.
 It is not difficult to check, that  $\check{x}_{1,+} = \check{x}_{2,+}$, $\check{x}_{1,-} = -\check{x}_{2,-}$.
Therefore the mapping  $(\hat x_1,\hat x_2) \mapsto (\check{x}_{1,-}^{-1}\check{x}_{1,+},\check{x}_{2,-}^{-1}\check{x}_{2,+})$
transforms the points of an ordered pair into the antipodal points on  $S^1$. 
The changing of a pair of the angle coordinates to an equivalent pair, which keeps the order of the points of the pair,
does not change the equivariant mapping. The changing of the order of points in the pair transforms the equivariant mapping to the antipodal mapping.
 The constructed equivariant mapping is the required equivariant mapping $(\ref{pi1a})$ for the stratum
of the type $\E_{b \times \bb}$.

Assume a point
$[(x_1,x_2)] \in K^{[r-s,i]}(k_1, \dots, k_s)$ belongs to an elementary stratum of the type  $\I_a$ 
(including the case, when a stratum is anti-diagonal). The mapping
$(\ref{pi1a})$ is determined by a transformation of the prescribed pair of the angle coordinates with the residue  $+\i$, 
which we denote (and the same time introduce an order of the pair) as  $(\check{x}_{1,+\i},\i\check{x}_{1,+\i})$. 
The mapping $(x_1,x_2) \mapsto (\check{x}_{1,+\i}^2,-\check{x}_{1,+\i}^2) $ transforms the points of the ordered pair into an antipodal
points on  $S^1$. This mapping is the required mapping
$(\ref{pi1a})$ for the elementary stratum   of the type  $\I_a$. 

Assume a point
$(x_1,x_2) \in K^{[r-s]}(k_1, \dots, k_s)$  belongs to an elementary stratum of the type $\I_d$. The mapping
$(\ref{pi1a})$ is determined by a transformation of the prescribed pair of the angle coordinates with the residue  $+\i$, 
which we denote by $[(\check{x}_{1,+\i},\i\check{x}_{1,+\i})]$. 
 The mapping $(x_1,x_2) \mapsto (\check{x}_{1,+\i})^2,-\check{x}_{1,+\i})^2 $ transforms the points of the ordered pair into an antipodal
points on $S^1$. This mapping is the required mapping
$(\ref{pi1a})$ for the elementary stratum  of the type $\I_d$.
Let us denote that the constructed mapping 
$(\ref{pi1a})$ on each elementary stratum of the type $\I_d$ is homotopic to the constant mapping.

Lemma $\ref{lemma32a}$ is proved.
\[  \]

\subsubsection*{Prescribed coordinate system}
A coordinate system
$\Omega(\hat \alpha)$ on an elementary stratum $\alpha$ of the type  $\I_a$, $\I_d$, or $\hat \alpha$ of the type $\E_{b \times \bb}$
is determined by a fixation and an order of spherical preimages
  $(\check{x_1},\check{x_2})$ of the marked point. At the marked point a prescribed coordinate system is well-defined using the reduction of a structured subgroup of the stratum. A prescribed coordinate system is extended along an arbitrary path on an elementary stratum.
A prescribed coordinate system is well-defined, because a reduction of the structure group to the corresponding subgroup is well-defined.

\subsubsection*{Admissible pairs of elementary strata}
For an arbitrary elementary strata $\beta \subset K^{[r-s,i]}(k_1, \dots, k_s)$ of the type $\I_{a}$ 
(correspondingly, of the type $\I_d$, or of the type $\E_{b \times \bb}$) of the space
$K_{\circ}$ let us consider a smallest elementary stratum $\alpha$, $\beta \ne \alpha$ of the same type, which is contained in the boundary  of 
$\alpha \subset Cl(\beta) \subset Cl(K^{[r-s,i]}(k_1, \dots, k_s))$ (in the closure of $\alpha$). 
Analogously, for an arbitrary elementary stratum
$\hat \beta \subset \hat K^{[r-s,i]}(k_1, \dots, k_s)$ of the type $\E_{b \times \bb}$  of the space
$\hat K_{\circ}$ let us consider a smallest elementary strata $\hat \alpha$ of the same type, which is contained in the boundary of  $\hat \beta$.
Let us write $\alpha \prec \beta$, or $\hat \alpha \prec \hat \beta$.

Let us consider the prescribed coordinate system  
$\Omega(\beta)$ on the stratum $\beta \subset K(k_1, \dots, k_s)$. Restrict the prescribed coordinate system
on an elementary stratum $\alpha \subset Cl( \beta)$ inside its boundary.  Assume that the restriction of the prescribed coordinate system on $\Omega(\beta) \vert_{\alpha}$ corresponds to the prescribed coordinate system $\Omega(\alpha)$ on the smallest stratum. In this case we shall call that the pair of the strata $(\alpha,\beta)$, $\alpha \prec \beta$, is admissible. In the opposite case, we call the pair $(\alpha,\beta)$, $\alpha \prec \beta$, is non-admissible.

\subsubsection*{Convention}
Strata of the different type, say $\alpha$ of the type
 $\I_a$, and $\beta$ of the type $\I_d$, $\alpha \prec \beta$, are admissible.

\subsection*{Reduction of the structured mapping of the polyhedron  $K_{\circ}$}
Let us consider the restriction of the structured mapping on the polyhedra
 $(\ref{Ia})$, $(\ref{Id})$, $(\ref{Ib})$:
$\eta_{\circ}: K_{\circ} \to K(\D,1)$. Define the groups $\I_a \int_{\chi^{[2]}} \Z$, $\I_c \int_{\chi^{[2]}} \Z$ analogously with $\I_{b \times \bb} \int_{\chi^{[2]}} \Z$, and the group $\E_{b \times \bb} \int_{\chi^{[2]}} \Z$. The groups are equipped by the natural epimorphims  $\Phi^{[2]}: \I_a \int_{\chi^{[2]}} \Z \to \D$, $\Phi^{[2]}: \I_c \int_{\chi^{[2]}} \Z \to \D$, $p_{\I_a}: \I_a \int_{\chi^{[2]}} \Z \to \Z$,
$\Phi^{[2]}: \E_{b \times \bb} \int_{\chi^{[2]}} \Z \to \E$, $p_{\E_{b \times \bb}}: \E_{b \times \bb} \int_{\chi^{[2]}} \Z \to \Z$. 
Let us prove the following result. 

\begin{lemma}\label{redK}
Assume that the deep $r$ of the stratification 
 $(\ref{stratJ})$ is odd. 

--1. The reduction 
\begin{eqnarray}\label{etKa}
\eta_a: (K_a,UQ_{antidiag}) \to (K(\I_a \int_{\chi^{[2]}} \Z,1),K(\I_a,1) \to (K(\D,1),K(\I_a,1)).
\end{eqnarray}
of the structured mapping 
$\eta_{\circ}$ over $UQ_{antidiag} \subset K_{a} \subset K_{\circ}$ is well-defined.

--2. The reduction 
\begin{eqnarray}\label{etKd}
\eta_d: K_d \to K(\I_c \int_{\chi^{[2]}} \Z,1) \to K(\D,1).
\end{eqnarray}
of the structured mapping 
 $\eta_{\circ}$ over $K_{d} \subset K_{\circ}$ is well-defined..

--3. The reduction 
\begin{eqnarray}\label{etKb}
\hat{\eta}_{b \times \bb}: (\hat K_{b \times \bb},U\hat{Q}_{antidiag}) \to (K(\E_{b \times \bb} \int_{\chi^{[2]}} \Z,1),K(\E_{b \times \bb},1) \to (K(\E,1),K(\E_{b \times \bb},1)).
\end{eqnarray}
of the structured mapping 
$\hat{\eta}_{\circ}$ over $U\hat{Q}_{antidiag} \subset \hat{K}_{b \times \bb} \subset K_{\circ}$ is well-defined.
\end{lemma}

\subsubsection*{Proof of Lemma $\ref{redK}$}
Let us prove Statement 1. 
Define a mapping
$$ p_{a}:  K_a \to S^1, $$ 
then describe the canonical  mapping
$\eta_a: K_a \to K(\I_a \int_{\chi^{[2]}} \Z,1)$,
then prove that the mapping
$ p_{a}$ coincides with the composition 
$$ p_{\I_a} \circ \eta_a: K_a \to K(\I_a \int_{\chi^{[2]}} \Z,1) \to S^1. $$

Denote by   $K^{(2)}_{a} \subset K_a$ the polyhedron, which is defined as the union of elementary strata of the deeps  0, 1 и 2 (we use upper indexes, which are agree with the indexes of the stratification  $(\ref{stratJ})$).  
Below in the proof we omit the lower index 
$a$ to short denotations. The inclusion $K^{(2)} \subset K$ induces an isomorphism of the fundamental groups.
Therefore it is sufficient to prove that the polyhedron $K^{(2)}$ satisfies the statement of the lemma.

Define the polyhedron
$M^{(1)}$ and 2-sheeted covering $\tilde M^{(1)} \to M^{(1)}$, which is included to the commutative diagram of 2-sheeted coverings:
 \begin{eqnarray}\label{M1}
\begin{array}{ccc}
\tilde M^{(1)} & \subset & \tilde{K}^{(1)}\\
\downarrow &  &  \downarrow   \\
M^{(1)} & \subset & K^{(1)},\\
\end{array}
\end{eqnarray}
moreover, the horizontal inclusions in the diagram are homotopy equivalences. 
Consider 2-sheeted covering $\tilde K_a \to K_a$, which corresponds to the subgroup  $\I_a \subset \D$ 
in the target group of the structured homomorphism. The right 2-sheeted covering
$\tilde K^{(1)} \to K^{(1)}$ is induced from the covering $\tilde K_a \to K_a$ by the standard inclusion $K^{(1)} \subset K_a$. The left covering and the right covering in the diagram $(\ref{M1})$ are described by 
similar angles for each elementary stratum, correspondingly.
The momenta for each strata of the deep $0$ and $1$ of the polyhedron $ K^{(1)}$ belong to $(r-1)$-dimensional and $(r-2)$-dimensional simplexes correspondingly. In $M^{(1)}$ the momenta are replaced by baricenters of the corresponding simplexes and by segments, which join the corresponding baricenters. As the result, the spaces of the diagram   
$(\ref{M1})$ are equipped with mappings onto a graph $\Gamma$, which is described below. This mappings are commuted with vertical mappings of the diagram. 

Let us describe the graph $\Gamma$. Elementary strata of the deeps
$0$ and $1$ of the polyhedron $K^{(1)}$ correspond to vertexes of the graph  $\Gamma$ of the two different types.
A pair of vertexes are joined by a edge, if the corresponding stratum of the deep 1 are included to the boundary
of the corresponding (open) stratum of the deep $0$.
Denote the mapping 
$f: K^{(1)} \to \Gamma$, for which elementary strata of the deep $0$ are mapped into neighborhoods of the corresponding
vertexes. Elementary strata of the deep $1$ are mapped into corresponding vertexes.

By the construction below the mapping $f$ depends no of angles and is constructed using momenta only.
Denote by $X^{(1)}$ the polyhedron of momenta for $K^{(1)}$.  The natural forgetful projection
$g: K^{(1)} \to X^{(1)}$ is well-defined. The mapping  $f$ coincides with the composition of the projection  $g$ with the projection $h^{-1}: X^{(1)} \to \Gamma$, which is standardly defined, this projection is a homotopy equivalence.  The inclusion 
$h: \Gamma \subset X^{(1)}$, which is inverse to the projection  $h^{-1}$ is well-defined. By the inclusion $h$ 
vertexes of $\Gamma$ are mapped into barycenters of the corresponding simplexes, edges are mapped into segments, which join the baricenters of the corresponding $(r-1)$-simplexes with the baricenters  of the corresponding $(r-2)$--simplexes in the boundary. 

Define $M^{1} \subset K^{(1)}$ as the preimage of $h(\Gamma)$ by the projection $g$. The last spaces in the diagram $(\ref{M1})$ are defined analogously. 

Let us define an equivariant mapping
 $ c^{(1)}:    \tilde M^{(1)} \to S^1$, which is classified the left vertical 2-sheeted covering in the diagram $(\ref{M1})$.  Therefore using homotopy equivalences, which correspond to horizontal mappings in the diagram
$(\ref{M1})$, the equivariant mapping
$ c^{(1)}: \tilde K^{(1)} \to S^1$ is defined. Then let us check that the equivariant mapping $ c^{(1)}$
is extended to an equivariant mapping  $ c^{(2)}:   \tilde K^{(2)} \to S^1$. 

Let us fix for each maximal stratum in  $\widetilde M^{(0)}$ a system of angles  $\A_+ = (\check{x}_{1,1}, \dots, \check{x}_{1,r};\check{x}_{2,1}, \dots, \check{x}_{2,r})$ such that the corresponding system of residues  $\{v_{1}^+ = \check{x}_{1,1}\check{x}_{2,1}^{-1}, \dots v_{r}^+ = \check{x}_{1,r}\check{x}_{2,r}^{-1}\}$  
satisfy the equation: 
$$ \prod_{i=1}^r v_i^+ = \i. $$
Evidently, this is possible, because $r$ is odd. 
The coordinate system at a point  $[(x_1,x_2)] \in \widetilde M^{(0)}$ is not uniquely defined, but up to a transformation in the group $\I_a$. A system of angles $\A_-$ is defines oppositely,  the corresponding residues $\{v_{1}^-,
 \dots v_{r}^-\}$  satisfy the opposite equation:
 $$ \prod_{i=1}^r v_i^+ = -\i. $$
 Evidently, the system of angles $\A_+$ is  corresponds to a reduction of the structured mapping on the considered maximal stratum into the subgroup $\I_a$, the system  $\A_-$ corresponds to the opposite reduction homomorphism. 

Let us  fix a lift of the each residue 
$v_i^+$ at the point  $\widehat x_1 \in \widetilde K^{0}$ to the real
$\widehat{x_1}_i^+$ with respect to the projection  $proj: \R \to S^1$, $x \mapsto x \pmod {2\pi}$, $x \in \R$. 
Analogously, let us  fix a lift of the each residue  $v_i^-$ at the antipodal point $\widehat x_2 \in \widetilde K^{0}$
to the antipodal real $\widehat{x_2}_i^-=-\widehat{x_1}_i^+$. 
Denote a skew-symmetric mapping
 $\widehat c^{(0)}:  \widehat M^{(0)} \to \R$ by the formula:
 $\widehat {c}^{(0)}(\widehat x_1, \widehat x_2) = \sum_i \frac{\widehat{x}_{1,i} - \widehat{x}_{2,i}}{2}$,
this mapping contains the values  $\pm \frac{\pi}{2} \pmod {2\pi}$. On the each elementary stratum this mapping is constant. Denote by 
  $\widehat c^{(0)}:  \widetilde M^{(0)} \to S^1$ an equivariant (locally-constant) mapping, which contains imaginary values $\pm \i \in S^1$, by the composition  of the mapping 
 $\widehat c^{(0)}$ with the covering  $proj: \R \to S^1$.  

Let us extend the mapping 
  $\widehat c^{(0)}$ to an equivariant mapping $\widehat c^{(1)}:  \widetilde M^{(1)} \to S^1$ and let us prove that this mapping is the classifying mapping  for the covering $\widetilde M^{(1)} \to M^{(1)}$. 

From a technical point of view, it is convenient to reformulate the problem to construct the equivariant mapping  
$\tilde c^{(1)}$ into the problem to construct an equivariant fiberwise monomorphism of a line bundle with the fiber
  $\R^1 \cup_0 \R^1$ (a pair of real lines, which are identified at the origins), which is equipped with a fiberwise involution, into a line complex bundle, also equipped with a fiberwise involution, the both bundles are over $M^{(1)}$. 
Define the complex line bundle  $p: E(\kappa) \to M^{(1)}$ (with a fiber $\C$) as the Hermitian bundle, which is classified
by the characteristic class   $\kappa \in H^1 (M^{(1)};\Z/2)$ of the covering  $\widetilde M^{(1)} \to M^{(1)}$.
The involution  $conj: E(\kappa) \to E(\kappa)$ in the fibers of $p$ is also defined by the conjugation.

Consider the line bundle $L(\kappa) \to M^{(1)}$, which is associated with $\kappa$. In a neighborhood of a point $x \in M^{(1)}$ of the base consider a non-ordered pair  $(\psi_+, \psi_-)$ of sections, each section is a monomorphism 
$L(\kappa) \to E(\kappa)$ over the neighborhood of $x$. Additionally, the following properties are satisfied:

---the sections $\psi_+$, $\psi_-$ are transformed to each other with respect to the conjugation 
$conj: \C \to \C$  in the target fibers, 

-- the sections $\psi_+$, $\psi_-$ are non-vanished, 

-- the sections $\psi_+$ $\psi_-$ are transformed to each other along each closed path  $l \subset M^{(1)}$, for which  $\kappa(l)$ is non-vanished, i.e. the path admits not a lift into a closed path on the covering  $\widetilde M^{(1)} \to M^{(1)}$.

The family of sections
$(\psi_+,\psi_-)$ determines an equivariant morphism of a bundle with the fiber  $\R^1 \cup_0 \R^1$/ Denote this bundle by 
$G(\kappa)$, let us define this bundle. Consider the line tautological bundle $\tilde L(\kappa)$ over $\widetilde M^{(1)}$, which is associated with 2-sheeted covering
$\widetilde M^{(1)} \to  M^{(1)}$ (this covering coincides not with the canonical covering) with additional involution in the fibers. This involution is induced by the standard identification of the two lines in each fiber.
The direct image of the bundle  $\tilde L(\kappa)$ over $M^{(1)}$ is the required bundle  $G(\kappa)$.
The bundle $G(\kappa)$ is equipped with the fiberwise involution, described above.  The natural morphism
$G(\kappa) \to L(\kappa)$ of the bundle is well-defined, this morphism identifies the pair of lines in the fiber, correspondingly to the involution.  

The pair of sections
$(\psi_+,\psi_-)$ with the prescribed properties determine the equivariant morphism of the bundles
$\psi: G(\kappa) \to E(\kappa)$ over $ M^{(1)}$. 

Let us prove that the existence of an equivariant morphism
$\psi: G(\kappa) \to E(\kappa)$ implies the existence of a fiberwise monomorphism
$\psi_L: L(\kappa) \to \C \times  M^{(1)}$ of the (real) line bundle  $L(\kappa)$ into the trivial (complex) line bundle $\C \times  M^{(1)} \to  M^{(1)}$. 

The real component of the morphism
$\psi$ determines the morphism 
$\Re(\psi_L): L(\kappa) \to \R(\pm) \times  M^{(1)}$. The imaginary part $\Im(\psi_+), \Im(\psi_-)$ of the section $\psi$ represent the same equivariant morphism
$G(\kappa) \to \R(\pm\i)$  of the line equivariant bundle with the fiber  $\R^1 \cup_0 \R^1$ 
into the equivariant line bundle $\R(\pm\i)$. In a neighborhood of a point in the base of the bundles the equivariant structure of the target bundle $\R(\pm\i)$ is well defined by the antipodal transformation of fibers. The equivariant structure in the source  $G(\kappa)$ is well-defined by the transposition of the pair of the lines.
Moreover, the equivariant structure in the image and in the preimage is classified by the same characteristic class
 $\kappa$. Therefore using tensor product with  the $\Z/2$--covering  $p: \widetilde M^{(1)} \to M^{(1)}$ with characteristic class  $\kappa$ we get in the image the trivial line bundle  $\R(\pm\i) \otimes p$ over $M^{(1)}$ with equivariant structure, by the fiberwise involution, and in the preimage we get the line bundle  $(\Z/2 \tilde \times G(\kappa)) \otimes p$, which is isomorphic to the equivariant bundle
$(\Z/2 \times L(\kappa))$ (or, to the pair of bundles $L(\kappa) \cup_0 L(\kappa)$, in which the origins of fibers are identified.
The equivariant structure in the preimage   $(\Z/2 \tilde \times G(\kappa))$ is well-defined by a permutation 
of fibers of the 2 isomorphic copies of the bundle $L(\kappa)$, the equivariant structure in the (trivial) target bundle $\R(\pm\i) \otimes p$ is well-defined by the antipodal mappings.

The constructed morphism of the equivariant bundles determine a morphism of real line bundles, denote this morphism by 
$\Im(\psi_L): L(\kappa) \to \R(\pm\i)$. The Whitney sum $\Re(\psi_L) \oplus \Im(\psi_L)$ represents the required
fiberwise monomorphism $\psi: L(\kappa) \to \C \times M^{(1)}$
(where $\C = \R(\pm) \oplus \R(\pm\i)$). Additionally,  $\psi: G(\kappa) \to E(\kappa)$ and $\psi_L$ are fiberwise monomorphisms. We get the required equivariant mapping (after we take associated spherical bundles)  $\widehat c^{(1)}:  \widetilde M^{(1)} \to S^1$. 

Let us construct the pair
$(\psi_+, \psi_-)$ of conjugated sections over $M^{(1)}$,  which satisfies the required properties.
Consider the vertexes $y_j \in \Gamma$, which correspond to strata of $K^{[1]}$ of the deep $1$ and vertexes 
$x_i \in \Gamma$,  which corresponds to strata of $ K^{[0]}$ of the deep $0$. Include each vertex  $y_j$ 
in an edge $\gamma_{i,j}$ of the graph $\Gamma$ from $x_i$ to $y_j$. 
Consider a small regular neighborhood  $U_j \subset \Gamma$ of the vertex $y_j \in \Gamma$ and join this neighbourhood
to  $\gamma_{i,j}$. The last edges in the vertex $y_j$ separate from $U_j$ by a larger neighborhood $V_j$.  
In the graph $\Gamma$ a subset $\Delta$ is defined, each connected component of this subset is a subgraph, which contains a unique vertex  $x_i$, which corresponds to a stratum of the deep $0$. With each such a vertex  $x_i$ 
edges is associated, some of the edges end near the vertex  $y_j$ and contain no points from the neighborhood  $V_j$,
all the last segments are joined by the opposite end points to a neighborhood  $U_{j'}$. The complement $\Gamma \setminus \Delta$ is a collection of short segments, each segment is contained inside circular neighborhood $\cup_j (V_j \setminus U_j)$.

Extend by locally constant functions the collection of monomorphism
$(\psi^{(0)}_+$, $\psi^{(0)}_-)$ с $ M^{(0)}$ (the source space of this family is the collection of vertexes $\cup_i \{ x_i\}$, which are preimage of the mapping  $g:  M^{(1)} \to \Gamma$) to a family of monomorphisms on the preimage 
of $\Delta \subset \Gamma$. Then define the extension on the space  $ M^{(1)}$.  Consider an arbitrary segment
$[a_j,b_i] \subset \Gamma \setminus \Delta$ with the left end point $a_j$ near the vertex  $y_j$,  the right end point
is bounded by a short segment from $x_i$. Denote by  $M^{(1)}[a_j,b_i]$ the domain in
$ M^{1}$, which is defined as the inverse image of $[a_j,b_i]$ by the mapping $g$.  The following skew-symmetric mappings $\widehat c[a]$,  $\widehat c[b]$ are well-defined, the source spaces of this mappings coincide with the preimage  of the end points of the segment in $\R$. Consider the difference $\widehat c[b] - \widehat c[a] = \frac{1}{2}\sum_i \widehat x(b)_{1,i} - \widehat x(a)_{1,i} - \widehat x(b)_{2,i} + \widehat x(b)_{2,i}$. This difference is well-defined $\pmod{\pi}$ and determines a skew-symmetric function 
 $\widehat c^{(1)}[a_j,b_i] \to \R$ with prescribed boundary conditions. 

The composition of the section
 $\psi^{(0)}_+$ with the conjugation  $conj: S^1 \to S^1$ coincides with  the section 
$\psi^{(0)}_-$, which is defined using the permutation of the preimages in the canonical covering (or, equivalently, using the involution in the fibers of $G(\kappa)$). Therefore, because $\widehat c[b] - \widehat c[a]$
is a skew-symmetric, the equivariant mapping
$\widehat c^{(1)}$ with boundary conditions on  $M^{(1)}[a_j,b_i]$ is well-defined, this mapping is transformed into complex conjugated mapping after the permutation of the inverse images in the canonical covering.

The pair 
 $(\psi^{(1)}_+,\psi^{(1)}_-)$  of conjugated local monomorphisms of the bundle  $L(\kappa)$ to the bundle  $E(\kappa)$ over $M^{(1)}$ is well-defined. Equivalently, we defined an equivariant section  $\psi$ of the bundle $G(\kappa)$ to the bundle $E(\kappa)$. The required conditions for the pair  $(\psi^{(1)}_+,\psi^{(1)}_-)$ are satisfied.
The required equivariant mapping 
$\widehat c^{(1)}:  \widetilde K^{(1)} \to S^1$ is well-defined.

Let us prove that there exists an extension
 $ \widehat c^{(2)}:   \widetilde K^{(2)} \to S^1$. To proof this fact, it is sufficient to extend 
the conjugated pair $(\psi^{(1)}_+,\psi^{(1)}_-)$ to a conjugated pair
$(\psi^{(2)}_+,\psi^{(2)}_-)$ over  $K^{(2)}$. Let us investigate elements of the fundamental group $\pi_1(K^{(1)})$,
which are in the kernel of the homomorphism $\pi_1(K^{(1)}) \to \pi_1(K^{(2)})$, this homomorphism is induced by the inclusion 
$K^{(1)} \subset K^{(2)})$.
Conjugated class of the corresponding element in the kernel is well-defined by a homotopy class of a free loop
in $ M^{(1)}$ (recall, that $K^{(1)}$ and $M^{(1)}$ are homotopy equivalent) as following. 
Consider a vertex  $x_{e,e} \in \Gamma$ of the type $0$ and a pair of vertexes   $y_{1,e},y_{2,e} \in \Gamma$ of the type 1, which are in a neighborhood of  $x_{e,e}$. Let us fix $2$ arbitrary elements  $a,b \in \I_a$ (denote the unite by $e \in \I_a$) and define a sequence of vertexes  $x_{e,e},x_{a,e},x_{e,b},x_{a,b}$ of type $0$ and a sequence of vertexes $y_{1,e},y_{1,b},y_{2,e},y_{2,a}$ of type $1$ of the graph  $\Gamma$ as following. 

Consider a lift of the segment 
$[x_{e,e},y_{1,e}]$ from the marked point  $\hat x_{e,e} \in  M^{(1)}$ over $x_{e,e}$ into a point  $\hat y_{1,e}$ 
over $y_{1,e}$ by the projection $g$, the projection keeps angles. Consider a path from  
$\hat y_{1,e}$ to $\hat x_{a,e}$, which has the same momenta as the inverse path from  $\hat y_{1,e}$ to $\hat x_{e,e}$,
and the same angles, except the only angle (of number $1$), for this angle the corresponding residue is multiplicatively changed in the point $\hat y_{1,e}$ by $-1 \in \I_a$.

Let us transform the $4$-gone
$[x_{a,e},y_{1,e},x_{e,e},y_{2,e},x_{e,b}]$ by multiplicative changing residues of the coordinates  by the elements $a^{-1},b^{-1}$ correspondingly.  We get a closed $8$-gone 
$[x_{a,b},y_{2,a},x_{a,e},y_{1,e},x_{e,e},y_{2,e},x_{e,b},y_{1,b},x_{a,a}]$. Assume that this polygon is contained in $M^{(1)}$, this means that the collection of residues contains real values and contains complex values simultaneously. 

Consider an arbitrary free path on
 $ M^{(1)}$, which is defined by this polygone. It is easy to see that the image of an arbitrary such path is contractible with respect to the composition $ M^{(1)} \subset  K^{(1)} \subset  K^{(2)}$.
Evidently, the collection of all considered paths determines conjugated classes of generators in the kernel of
the homomorphism 
$\pi_1( K^{(1)}) \to \pi_1( K^{(2)})$.

Let us check that the each section
$\psi_+^{(1)}$, $\psi_-^{(1)}$
restricted to the $8$-gone path is contractible.
Assume that integer lifts
$\widetilde x_{1,i}$
$\widetilde x_{2,i}$ of the each regular angle (there exist $r-2$ regular angles) along the $8$-gone are constant. 
Then the proof is obvious. The considered assumption gives no loss of a generality, because in the general case
the sum of jumps of the each regular coordinate along the path equals to zero. The mapping
$ p_{a}:  K_a \to S^1 $ is constructed.  

Define the mapping $\eta_a$. Let us consider the homomorphism
$\pi_a \circ \eta: \pi_1(K_a) \to \Z/2$, where $\eta: \pi_1(K_a) \to \D$, and $\pi_a: \D \to \Z/2$ is the projection
with the kernel $\I_a \subset \D$.  It is easy to verify that
$p_a \cong \pi_a \pmod{2}$. Therefore $p_a$ determines a lift of the homomorphism  $\eta$ (the target group of $\eta$
is $\D \cong \I_a \int_{\chi^{[2]}} \Z/2$) to a homomorphism $\eta_a$ with the image  $\I_a \int_{\chi^{[2]}} \Z$.  The mapping $\eta_a$ classify this homomorphism (the classifying mapping and the corresponding homomorphism have the same denotations).

Let us prove that the following boundary conditions are satisfied: the restriction of the mapping
 $\eta_a$ on $UQ_{antidiag}$ taks the value in  $K(\I_a,1) \subset K(\I_a \int_{\chi^{[2]}} \Z,1)$.

Consider an arbitrary closed path  $\lambda$, which consists of a finite number of segments in a general position in $UQ_{antidiag}$. Such a path in its free homotopy class is represented (we give no a new denotation)
by a path, which intersects no with
 $K^{(2)}$ and intersects with 
$K^{(1)}$ by a finite set of points $\{l_{1,1},l_{1,2},l_{2,1},l_{2,2}, \dots , l_{k,1},l_{k,2}\}$. Moreover, for an arbitrary $j$ $(l_{j,1},l_{j,2})$ is a (short) interval where the two residue are equal to $-\i$, the last are equal to $+\i$. Along an arbitrary (long) last interval of $\lambda$ the only one corresponding coordinate has the residue  $-\i$, the last coordinates have residues $+\i$. It is sufficient, without loss of a generality, to assume that $\eta_{\circ}(\lambda)$ is a unite element in $\I_a$, because, if the boundary conditions are not satisfied, the image of the homomorphism is an infinite group (the free cyclic group), but the group  $\I_a$ itself is a finite group.

Let us investigate  $\psi_-$, this section is a constant on an arbitrary long path and 
have a jump on each short segment of the path. For each angle  an $\R$-lift and the  angle-rotation value are 
well-defined. Because the  angle-rotation value on the each short segment is equal to zero, the number of rotation of the each coordinate is divided by $\pi$. By assumption, $\eta_{\circ}(\lambda)$ is the trivial element in  $\I_a$, therefore the angle-rotation value for the each coordinate is divided  by $2\pi$.

Assume that along the path $\lambda$ angles are unchanged. Then the total sum of jumps of lifts of the each coordinates 
on short segments of the path $\lambda$ equals to zero. In this case the total changing of 
 $\psi_-$ along $\lambda$ (this changing coincides with the sum of jumps of corresponding differences of lifts of angles on the collection of the short paths) equals to zero, because at the start-point and at the end-point of the path the values of the lifts of each angle coincide, and therefore, the sum of jumps of residues of each angle is unchanged.   The assumption that each angle along the path gives no a lost of a generality, because there exists 
 an axillary closed  path, which intersects no with   $K^{(1)}$, for which the angle-rotation value of an arbitrary prescribed angle is a given $2\pi k$, and  last coordinates  are unchanged. 
 For the such  path the boundary values, evidently, are satisfied. As the result, we may change by a connected sum the given path $\lambda$ with some number of axillary paths, in this case the number of angle-rotation of the each angle along the each ling interval becomes trivial.    

Statement 1 is proved.

Statements 2 and 3 are proved analogously after the structure mappings are correspondingly changed as in the diagram $(\ref{140})$.

\subsection*{Space $Y_{a}$}

Let $\alpha$, $\beta$ be elementary strata of $\Sigma_{a}$ in the sense of the stratification $(\ref{compstrat})$.
Assume that $\alpha \prec \beta$ and define an elementary $\varepsilon$-cone of a smallest stratum 
 $\alpha$ into  the closure of $\beta$, which is defined as the open cone of a small height  $\varepsilon$, $\varepsilon <<1$,
 over the interior of the closure of the union of  all lower-dimensional $\varepsilon$-cones, which are 
inside  the closure of $\beta$. The structure of an elementary $\varepsilon$-cone corresponds to the Euclidean structure in the $r$-simplex,
given by the corresponding momenta coordinates. The elementary cone of the strata $\alpha$ in $\beta$ is denoted by   
  $Con(\alpha,\beta) \subset \beta$.

For each pair of strata (admissible, or non-admissible)
$\alpha \prec \beta$ consider the elementary  $\varepsilon$--cone  $Con(\alpha,\beta)$ and define:  
\begin{eqnarray}\label{betaodot}
\beta^{\odot} = \beta \setminus (\bigcup_{i} Con(\alpha_i,\beta;\varepsilon) \cup \bigcup_{j} Con(\alpha_j,\beta)),
 \end{eqnarray} 
where a pair $\alpha_i \prec \beta$ is non-admissible, a pair $\alpha_j \prec \beta$ is admissible. 
 
 Define the space
$Z^{\odot}\subset \Sigma_{a}$ by the formula: 
\begin{eqnarray}\label{Y}
Z^{\odot} = \bigcup_{\beta} \beta^{\odot}.
\end{eqnarray}  

Define the space 
$Y_{a}^{\odot} \subset \Sigma_{a}$ by the formula:
\begin{eqnarray}\label{Y}
Y_{a}^{\odot} = \Sigma_{a} \setminus Z^{\odot}.
\end{eqnarray} 
 
  For an arbitrary elementary  
$\varepsilon$--cone $Con(\alpha,\beta)$, which is constructed for a non-admissible pair of strata  (for an admissible pair of strata correspondingly) $\alpha \prec \beta$ define  the down- reduction, which is a subdomain $Con^{\downarrow}(\alpha,\beta)  \subset Con(\alpha,\beta)$.
 
Consider elementary strata 
$\alpha_i$ such that $\alpha_i \prec \alpha \prec \beta$, assuming that the pair 
$\alpha_i \prec \beta$ is admissible (is non-admissible, correspondingly).
Define  
\begin{eqnarray}\label{prec2}
Con^{\downarrow}(\alpha,\beta) = Con(\alpha,\beta) \setminus \bigcup_{i} Con(\alpha_i,\beta).
\end{eqnarray}

Define the space 
\begin{eqnarray}\label{ydown}
y^{\downarrow}  =
\bigcup_{\alpha,\beta} Con(\alpha,\beta) \setminus  Con^{\downarrow}(\alpha,\beta),
\end{eqnarray}
where the union is taken over all pairs of strata  
$\alpha$, $\beta$ such that the pair $\alpha \prec \beta$ is non-admissible.

For an arbitrary elementary down-reduced $\varepsilon$-cone 
 $Con^{\downarrow}(\alpha,\beta)$, which is constructed from an non-admissible pair of strata $\alpha \prec \beta$,  define down-up-reduced $\varepsilon$-cone, which is denoted by
$Con^{\downarrow\uparrow}(\alpha,\beta)  \subset Con^{\downarrow}(\alpha,\beta)$.

Consider elementary strata
$\alpha_j$ such that $\alpha \prec \alpha_j \prec \beta$ assuming that the pair 
$\alpha_j \prec \beta$ is admissible. Define
\begin{eqnarray}\label{prec22}
Con^{\downarrow\uparrow}(\alpha,\beta) = Con^{\downarrow}(\alpha,\beta) \setminus \bigcup_{j} Con(\alpha_j,\beta).
\end{eqnarray}

Define the subspace
\begin{eqnarray}\label{yup}
y^{\uparrow} = \bigcup_{\alpha, \beta} Con^{\downarrow}(\alpha,\beta) \setminus Con^{\downarrow\uparrow}(\alpha,\beta), 
\end{eqnarray}
where the union is taken over all non-admissible pairs of strata
 $\alpha \prec \beta$.

Define the space
$$yy^{\downarrow} = y^{\downarrow} \setminus (y^{\uparrow} \cap y^{\downarrow}) .$$
The following formula is satisfied:
$$yy^{\downarrow} \cup y^{\uparrow} = y^{\downarrow} \cup y^{\uparrow},$$
where the spaces in the left side of the formula are disjoint.

It is not difficult to see that the space
$yy^{\downarrow}$  is well-defined by the following formula: 
\begin{eqnarray}\label{zzdown}
yy^{\downarrow} = \bigcup_{\alpha,\beta} [Con(\alpha,\beta) \setminus Con^{\downarrow}(\alpha,\beta)] \setminus \bigcup_k Con(\alpha_k,\beta),
\end{eqnarray}
where the pair
$\alpha \prec \beta$ is non-admissible, and the following flag of strata is well-defined:
$\alpha_j \prec \alpha \prec \alpha_k \prec \beta$, where the pair 
$\alpha_k \prec \beta$ is admissible, and 
$\alpha_j$  is a stratum, such that the pair  $\alpha_j \prec \beta$ is admissible, which 
is used in the construction of the down-reduced cone
$Con^{\downarrow}(\alpha,\beta)$.

Let us clarify the formula  $(\ref{zzdown})$  as following. Let $x \in yy^{\downarrow}$ 
belongs to an elementary stratum $\beta$. Consider the maximal flag of strata 
\begin{eqnarray}\label{flagstrdown}
\alpha_{j} \prec \dots  \prec \alpha \prec  \dots  \prec \beta,
\end{eqnarray}
such that the point $x$ belongs to the cones of all this strata in the stratum $\beta$. 
Let us say, the point $x$ belongs to 
$yy^{\downarrow}$, if all strata, more deeper then  $\alpha$ (the stratum $\alpha$ itself is not included)  
are admissible with respect to $\beta$, the stratum $\alpha$ is non-admissible with respect to $\beta$,
all the last strata, less deeper then the stratum $\alpha$ are arbitrary (admissible, or non-admissible) with respect to $\beta$.
  
Define a subspace
 $$yy^{\downarrow}_{\omega \longrightarrow 1} \subset yy^{\downarrow}$$ 
 by the following formula: 
\begin{eqnarray}\label{zdown+}
yy^{\downarrow}_{\omega \longrightarrow 1} = yy^{\downarrow} \setminus \bigcup_{i} Con(\alpha_i,\beta),
\end{eqnarray}
where $yy^{\downarrow}$ is represented as in $(\ref{zzdown})$, $\alpha_j \prec \alpha_i \prec \alpha \prec \beta$,
where $\alpha_i \ne \alpha_j$ and the pair  $\alpha_i \prec \beta$ (and, also, the pair $\alpha_j \prec \beta$) is admissible.

Let us clarify the  formula $(\ref{zdown+})$ as following.  
In the maximal flag of strata $(\ref{flagstrdown})$ all cones of more deeper strata with respect to $\alpha$,
which are admissible with respect to $\beta$ are excluded, except the most deeper stratum  $\alpha_j$.

A point from the space 
$yy^{\downarrow}_{\omega \longrightarrow 1}$ belongs to a maximal flag of strata
\begin{eqnarray}\label{rem}
 \alpha_i \prec \alpha \prec \dots \prec \beta, 
\end{eqnarray} 
where the pair 
 $\alpha_i \prec \beta$ is admissible, the pair  $\alpha \prec \beta$ and the pairs of more deep strata in the sequence with  $\beta$ are non-admissible.

Let us define a space 
$yy^{\downarrow+}_{\omega \longrightarrow 1}$ and an inclusion  $$yy^{\downarrow}_{\omega \longrightarrow 1} \subset yy^{\downarrow+}_{\omega \longrightarrow 1}.$$ 
Define 
a space $yy^{\downarrow+}_{\omega \longrightarrow 1}$ as the closure of the subspace
$yy^{\downarrow}_{\omega \longrightarrow 1}$ in $yy^{\downarrow}$.

The space $yy^{\downarrow+}_{\omega \longrightarrow 1}$ is defined explicitly as following. Join all lateral boundaries of elementary strata,
of the space 
$yy^{\downarrow}_{\omega \longrightarrow 1}$ using the formula  $(\ref{rem})$. 
To the space $yy^{\downarrow}_{\omega \longrightarrow 1}$ points from 
$y^{\uparrow}$,  and $Z^{\downarrow\uparrow}$ are added  (see the formula  $(\ref{cdotZo})$ below), 
the points belong to lateral boundaries of cones
$Con(\alpha_j,\beta)$ of most deep strata in 
$(\ref{flagstrdown})$. This is sufficiently to prove separately for each elementary strata.

Let us clarify the formula $(\ref{yup})$ as following. Assume a point $x \in y^{\uparrow}$ belongs to an elementary stratum $\beta$. Consider the maximal flag of the strata: 
\begin{eqnarray}\label{flagstrup}
\alpha_i \prec \dots \prec \alpha  \prec \alpha_j \prec \dots \prec \beta,
\end{eqnarray}
such that the point $x$ belongs to the cones of all this strata in the stratum $\beta$. 
Let us say, the point $x$ belongs to 
$y^{\uparrow}$, if all strata, more deeper then  $\alpha$ (the stratum $\alpha$ itself is not included)  
are non-admissible with respect to $\beta$, the stratum $\alpha$ is admissible with respect to $\beta$,
all the last strata, less deeper then the stratum $\alpha$ are arbitrary (admissible, or non-admissible) with respect to $\beta$.

Define the subspace
 $y^{\uparrow}_{\omega \longrightarrow 1} \subset y^{\uparrow}$
 by the following formula: 
 \begin{eqnarray}\label{zup+}
y^{\uparrow}_{\omega \longrightarrow 1} = y^{\uparrow} \setminus \bigcup_{i'} Con(\alpha_{i'},\beta),
\end{eqnarray}
where $y^{\uparrow}$ is represented as in $(\ref{yup})$, $\alpha_i \prec \alpha_{i'} \prec \alpha \prec \beta$, 
where $\alpha_i \ne \alpha_{i'}$ and the pair  $\alpha_{i'} \prec \beta$ (and, also, the pair $\alpha_i \prec \beta$) is non-admissible. Let us clarify the  formula $(\ref{zup+})$  as following.  
In the maximal flag of strata $(\ref{flagstrup})$ all cones of more deep strata with respect to $\alpha$,
which are non-admissible with respect to $\beta$ are excluded, except the most deeper stratum  $\alpha_i$.

Define a space
$y^{\uparrow+}_{\omega \longrightarrow 1}$ and an inclusion $$y^{\uparrow}_{\omega \longrightarrow 1} \subset y^{\uparrow+}_{\omega \longrightarrow 1}.$$ 
Define the space $y^{\uparrow+}_{\omega \longrightarrow 1}$ as a subspace in 
$\Sigma_a$, which is obtained by joining of the lateral boundaries of elementary strata, which are used to define the space  $y^{\uparrow}_{\omega \longrightarrow 1}$ with respect to the formula  $(\ref{zdown+})$,
the boundaries of strata coincide with the lateral boundaries of corresponding cones; the natural inclusions of the subspaces
are well-defined.  It is easy to see that several extra points from the spaces $(\ref{y^ad})$)  (the formula $(\ref{cdotZo})$, see below) and $yy^{\downarrow}$ are added to the space 
$y^{\uparrow}_{\omega \longrightarrow 1}$.  

Define the subspace
\begin{eqnarray}\label{cdotZo}
Z^{\downarrow\uparrow} = \bigcup_{\alpha, \beta}  Con^{\downarrow\uparrow}(\alpha,\beta),
\end{eqnarray}
where the union is taken over all pairs of strata
 $\alpha \prec \beta$ such that the pair $\alpha \prec \beta$ is non-admissible.

Define the following sequence of subspaces
\begin{eqnarray}\label{zepochkaINT}
Z^{\odot} \subset Z_1 \subset Z'_a,
\end{eqnarray}
by the formula $Z_1 \setminus Z^{\odot} = Z^{\downarrow\uparrow}$, where $Z^{\odot}$ is given by the formula  $(\ref{Y})$,
\begin{eqnarray}\label{Z'a}
Z'_a \setminus Z_1 = y^{\uparrow}_{\omega \longrightarrow 1} \cup yy^{\downarrow}_{\omega \longrightarrow 1}.
\end{eqnarray}

The complements of the spaces in $\Sigma_a$ determine the sequence of extensions
\begin{eqnarray}\label{zepochkaEXT}
Y^{\odot} \supset Y_1 \supset Y'_a.
\end{eqnarray}

Let us consider a pair of admissible strata
 $\alpha \prec \beta$. Denote 
$$Con^{ad}(\alpha,\beta)=Con^{\downarrow}(\alpha,\beta) \setminus \bigcup_{\alpha_1} Con(\alpha_1,\beta),$$
where $\alpha \prec \alpha_1 \prec \beta$, and the pair  $\alpha_1 \prec \beta$ is non-admissible. In this formula the cone  $Con^{\downarrow}(\alpha,\beta)$ is defined by an analogous formula as the formula
$(\ref{prec2})$, but the admissible pairs of strata are changed into non-admissible, and oppositely.

It is easy to verify the following equations:
\begin{eqnarray}\label{Z}
Y'_a=y^{ad} \cup y_a^{\uparrow}  \cup yy_a^{\downarrow},
\end{eqnarray}
where
\begin{eqnarray}\label{y^ad}
y^{ad} = \bigcup_{\alpha,\beta} Con^{ad}(\alpha,\beta),
\end{eqnarray}
the union is taken over all admissible pairs of strata $\alpha \prec \beta$;
\begin{eqnarray}\label{Xdonwup}
y^{\uparrow}_a=  y^{\uparrow} \setminus y^{\uparrow}_{\omega \longrightarrow 1},
\end{eqnarray}
\begin{eqnarray}\label{Xdonwup}
yy^{\downarrow}_a=  yy^{\downarrow} \setminus yy^{\downarrow}_{\omega \longrightarrow 1}.
\end{eqnarray}

We need to extend the space
$Z'_a$ and define an extension  $Z'_a \subset Z_a$. Let us join to the subspace $Z'_a \setminus Z_a$ several points  from the spaces $y^{ad} \subset Y'_a$, $y^{\uparrow}_a \subset Y'_a$, $yy^{\downarrow}$. Replace the space $Z'_a \setminus Z_1$, which is defined by the formula $(\ref{Z'a})$, into the space 
$Z_a \setminus Z_1$, which is defined to replace the space in the right side of the formula into the space $y^{\uparrow+}_{\omega \longrightarrow 1} \cup yy^{\downarrow+}_{\omega \longrightarrow 1}$. 
To the space $Z'_a$ points from $Y'_a$ are added.
 
We get the sequence of the inclusions: 
$$Z^{\odot} \supset Z_1 \supset Z_a.$$

Let us take the complements in $\Sigma_{a}$, we get the corresponding sequence of extensions:
$$Y^{\odot} \subset Y_1 \subset Y_a,$$
The space $Y_a$ is the required space. Because of this correction, the formula
$(\ref{Z})$ for the smallest subspace is changed as following: 
\begin{eqnarray}\label{ZZ}
Y_a=Y^{ad} \cup Y_a^{\uparrow} \cup Y_a^{\downarrow},
\end{eqnarray}
$Y^{ad}\subset y^{ad}$, $Y_a^{\uparrow} \subset y_a^{\uparrow}$, $Y_a^{\downarrow} = yy^{\downarrow}_a$. 
Now in each elementary stratum $\beta$ the subspace
$Y^{ad}$ have no common points in its closure with the space $Y_a^{\uparrow}$, and also the subspaces
$Y_a^{\uparrow}$, $Y_a^{\downarrow}$ have no common points, which belong to their closures.

\subsubsection*{Structured mapping on the space $Y_a$}

Consider the mapping
$\eta:  \Sigma_{a} \to K(\D,1)$, which is defined by the formula  $(\ref{eta})$.
Consider the restriction of this mapping to the subspace
\begin{eqnarray}\label{etaZSiga}
\eta_a:  Y_a  \to K(\D,1).
\end{eqnarray}

\begin{lemma}\label{redZ}

--1. The mapping  $(\ref{etaZSiga})$ admits a reduction, which is given by a mapping 
\begin{eqnarray}\label{muZSiga}
\mu_{a}:  Y_{a} \to K(\I_a,1),
\end{eqnarray}
$i_{\I_a,\D} \circ \mu_{a} = \eta_{a\circ}$.

--2. The canonical 2-sheeted covering over the subspace  $K_{d} \subset Y_{\circ}$ is classified by a mapping into the circle.

--3. The restriction of the structured mapping to the space  $\hat K_{b \times \bb \circ}$
admits a reduction, which is given by the mapping 
\begin{eqnarray}\label{muZSigb}
\mu_{b \times \bb}:  (\hat K_{b \times \bb \circ},
RQ_{diag}) \to (K(\E_{b \times \bb} \int_{\chi} \Z,1),K(\E_{b \times \bb},1))
\end{eqnarray}
\end{lemma}

\subsubsection*{Proof of Lemma $\ref{redZ}$}

Proof of Statement 1.

By the construction the following sequence of inclusions is well-defined: 
$$Y_a \subset Y'_a \subset Y_1 \subset Y^{\odot}.$$

Consider the space 
$$Y_1 \setminus  Y_a  = y^{\uparrow+}_{\omega \longrightarrow 1} \cup yy^{\downarrow+}_{\omega \longrightarrow 1},$$
and the subspace  
$$Y_a^{\uparrow} \subset Y_a.$$

Define the space   $\Omega^{\uparrow}_+(\beta)$ as the intersection of the space  $y^{\uparrow+}_{\omega \longrightarrow 1}$ with the elementary stratum $\beta$. Analogously, define the subspace  $\Omega^{\downarrow}_+(\beta)$  in the space  $yy^{\downarrow+}_{\omega \longrightarrow 1}$ and define the subspace 
$\Omega^{\downarrow\uparrow}(\beta)$ in the space
 $Y_a^{\uparrow}$.



By construction, the following formula are satisfied:
$$y^{\uparrow+}_{\omega \longrightarrow 1} = \bigcup_{\beta} \Omega_+^{\uparrow}(\beta),$$
$$yy^{\downarrow+}_{\omega \longrightarrow 1} = \bigcup_{\beta} \Omega_+^{\downarrow}(\beta),$$
$$Y^{\uparrow}_a = \bigcup_{\beta} \Omega^{\downarrow\uparrow}(\beta).$$
We need descriptions of the spaces using intersections of cones of strata, which belong to the neighborhood of $\beta$. 

Define a subspace $A^{\uparrow}_+(\alpha_j,\alpha_k;\alpha,\beta) \subset \beta$. Assume that in the sequence
$\alpha \prec \alpha_j \prec \alpha_k \prec \beta$
the pair $\alpha \prec \beta$  is non-admissible, the pair $\alpha_j \prec \beta$ is admissible, the pair $\alpha_k \prec \beta$ is arbitrary.
The case $\alpha_j = \alpha_k$ is included. 

Consider the maximal flag  of cones between 
$\alpha$ and  $\beta$:
\begin{eqnarray}\label{flag}
\alpha  \prec \dots \prec \alpha_{j} \prec \dots \prec \alpha_{k}  \prec \dots \prec \beta.
\end{eqnarray} 

Define the subspace $A^{\uparrow}_+(\alpha_j,\alpha_k;\alpha,\beta)$ as a subspace of points in  $\beta$, which belong to the intersection of the $3$ cones $\overline{Con}(\alpha,\beta)$, $Con(\alpha_j,\beta)$, $Con(\alpha_k,\beta)$ and
are not inside another cone, where by $\overline{Con}$ is denoted the union of the corresponding cone with its lateral boundary.

Define 
 $A^{\downarrow}_+(\alpha,\alpha_j;\beta)$ as a subspace in  $\beta$, which is inside the intersection of the cones $\overline{Con}(\alpha_j,\beta)$, $Con(\alpha,\beta)$, but not inside an arbitrary (open) cone, which  is below  $\alpha$ and above $\alpha_j$, where by $\overline{Con}$ 
the union of the corresponding cone with its lateral boundary is denoted. 

Define a subspace  $A^{\downarrow\uparrow}_+(\alpha,\beta) \subset \beta$, assuming that the pair  $\alpha \prec \beta$ is non-admissible.

\begin{eqnarray}\label{Adownup}
A^{\downarrow\uparrow}_+(\alpha,\beta) = 
Con^{\downarrow}(\alpha,\beta) \setminus Con^{\downarrow\uparrow}(\alpha,\beta)
\setminus [\Omega^{\uparrow}_+(\beta) \cup \Omega^{\downarrow}_+(\beta)]
\end{eqnarray}

\begin{lemma}\label{strati}
The following formula are satisfied:

1.  $\Omega^{\uparrow}_+(\beta) = \bigcup_{j,k,\alpha} A^{\uparrow}_+(\alpha_j;\alpha,\beta)$.

2.  $\Omega^{\downarrow}_+(\beta) = \bigcup_{i,j,\alpha} A^{\downarrow}_+(\alpha_j;\alpha,\beta)$.

3.  $\Omega^{\downarrow\uparrow}_+(\beta)= \bigcup_{i} A^{\downarrow\uparrow}_+(\alpha_i,\beta)$,
where the pair
$\alpha_i \prec \beta$ is non-admissible. The spaces  $A^{\downarrow\uparrow}_+(\alpha_{i_1},\beta_1)$,
$A^{\downarrow\uparrow}_+(\alpha_{i_2},\beta_2)$ intersect each other, if and only if $\alpha_{i_2} \prec \alpha_{i_1} \prec \beta_1$ 
(the case $\alpha_{i_1} \prec \alpha_{i_2} \prec \beta_2$ is obtained by a permutation of indexes and also possible),
moreover, the intersection is inside the union of the cones 
$$ \bigcup_{j_1} Con(\alpha_{j_1},\beta_1), \qquad \bigcup_{j_2} Con(\alpha_{j_2},\beta_2),$$ 
where $\alpha_{i_1} \prec \alpha_{j_1} \prec \beta_1$, $\alpha_{i_1} \prec \alpha_{j_2} \prec \beta_2$, the pairs  $\alpha_{j_1} \prec \beta_1$, $\alpha_{j_2} \prec \beta_2$ are admissible.

\end{lemma}

\subsubsection*{Proof of Lemma  $\ref{strati}$}


Let us prove the formula 2. By definition, the following inclusion is well-defined:
\begin{eqnarray}\label{decompos}
\Omega_+^{\downarrow}(\beta) \subset
\bigcup_{i,j} \overline{Con}(\alpha_i,\beta) \cap Con(\alpha,\beta),
\end{eqnarray}
where $\alpha_i \prec \alpha \prec \beta$, and the both neighbor pairs in the sequence are non-admissible.

Divide the intersection  $Con(\alpha,\beta) \cap Con(\alpha_i,\beta)$ into 2 subdomains. The first subdomain
is defined as the complement to the lower reduction of the cone
$Con(\alpha,\beta)$ by the cone $Con(\alpha_i,\beta)$ and, additionally, points of the subdomain have 
to be inside a deeper cone
$Con(\alpha_{i_1},\beta)$, such that the pair $\alpha \prec \beta$ is admissible.
Such points belong not $\Omega_+^{\downarrow}(\beta)$  (such points belong to $Y^{\downarrow}_a$). 
The second subdomain with the lateral boundary of the cone 
$Con(\alpha,\beta)$ is included into $\Omega_+^{\downarrow}(\beta)$.  Points from the intersection  $Con(\alpha_{i_1},\beta) \cap Con(\alpha,\beta)$ belong to 
$\Omega^{\downarrow}(\beta)$, or, belong to a deeper cone of a stratum, which is admissible with respect to $\beta$.
In the last cases such a point belong no to 
 $\Omega^{\downarrow}(\beta)$. (Let us remark, that a length of the sequence of the cones is bounded by the deep of
strata in the stratification of the space 
 $\Sigma_a$.) 
Now, let us note, that the each subspace in 
$\Omega_+^{\downarrow}(\beta)$ admits additional
fraction using the cones in the maximal flag of strata, which are in the boundary of $\alpha$.
The formula 2 is proved.

Let us prove the formula 3. The inclusion
$$\Omega_+^{\downarrow\uparrow}(\beta) \subset
\bigcup_{\alpha} Con^{\downarrow}(\alpha,\beta) \setminus Con^{\downarrow\uparrow}(\alpha,\beta)$$
is satisfied by definition. 

Assume that $\beta_1=\beta_2=\beta$.
Let us investigate pairwise intersections in the right side of the formula.
Assume a point
 $x \in \Omega^{\downarrow\uparrow}(\beta)$,
belongs to $Con^{\downarrow}(\alpha,\beta)$ and, moreover, belongs to 
$Con(\alpha',\beta)$ for a more deeper strata  $\alpha'$, $\alpha' \prec \alpha \prec \beta$, where the pair
$\alpha' \prec \beta$ is non-admissible. Then $x$ belongs not to any cone  $Con(\alpha_{j},\beta)$, where
$\alpha' \prec \alpha_{j} \prec \alpha \prec \beta$, and the pair $\alpha_j \prec \beta$ is non-admissible, because
the cone  $Con(\alpha_{j},\beta)$ is eliminated from  $A^{\downarrow\uparrow}(\alpha,\beta)$. 
The case $\alpha' \prec \alpha \prec \alpha_{j} \prec \beta$ is possible, the cone  $Con(\alpha_{j},\beta)$
is eliminated from the both domains $A^{\downarrow\uparrow}(\alpha_j,\beta)$, $A^{\downarrow\uparrow}(\alpha,\beta)$.
Therefore, the intersection of the domains $A^{\downarrow\uparrow}(\alpha_j,\beta)$ and $A^{\downarrow\uparrow}(\alpha,\beta)$ is inside cones of strata, which are admissible with respect to $\beta$.

For the case $\beta_1 \ne \beta_2$ the proof is analogous. The formula $3$ is proved.

Lemma $\ref{strati}$ is proved.
\[  \]

Evidently, a reduction mapping  
\begin{eqnarray}\label{muZSigaR}
\mu_{a}:  Y^{ad} \to K(\I_a,1),
\end{eqnarray}
$i_{\I_a,\D} \circ \mu_{a} = \eta_{a\circ}$ is well-defined. Moreover, the restriction of this mapping to an arbitrary
elementary stratum 
$\beta$  coincides with the restriction of the classifying mapping. 
An arbitrary segment in  $Y^{ad}(\varepsilon)$, which joins a point on a stratum $\beta$ with a point 
on a deeper stratum  $\alpha$ preserves the reduction mapping,  because the pair 
$\alpha \prec \beta$ is admissible.  Or, the prescribed coordinate systems on  $\beta$ and $\alpha$ 
are agree. 

Let us prove the following.

\begin{lemma}\label{redukzija}

1. Prescribed coordinate systems on each pair of elementary strata $\alpha$, $\beta$ of the space 
$Y_a^{\downarrow}$, $\alpha \prec \beta$ are agree (the pair 
$\alpha \prec \beta$ is admissible). In particular, a reduction mapping of structured mapping
is well-defined:
\begin{eqnarray}\label{mudown}
\mu_a: Y_a^{\downarrow} \to K(\I_a,1).
\end{eqnarray}

2. 
Prescribed coordinate systems on each pair of elementary strata
 $\alpha$, $\beta$ of the space
$Y_a^{\uparrow}$, $\alpha \prec \beta$ 
are agree (the pair 
$\alpha \prec \beta$ is admissible). In particular, a reduction mapping of the structured mapping
is well-defined:
\begin{eqnarray}\label{muup}
\mu_a: Y_a^{\uparrow} \to K(\I_a,1).
\end{eqnarray}

3.  For an arbitrary pair of elementary strata 
$\alpha \subset  Y_a^{ad}$, 
$\beta \subset  Y^{\downarrow}_a$ the only possibility $\alpha \prec \beta$ is possible, 
where this pair is admissible.

4.  For an arbitrary pair of elementary strata 
$\alpha \subset  Y_a^{ad}$, 
$\beta \subset  Y^{\uparrow}_a$ the only possibility $\alpha \prec \beta$ is possible,
where this pair is admissible.

5.  For an arbitrary pair of elementary strata   
$\beta_1 \subset  Y_a^{\downarrow}$, 
$\beta_2 \subset  Y^{\uparrow}_a$ the following possibilities are possible:  $\beta_1 \prec \beta_2$, 
$\beta_2 \prec \beta_1$, in the each case the pair is non-admissible.

\end{lemma}

Statement 1 of Lemma  
 $\ref{redZ}$ is deduced from Lemma $\ref{redukzija}$. The reduction of the characteristic mapping 
is well-defined on the corresponding subspace by the formula 
 $(\ref{muZSigaR})$, $(\ref{mudown})$ $(\ref{muup})$.  The reduction mapping on the space $Y_a^{\uparrow}$
 is  conjugated  to the mapping  $(\ref{muup})$. Using Statements  3,4,5 of Lemma $\ref{redukzija}$, the reduction mappings on the corresponding subspaces are glued into the common reduction mapping $(\ref{muZSiga})$ on the space $Y_a$. 
 
 Statements 2 and 3 of Lemma  $\ref{redZ}$ are deduced from Statements 2 and 3  of Lemma $\ref{redK}$  correspondingly. 
 Lemma $\ref{redZ}$ is proved.

\subsubsection*{Proof of Lemma $\ref{redukzija}$}

Let us prove Statement 1, using the inclusion 
$Y^{\downarrow}_a \subset yy^{\downarrow}$.
The space $yy^{\downarrow}$ is defined by the formula $(\ref{zzdown})$. Let us prove the statement for this
greater space.

Assume $x \in \beta_1 \cap yy^{\downarrow}$. Let us consider the minimal  $\varepsilon$--cone $Con(\alpha_{min},\beta_1;\varepsilon)$, which contains the point $x$. Because the flag
$(\ref{flagstrdown})$ is well-defined, the pair 
$\alpha_{min} \prec \beta_1$ is admissible. 

Assume that $x$ belongs to the closure of $\beta_0$. Then 
$\alpha_{min}$ is the minimal stratum, for which the point  $x$ belongs to the cone $Con(\alpha_{min},\beta_0;\varepsilon_1)$ in the neighborhood of $\beta_0$.
Analogously with the previous arguments, one may concludes that the pair
 $\alpha_{min} \prec \beta_0$ is admissible. It is proved that the pair $\beta_1 \prec \beta_0$
 is admissible. Statement 1 is proved.  

Prove Statement 2.
The proof is followed from Statement 3 of Lemma   $\ref{strati}$. Let us give a straightforward proof.
Let us prove that an arbitrary elementary strata 
$\Omega^{\downarrow\uparrow}(\beta_0)$, 
$\Omega^{\downarrow\uparrow}(\beta_1)$ of the space 
$y^{\downarrow\uparrow}$,  $\beta_1 \prec \beta_0$ are such that a less deep stratum contains not a more deed
stratum in its closure, or, the pair  $\beta_1 \prec \beta_0$ is admissible. 

It is sufficient to prove that the domain 
 $A^{\downarrow\uparrow}(\alpha_1,\beta_1)$ is not contained in the closure of the domain  $A^{\downarrow\uparrow}(\alpha_0,\beta_0)$, assuming that the pair $\beta_1 \prec \beta_0$ is non-admissible.
Moreover, without loss of a generality, let us assume that  $\beta_1=\alpha_0$,
in the opposite case we may take $A^{\downarrow\uparrow}(\beta_1,\beta_0)$ instead of  $A^{\downarrow\uparrow}(\alpha_0,\beta_0)$, using the fact that the inclusion
$A^{\downarrow\uparrow}(\alpha_0,\beta_0) \cap A^{\downarrow\uparrow}(\alpha_1,\beta_1) \subset A^{\downarrow\uparrow}(\beta_1,\beta_0)$ is well-defined.

An arbitrary subspace
 $A^{\downarrow\uparrow}(\alpha_1,\beta_1=\alpha_0)$ is inside the cone $Con(\alpha_1,\beta_1)$ and intersect not with
$A^{\downarrow\uparrow}(\alpha_0,\beta_0)$, because points of the cone $Con(\alpha_1,\beta_0)$ is eliminated from
$A^{\downarrow\uparrow}(\alpha_0,\beta_0)$. A point in the lateral boundary of the cone  $Con(\alpha_1,\beta_0)$
could belongs to the closure of $A^{\downarrow\uparrow}(\alpha_0,\beta_0)$, but, an arbitrary point 
belongs to
$y^{ad}$, or, belongs to  $z^{\downarrow\uparrow}$ and is not in
$y^{\downarrow\uparrow}$. Statement 2 is proved.

Statement 3 is proved analogously to Statement 1.

Statements 4-5 are proved by explicit arguments.




 Lemma $\ref{redukzija}$ is proved.
\[  \]


Define a CW-complex
\begin{eqnarray}\label{cdotZup}
CZ_a \supset Z_a
\end{eqnarray}
as a cell closure of the space  
$(\ref{cdotZo})$. In the space  $(\ref{cdotZup})$ all open cells of the subspace
$(\ref{cdotZo})$ are replaced by its closures, and the glued mappings are extended by continuity.  

The "`resolution"' mapping 
\begin{eqnarray}\label{rezup}
R: CZ_a \to Z_a \subset \Sigma_a
\end{eqnarray}
is well-defined. The restriction of the mapping $R$ to the subspace 
 $(\ref{cdotZo})$ is the inclusion.

\begin{lemma}\label{Ycirc}

The canonical $2$-sheeted covering
 $\overline{CZ}_a \to CZ_a$ over the space
 $(\ref{cdotZup})$ is induced by an equivariant mapping
\begin{eqnarray}\label{rez}\label{F} 
 \bar F_a: \overline{CZ}_a \to \bar P, 
\end{eqnarray} 
where $\bar P$ is a $11$-dimensional CW-complex, equipped with a free involution  $T_{\bar P}$. 
 \end{lemma}

Lemma $\ref{Ycirc}$ is deduced from the following lemma. 

\begin{lemma}\label{Z_a}

1. The restriction of the canonical $2$--sheeted covering $\bar{Z}_a \to Z_a$ on the subspace 
$yy^{\uparrow+}_{\omega \longrightarrow 1} \subset Z_a$ is induced by an equivariant mapping
\begin{eqnarray}\label{rez}\label{Fup} 
 \bar F^{\uparrow}: \bar y^{\uparrow+}_{\omega \longrightarrow 1} \to \bar P^{\uparrow}, 
\end{eqnarray} 
 where $\bar P^{\uparrow}$ is a $3$-dimensional CW-complex, equipped with a free involution  $T_{\bar P^{\uparrow}}$. 
 
2. The restriction of the canonical $2$-sheeted covering  $\bar{Z}_a \to Z_a$ on the subspace 
$yy^{\downarrow+}_{\omega \longrightarrow 1} \subset Z_a$
is induced by an equivariant mapping
\begin{eqnarray}\label{rez}\label{Fdown} 
 \bar F^{\downarrow}: \bar y^{\downarrow+}_{\omega \longrightarrow 1} \to \bar P^{\downarrow}, 
\end{eqnarray} 
 where $\bar P^{\downarrow}$ is a $3$-dimensional CW-complex,  equipped with a free involution   $T_{\bar P^{\downarrow}}$. 
 
 3. The restriction of the canonical $2$-sheeted covering  $\bar{Z}_a \to Z_a$ on the subspace 
$Z_1 \subset Z_a$ is induced by an equivariant mapping
\begin{eqnarray}\label{rez}\label{F1} 
 \bar F_1: \bar Z_1 \to \bar P_1, 
\end{eqnarray} 
where  $\bar P_1$ is a $3$-dimensional CW-complex,  equipped with a free involution   $T_{\bar P_1}$. 
 
\end{lemma}

\subsubsection*{Proof of Lemma $\ref{Z_a}$}

Let us prove Statement 1.
Using Statement 1 of Lemma 
$\ref{strati}$, define the mapping $\bar F^{\uparrow}$ for each canonical covering over each elementary stratum of $\Omega_+^{\uparrow}(\beta)$, then, let us prove that the constructed mappings are glued into the
common mapping, the sorce of this mapping is the canonical covering over the space $y^{\uparrow+}_{\omega \longrightarrow 1}$. Let us decompose the space  
$\Omega_+^{\uparrow}(\beta)$ into 3 subdomain, each subdomain is itself decomposed into disjoint components,
the components have no common boundary points, which belong to the space $y^{\uparrow+}_{\omega \longrightarrow 1}$. 
By the required subdomains let us take interior and boundary points of the family of spaces  $A^{\uparrow}_+(\alpha_j,\alpha_k;\alpha,\beta)$.
Then let us mark an each domain by a label, correspondingly to an elementary stratum, and let us check that 
for neighbor strata a corresponding domain
 $\beta_1 \prec \beta$ is marked by the same label.

Consider the case $\alpha_j \ne \alpha_k$. Let us mark points  of the space $A^{\uparrow}_+(\alpha_j,\alpha_k;\alpha,\beta)$, which are inside the open cone  $Con(\alpha,\beta)$ by the label of the stratum $\alpha$. Let us say that such points belongs to the domain of type $A$.
Let us say the last points belong to the domain of type $B$.
Points of the type $B$ belong to the lateral boundary of the cone:
$\overline{Con}(\alpha,\beta)\setminus Con(\alpha,\beta)$.

Points of the domain $B$ are divided into subdomains, which are intersections of cones of the stratum (maximal)  
$\alpha_j$ with lateral boundary of the cone $Con(\alpha,\beta)$ of the stratum (minimal) $\alpha$.

Obviously, points of the domain $B$ on lateral boundaries of the cone 
$Con(\alpha,\beta)$, which are not inside an open cone $Con(\alpha_j,\beta)$ belong not to the space $y^{\uparrow+}_{\omega \longrightarrow 1}$, 
the points belong to $\beta^{\odot}$. Points of $B$ are divided, therefore, into disjoint subdomains in the maximal cone 
  $Con(\alpha_i,\beta)$ of the flag $(\ref{flag})$ (above the stratum $\alpha$).
Let us mark points of a domain in the boundary of $Con(\alpha,\beta)$ by the mark $\alpha_i$.

Let us consider the domains $A$ inside the union of all domains 
$A^{\uparrow}_+(\alpha_j;\alpha,\beta)$. Assume that two of domains of this union are marked by the same 
$\alpha$. In this case the mapping into the circle, which is defined by the reduction of the classifying mapping of the canonical covering, 
are induced by the reduction mapping for the stratum 
$\alpha$. The classifying mapping of the union of domains $A$ are induced by the mapping into the disjoint union of all strata of $Y_a$ and 
by the mapping into the circle, which is classified the canonical covering over disjoint union of elementary strata.

Assume that a subdomain $B$ in $\beta$ 
contains a point, which belongs to the closure of another domain of the same type, the both subdomains are marked  by the same stratum. 
In the case a subdomain 
 $B$ in $A^{\uparrow}_+(\alpha_j;\alpha,\beta)$ is marked by $\alpha_i$, this subdomain can contain a subdomain  $A^{\uparrow}_+(\alpha_{j_1};\alpha,\beta_1)$, which also is marked by $\alpha_i$.


The space of points in subdomains
$B$ in the union of the domains $A^{\uparrow}_+(\alpha_j;\alpha,\beta)$,
which are in  $\Omega^{\uparrow}_+(\beta)$, admits a mapping into the circle, which classify the canonical covering. 
This mapping is induced using marks by the mapping into a disjoint union of all strata of the polyhedra   
$\Sigma_{\circ}$ by composition with the mapping into the circle.

The mapping on 
the space $yy^{\uparrow+}_{\omega \longrightarrow 1}$ is well-defined by the equivariant join construction from the $3$ mappings, constructed above. Statement 1 is proved. 

The Statement 2 is analogical for the space
$\Omega^{\downarrow}_+(\beta)$.
Statement 2 is proved.

Let us prove Statement 3.
It is easy to see that an elementary stratum of the space 
$Z_1$ can contains elementary strata in its closure only, in the case the strata are admissible. 
Using the join construction over the full collection of canonical coverings, we get a mapping into $3$-dimensional CW-complex. 
Statement 3 is proved.

Lemma  $\ref{Ycirc}$ is proved.

\subsubsection*{Thickenings}

Coordinates in an elementary cone, constructed from a pair of strata 
$\alpha \prec \beta$,  complete by all last momentum coordinates (and all possible corresponding angle coordinates), which degenerate on $\beta$ and contain in the corresponding face 
of the coordinate $r$-simplex. Let us call the coordinates in $\beta$ principle coordinate and the additional coordinate correspondingly.  
By definition, an arbitrary additional momentum coordinate belongs to the interval $(0,\varepsilon_1)$, $0 < \varepsilon_1 << \varepsilon$.
Denote the thickening by  
 $C(\alpha,\beta;\varepsilon,\varepsilon_1)$.

Consider the space 
$Y_a(\varepsilon,\varepsilon_1)$. Denote the complements to $Y_a$ in $\Sigma_a$,
by $Z_a$ (a cell closure of the space  $Z_a$ is defined above in Lemma $\ref{Ycirc}$). The subspace $Z_a \subset \Sigma_a$ 
is not open and not closed.
The space 
 $Z_a$ is glued from domains inside elementary strata, this space is defined above such that an arbitrary domain inside the corresponding stratum 
contains no points of the boundary. Define a thickening of each domain inside the corresponding strata by additional coordinates. The result denote by  
 $Z^{REG}_a(\varepsilon_1)$. The space  $Z^{REG}_a(\varepsilon_1)$ is defined as the union of open subdomains, therefore this space is open.
 Moreover, the natural inclusion   $Z_a \subset Z^{REG}_a(\varepsilon_1)$ is well-defined. Consider the inclusion 
$Z_a \subset CZ_a$, which is induced from the inclusion $(\ref{cdotZup})$.

\begin{lemma}\label{inclusion}
The natural inclusion 
$Z_a \subset Z^{REG}_a(\varepsilon_1)$ admits the spitted mapping:
\begin{eqnarray}\label{splitting}
Z^{REG}_a(\varepsilon_1) \to CZ_a,
\end{eqnarray}
which is induced the isomorphism of the structural mappings.
\end{lemma}

\subsubsection*{Proof of Lemma $\ref{inclusion}$}
The space 
$CZ_a$ is a CW-complex, therefore a regular neighborhood of an inclusion  $CZ_a \subset \R^N$, $N \ge (dim(CZ_a))$, which is denoted by $UZ_a$,
admits a retraction  
$p: UZ_a \to CZ_a$. For a sufficiently small $\varepsilon_1 > 0$ the inclusion  $I: Z^{REG}_a(\varepsilon_1) \subset UZ_a$ is well-defined, 
and the restriction of this inclusion on
$Z_a$ coincides with the composition $Z_a \subset CZ_a \subset UZ_a$. The required mapping  $(\ref{splitting})$ coincides with the composition 
$p \circ I: Z_a \to 
CZ_a$. Moreover, the mapping  $(\ref{splitting})$ induces the isomorphism of the structured homomorphisms, 
because the structured homomorphism is well-defined by prescribed coordinate systems of elementary strata of the considered spaces.
Lemma 
$\ref{inclusion}$ is proved. 
\[  \]

Denote the complement 
$\Sigma_a \setminus Z^{REG}_a(\varepsilon_1)$ by  $Y_a^{REG}(\varepsilon_1)$. 
The subspace
\begin{eqnarray}\label{Yareg}
Y_a^{REG}(\varepsilon_1) \subset \Sigma_a
\end{eqnarray}
is closed. 
The inclusion
$Y_a^{REG}(\varepsilon_1) \subset Y_a$ is well-defined. In particular, one can replace in the formula $(\ref{muZSiga})$ the space $Y_a$ 
by the closed subcomplex 
$Y_a^{REG}(\varepsilon_1)$, and one can replace in the formula   $(\ref{F})$ the space  
$\overline{CZ}_a$ by its canonical 2-sheeted covering over  $Z^{REG}_a(\varepsilon_1)$, using the splitting  $(\ref{splitting})$ 
from Lemma 
 $\ref{inclusion}$. 

Analogously, define closed subspaces  
\begin{eqnarray}\label{Ybbreg}
Y_{b \times \bb}^{REG}(\varepsilon_1) \subset \Sigma_{b \times \bb}, \quad \hat Y_{b \times \bb}^{REG}(\varepsilon_1) \subset \hat \Sigma_{b \times \bb}
\end{eqnarray}
and obvious generalization of the formula 
$(\ref{splitting})$ is well-defined.

\subsubsection*{Definition of spaces  $R\Sigma_a$, $R\hat K_{b\times \bb\circ}$ in Lemma $\ref{lemma280}$}
Define the space
\begin{eqnarray}\label{RSig}
R\Sigma_a  \subset Y_{\circ},
\end{eqnarray}
as a thin regular neighborhood of the subspace  
 $(\ref{Yareg})$. 
 
Define the space 
\begin{eqnarray}\label{hatRSig}
RK_{b\times \bb\circ} \subset  Y_{\circ}
\end{eqnarray}
 by strata of the type  $\I_{b \times \bb}$ (see $(\ref{Ib})$, $(\ref{hatIb})$).
The space 
 $(\ref{hatRSig})$ is the standard 2-sheeted covering space, the base of this space denote by  $R\hat K_{b \times \bb}$.

All maps, which are included in the diagram 
$(\ref{16.23})$, in particular, the mapping  $pr$, $p \hat r$, are standardly defined.


\subsubsection*{Resolution maps $\phi_a: R\Sigma_a \to K(\I_a,1)$
$\hat \phi_{b \times \bb}: \hat R\Sigma_{b \times \bb\circ} \to 
K(\E_{b \times \bb} \int_{\chi^{[2]}} \Z,1)$. 
Proof of Lemma  $\ref{lemma280}$}

Let us consider the restriction 
\begin{eqnarray}\label{etaaRK}
\eta_{\circ} \vert_{Y_a}: Y_a \to K(\D,1),
\end{eqnarray} 
(recall, a homotopy equivalence  $R\Sigma_a\cong Y_a$ is well defined)
of the structure mappings on the subpolyhedron
 $(\ref{RSig})$. By Statement 1 of Lemma 
$\ref{redZ}$ the reduction of the mapping  $(\ref{etaaRK})$ is well defined:
\begin{eqnarray}\label{etaaaRK}
\phi_{Y_a}: Y_{a} \to K(\I_a,1),
\end{eqnarray} 
i.e the following equality is satisfied: $\eta_{\circ} \vert_{Y_a}= i_{\I_a} \circ \phi_{Y_a}$.
The resolution mapping
$\hat \phi_{b \times \bb}: \hat R\Sigma_{b \times \bb\circ} \to 
K(\E_{b \times \bb} \int_{\chi^{[2]}} \Z,1)$ is constructed in Statement 3 of Lemma  $\ref{redZ}$
Lemma $\ref{lemma280}$ is proved.
\[  \]

\subsubsection*{Proof of Statement 2 of Lemma $\ref{redZ}$} 

First of all, let us define the analog
$Y_{d}^{REG}(\varepsilon_3)$ of the closed space  $(\ref{Yareg})$ and of its complement  $Z^{REG}_d(\varepsilon_1)$ 
in the polyhedra $K_{d}$. Then, analogously to  Lemma  $\ref{Ycirc}$,  prove that the canonical covering over the space  $Y_{d}^{REG}(\varepsilon_1)$ 
is induced by a mapping into a $3$-polyhedron. The restriction of the structure mapping to the complement of the space
$(\ref{Yareg})$ in the polyhedra
$K_{d}$ admits a reduction to a mapping into $K(\I_c,1)$. Therefore, the canonical covering over the space $Z^{REG}_d(\varepsilon_1)$ 
is trivial. The canonical covering over the space $K_d$ is induced by a mapping into 4-polyhedron.
Lemma 
$\ref{redZ}$ is proved.

\subsubsection*{Proof of Lemma
 $\ref{osnlemma1}$}
 
 Using the dimensional restriction 
$(\ref{dimdimdim})$, let us consider an axillary mapping  $(\ref{c})$ and the composition 
$F \circ c: \RP^{n-k} \to J \subset \R^{n-17}$.  Consider formal (equivariant) mappings 
$(F \circ c)^{(2)}$, $c^{(2)}$, which are defined as formal extensions of the corresponding mappings.
Polyhedra of formal self-intersections of the mappings
$(F \circ c)^{(2)}$ and
$c^{(2)}$ coincides.  The required equivariant deformation of the formal mapping $(F \circ c)^{(2)}$ into the formal mapping $d^{(2)}$ define as vertical lift of the canonical covering over the subspace  $\Sigma_{\circ}$ by the following $2$ steps.

--1. Points of the polyhedron   $K_{d} \subset \Sigma_{\circ}$  of a formal self-intersection of the mapping  $(F \circ  c)^{(2)}$ are eliminated by a vertical lift in the codimension $2$. This is possible by Proposition 2 of Lemma  
 $\ref{redZ}$.

--2. Points of polyhedra  $Z^{REG}_a(\varepsilon_1)$ are eliminated by a vertical lift in the codimension $15$, this is possible by Lemma  $\ref{Ycirc}$, because the mapping  $(\ref{splitting})$ induces an equivalence of the structured mappings by Lemma  
 $\ref{inclusion}$.

The formal (equivariant) mapping 
  $d^{(2)}$ is well-defined.
 
Let us check the following two conditions, which are formulated in Lemma 
$\ref{lemma20}$.  The first condition for the component $\N_a$ is satisfied, namely, the mapping $\mu_a$, restricted to the component $\N_a$, 
determines a reduction of the structured mapping  $\eta_{\circ}$.

Let us prove two conditions in the statement of $\ref{lemma20}$. Condition 1 is, obviously, well proved,
namely, the restriction of the mapping  $\eta_{\circ}$
to the 
marked component $N_a$ admits a cyclic reduction, given by  $\mu_a$.

Let us prove Condition 2 in  $\ref{lemma20}$, which is formulated for the component
$N_{b \times \bb\circ}$.  For the convenience let us write-down this condition: 
\begin{eqnarray}\label{usll}
0=(p_{\I_c,\I_d} \circ \bar \eta)_{\ast}([\bar N_{b \times \bb}]) \in  H_{n-2k}(K(\I_d,1);\Z/2).
\end{eqnarray}

Assume that the polyhedron 
 $N_{b \times \bb\circ}$ is closed (let us remain that in this case the lower index $\circ$ in omitted) 
 and the mapping  $\eta$ admits a reduction
\begin{eqnarray}\label{redetabb}
\eta_{b \times \bb}: N_{b \times \bb} \to K(\I_{b \times \bb},1).
\end{eqnarray}
In this case the formula 
 $(\ref{usll})$ is satisfied, because the composition  
$$\bar \eta_{b \times \bb}: \bar N_{b \times \bb} \to K(\I_d,1)$$ 
is the composition of a mapping $N_{b \times \bb} \to K(\I_d,1)$
with the standard 2-sheeted covering
$$ \bar N_{b \times \bb} \to N_{b \times \bb} \to K(\I_{b \times \bb},1) \to K(\I_d,1),$$ 
where the mapping $K(\I_{b \times \bb},1) \to K(\I_d,1)$ is induced by the homomorphism 
$\I_{b \times \bb} \to \I_d$ with the kernel $\I_b \subset \I_{b \times \bb}$.

Assume that the polyhedron
 $N_{b \times \bb\circ}$ is not closed, and the mapping 
$\eta_{\circ}$ admits a reduction  $(\ref{redetabb})$ with the prescribed boundary conditions. 
The formula 
 $(\ref{usll})$ is rewritten as follows:
\begin{eqnarray}\label{usll2}
0=(p_{\I_c,\I_d} \circ \bar \eta_{b \times \bb \circ,\circ})_{\ast}([C\bar N_{b \times \bb\circ}]) \in  H_{n-2k}(K(\I_d,1);\Z/2).
\end{eqnarray}
The difference between the formulas 
 $(\ref{usll2})$ and $(\ref{usll})$ is following: if the polyhedron 
$N_{b \times \bb \circ}$ is non-closed, then the polyhedron  $\bar N_{b \times \bb\circ}$
is also non-closed. Therefore the polyhedron
$\bar N_{b \times \bb\circ}$ have to be compactified into a closed by a gluing of the cone of the canonical 2-sheeted cover 
$\bar N_{b \times \bb \circ} \to  N_{b \times \bb \circ}$ over the boundary. 
The result is a closed polyhedron, which is denoted in the formula
$(\ref{usll2})$ by $пїЅ\bar N_{b \times \bb \circ}$. The polyhedron
$пїЅ\bar N_{b \times \bb \circ}$ is the covering space of the 2-sheeted covering 
$пїЅ\bar N_{b \times \bb \circ} \to C N_{b \times \bb \circ}$, which corresponds to the subgroup 
$\I_{\bb} \subset \I_{b \times \bb}$ of the index 2. Therefore, as in the previous case, the cycle
$p_{\I_c,\I_d} \circ \bar \eta_{b \times \bb \circ}: C\bar N_{b \times \bb \circ} \to K(\I_d,1)$ is a boundary.

Let us consider a general case: the polyhedron  $N_{b \times \bb \circ}$ is non-closed and the mapping
$\eta_{\circ}$ admits a reduction 
$$\eta_{b \times \bb\circ}: N_{b \times \bb\circ} \to K(\I_{b \times \bb} \int_{\chi^{[2]}} \Z,1)$$
with prescribed boundary conditions.

By the assumption the following mapping
$$\hat \eta_{b \times \bb\circ}: \hat N_{b \times \bb\circ} \to K(\E_{b \times \bb} \int_{\chi^{[2]}} \Z,1)$$
is well-defined. Consider the 2-sheeted covering over the structure mapping, which we denote by 
$$\tilde \eta_{b \times \bb\circ}: C\widetilde{N}_{b \times \bb \circ} \to K(\E_{d} \times \Z,1).$$ 

Let us recall, that respectively to the diagram 
$(\ref{140})$, the 2-sheeted covering mapping $\tilde \eta_{b \times \bb \circ}$ over $\eta_{b \times \bb \circ}$ 
is totally defined by the subgroup of the index 2: 
\begin{eqnarray}\label{Ebb}
\E_{d} \times \Z \subset \E_{b \times \bb} \int_{\hat \chi^{[2]}} \Z.
\end{eqnarray}

The formula
$(\ref{usll2})$ is equivalent to the following condition: the homology class
\begin{eqnarray}\label{usll4}
(p_{\E_{d} \times \Z,\E_d} \circ \tilde \eta_{b \times \bb\circ})_{\ast}([C\widetilde{N}_{b \times \bb \circ}]) \in  H_{n-2k}(K(\E_d,1);\Z)
\end{eqnarray}
is even.

By the representation 
$\E_{b \times \bb} \int_{\hat \chi^{[2]}} \Z \to \Z/2^{[3]}$ the universal  4-bundle over $K(\E_{b \times \bb} \int_{\hat \chi^{[2]}} \Z,1)$
is well-defined, denote this bundle by 
$\hat \tau_{b \times \bb}$. The bundle 
\begin{eqnarray}\label{tauN}
\hat \eta_{b \times \bb\circ}^{\ast}(\hat \tau_{b \times \bb})
\end{eqnarray}
over $\hat N_{b \times \bb\circ}$ is well-defined.  

Denote by
\begin{eqnarray}\label{NNbb}
\widehat {NN}_{\circ} \subset \hat N_{b \times \bb \circ}
\end{eqnarray}
the 3-dimensional subpolyhedron, generally speaking, with boundary, as a homology Euler class
of the Whitney sum of 
$\frac{n-2k-3}{4}$ copies of the bundle $(\ref{tauN})$. The condition
 $(\ref{usll4})$ is equivalent to the following: the homology class 
\begin{eqnarray}\label{usll5}
(p_{\E_{d} \times \Z,\E_d} \circ \tilde \eta_{b \times \bb\circ})_{\ast}([C\widetilde {NN}_{\circ}]) \in  H_{3}(K(\E_d,1);\Z)
\end{eqnarray}
is even.

Consider the mapping
$\widehat {NN}_{\circ} \to K(\E_{b \times \bb} \int_{\chi^{[2]}} \Z,1) \to K(\Z,1)$.
Without loss of the generality, the inverse image by this mapping of the marked point of
$S^1 = K(\Z,1)$ is a closed 2-dimensional subpolyhedron, denoted by  
\begin{eqnarray}\label{LL}
\widehat {LL} \subset \widehat {NN}_{\circ}.
\end{eqnarray}
This polyhedron is ${\rm{PL}}$--homeomorphic to an oriented surface, which is equipped with a mapping
\begin{eqnarray}\label{fLL}
\hat f: \widehat {LL} \longrightarrow K(\E_{b \times \bb},1).
\end{eqnarray}
Let us use the following isomorphism: 
$H_2(K(\E_{b \times \bb},1);\Z)=\Z/2$. 

Let us prove that there exists a closed oriented 3-manifold
$\widehat {NN}$, its submanifold as in the formula  $(\ref{LL})$ and a mapping  
\begin{eqnarray}\label{hatF}
\hat F: \widehat {NN} \to K(\E_{b \times \bb} \int_{\chi^{[2]}} \Z,1), 
\end{eqnarray}
for which the following two conditions are satisfied:

--1. The image of the fundamental class by the mapping 
$(\ref{fLL})$ determines the generator of the group  $H_2(K(\E_{b \times \bb},1);\Z)$.

--2. The image of the fundamental class by the mapping 
$$\tilde F: \widetilde {NN} \to K(\E_d \times \Z,1) \to K(\E_d,1) = K(\Z/4,1) $$
is an even (or the trivial) element in the group 
$H_3(\E_d;\Z)$. 

Let us consider 2-torus $\widehat {LL}$, which is the the 2-skeleton of the standard cell decomposition of the space
$(\RP^{\infty} \times \RP^{\infty})/T_{\i} \supset (\RP^1 \times \RP^1)/T_{\i} = \widetilde {LL}$, where 
$T_{\i}:\RP^{\infty} \times \RP^{\infty} \to \RP^{\infty} \times \RP^{\infty}$ is the diagonal involution, 
which is defined by the standard involution $\i: \RP^{\infty} \to \RP^{\infty}$.  We may visualized the space $K(\E_d,1)$  as the space
$(\RP^{\infty} \times \RP^{\infty})/T_{\i} \setminus diag(\RP^{\infty})$.
By this construction the involution
$\hat \chi^{[2]}: K(\E_d,1) \to K(\E_d,1)$, which corresponds to the automorphism  $(\ref{hatchiE})$
is defined by the formula:  $x \times y \mapsto y \times x$. 

Define the (orientation preserving) involution $\hat{\chi}: \widehat {LL} \to \widehat {LL}$, which permutes the factors and reverses the diagonal. 
Define the mapping $\hat f: \widehat {LL} \to  K(\E_{b \times \bb},1)$ $(\ref{fLL})$, which transforms the diagonal 
generator 
$\i \in H_1(\widehat {LL};\Z)$ to the element $ab \in \H_{b \times \bb}$ (this element is represented by the sum of the diagonal loop with the generic loop of the first factor). Obviously, the mapping  $\hat f$ 
commutes up to homotopies with the involutions $\hat \chi$, $\hat \chi^{[2]}$ in the source and target spaces of the mapping
$\hat f$. Let us call the considered property Gluing Condition.

 Let us define the manifold
$\widehat {NN}$ as an oriented 3-manifold by the cylinder of the involution $\hat \chi: \widehat {LL} \to 
\widehat {LL}$.
The mapping  $(\ref{hatF})$ is well-defined by a fibered family over $S^1$ of mappings of
2-tori in the space
$K(\E_{b \times \bb},1)$ (the source and the target space of $(\ref{hatF})$ is the total spaces of fibrations
over $S^1$). By Gluing Condition the mapping $(\ref{hatF})$ is well-defined. This mapping satisfies Condition 1.

Let us check Condition 2. Consider the following composition:
\begin{eqnarray}\label{NNZ2}
p_{\E_d,\Z/2} \circ \tilde F: \widetilde {NN} \to K(\E_{d} \times \Z,1) \to K(\E_d,1) \to K(\Z/2,1),
\end{eqnarray}
where the mapping
$p_{\E_d,\Z/2}: K(\E_d,1) \to K(\Z/2,1)$ is induced by the epimorphism   $\E_d \to \Z/2$ with the kernel $\I_d \subset \E_d$. It is well-known, that the cellular mapping
 $p_{\E_d,\Z/2}$ transforms the standard 3-skeleton $S^3/\i \subset K(\E_d,1)$ into the standard 3-skeleton 
$\RP^3 \subset K(\Z/2,1)$ with degree 2.

Assuming Condition 2 is not satisfied and the mapping 
 $(\ref{hatF})$ determines the generic homology class, then the mapping
$(\ref{NNZ2})$ is not homotopic to zero. Assume that the mapping $(\ref{NNZ2})$ is cellular. 
Then the image of this mapping coincides with the standard 3-skeleton  $\RP^3 \subset K(\Z/2,1)$
and the degree of the mapping $(\ref{NNZ2})$ is equal to $2$ modulo 4. 
 
The mapping $(\ref{NNZ2})$ is a 2-sheeted covering over the mapping 
\begin{eqnarray}\label{hatNNZ2}
\widehat {NN} \to K(\E_{b \times \bb} \int _{\chi^{[2]}} \Z,1) \to K(\Z/2 \times \Z,1) \to K(\Z/2,1).
\end{eqnarray}
By the construction, the mapping 
 $(\ref{hatNNZ2})$ is homotopic to a mapping into the standard 2-skeleton $\RP^2 \subset K(\Z/2,1)$. 
This implies that image of the fundamental class by the mapping  
$(\ref{hatNNZ2})$, and by the mapping 
 $(\ref{NNZ2})$ is the trivial homology class. This prove that the degree of the mapping
 $(\ref{NNZ2})$ is equal to $0$ modulo 4. The mapping $\hat F$ satisfies Condition 2.

To prove Condition 
$(\ref{usll5})$ we may assume that the image of the fundamental class by the mapping $(\ref{fLL})$
is the trivial homology class. Therefore it is sufficiently to prove Condition $(\ref{usll5})$, assuming, that the surface $\hat LL$ is empty. In this case the mapping  $\hat \eta_{b \times \bb\circ}$ admits a reduction into the
subspace
 $K(\E_{b \times \bb},1) \subset K(\E_{b \times \bb}\int_{\chi^{[2]}} \Z,1)$.
Condition $(\ref{usll5})$ is reformulated analogously to Condition  $(\ref{usll2})$, which was proved above.
Condition 2 from [Lemma 27, A1] is proved. Lemma 
 $\ref{osnlemma1}$A  is proved.

\subsubsection*{A sketch of the proof of Lemma  $\ref{osnlemma2}$}

The proof is analogous to the proof of the main result of the paper
 $\cite{A1}$. Let us consider an axillary mapping
$p_1: S^{n-2k+n_{\sigma-1}+1}/\i \to J_1$,
given by the
formula  $(\ref{p_1})$, define by $C_{p_1}$ the cylinder of this
mapping. The projections  $\pi_I: C_{p_1} \to [0,1]$, $\pi_J:
C_{p_1} \to J_1$ are well defined, denote the Cartesian product of
this mappings by $F_1: C_{p_1} \to J_1 \times [0,1]$.

\begin{eqnarray}\label{A}
\begin{array}{ccc}
C_{\tilde p_1} & \longrightarrow & C_{p_1} \\
\downarrow \tilde \pi_I & & \swarrow \pi_I  \\
I & & \\
\end{array}
\end{eqnarray}

\begin{eqnarray}\label{B}
\begin{array}{ccc}
C_{\tilde p_1} & \longrightarrow & C_{p_1} \\
\downarrow \tilde \pi_J & & \swarrow \pi_J  \\
J_1 & & \\
\end{array}
\end{eqnarray}

Consider the inclusion  $I_J: J_1 \times [0,1] \subset \R^n \times
[0,1]$ and define the mapping  $I_J \circ \tilde F_1: C_{\tilde
p_1} \to \R^n \times [0,1]$, $I_J \circ F_1: C_{p_1} \to \R^n
\times [0,1]$. Consider the mapping $\tilde f_1: C_{\hat p_1} \to
\R^n \times [0,1]$ which was defined by a small generic alteration
of the mapping  $I_J \circ \tilde F_1$. The mapping $\tilde f_1$
will be taken to be coincided on the bottom of the cylinder $J_1
\subset C_{\tilde p_1}$ with the embedding  $I_J: J_1 \subset \R^n
\times \{0\}$. Moreover, the composition $p_{[0,1]} \circ \tilde
f_1 : C_{\tilde p_1} \to [0,1]$ to be coincided with  $\tilde
p_{I}$, where $p_{I}: \R^n \times [0,1] \to [0,1]$ is the
projection on the second factor. The mapping $f: C_{p_1} \to \R^n
\times [0,1]$ is also defined such that $\tilde f_1 = f_1 \circ
r_1$.

Denote by $\bar Q_1 \subset C_{p_1}$ the polyhedron of
self-intersection points of the mapping  $f_1$, defined as the
closure of the corresponded spaces by the formula:
$$ \bar Q_1 = Cl\{ x \in C_{p_1} : \exists y \in C_{p_1}, x \ne y,
f(x) = f(y) \}. $$ Because $n-4k=n_{\sigma}$, $\dim(\bar {\tilde
Q}_1)=n_{\sigma+1}+1$.

Denote by $\bar {\tilde Q}_1 \subset C_{\tilde p_1}$ the
polyhedron of self-intersection points of the mapping  $\tilde
f_1$, this polyhedron is defined as the closure of the
corresponded subspaces by the formula
$$ \bar {\tilde Q}_1 = Cl\{ x \in C_{\tilde p_1} : \exists y \in C_{\tilde p_1}, x \ne y, \tilde
f_1(x) = \tilde f_1(y) \}. $$ Because $n-4k=n_{\sigma}$, we get
$\dim(\bar Q_1)=n_{\sigma}+1$.

Consider the stratification $J_1^{[2]} \subset J_1^{[1]} \subset
J_1$ of the join. Denote by $\bar Q_{J_1}$ the intersection $\bar
Q \cap J_1$. Denote by  $\bar {\tilde Q}_{J_1}$ the intersection
$\bar {\tilde Q}_1 \cap J_1$. The polyhedron  $\bar Q_{J_1}$ has
the codimension $n_{\sigma+1}$. Because the codimension of
$J^{[2]}_1 \subset J_1$ is equal to $n_{\sigma+1}+1$, the
polyhedron  $\bar Q_{J_1} \subset J_1$ is outside a regular
neighborhood of the stratum $J^{[2]}_1$. The polyhedron  $\bar
{\tilde Q}_{J_1}$ has the codimension $n_{\sigma}$. Because the
codimension of  $J^{[1]}_1 \subset J_1$ is equal to
$n_{\sigma}+1$, the polyhedron $\bar {\tilde Q}_{J_1} \subset J_1$
is outside a regular neighborhood of the stratum  $J^{[1]}_1$.
Define the polyhedron $\bar {\tilde Q}_{J_1}(\varepsilon)$ as the
set of points from   $\bar {\tilde Q}_{J_1}$ which are mapped with
respect to the projection $\tilde \pi_I$ into a small positive
$\varepsilon \in I$.

Define the involution  $T_{\bar {\tilde Q}}: \bar {\tilde Q} \to
\bar {\tilde Q}$ which permutes points of self-intersection on the
canonical covering. The involution   $T_{\bar {\tilde Q}}$ keeps
the values of the mapping $\tilde \pi_I$. The polyhedron $\bar
{\tilde Q}_{J_1}(\varepsilon)$ is invariant with respect to the
involution   $T_{\bar {\tilde Q}}$. Denote by $T_{\bar {\tilde
Q}}(\varepsilon)$ the restriction of the considered involution on
the polyhedron $\bar {\tilde Q}_{J_1}(\varepsilon)$, this
restriction is a free involution.

Define the mapping $d_1: S^{n-2k}/\i \to \R^n \times
\{\varepsilon\} = \R^n$ as the restriction of the mapping  $\tilde
f_1$ on $S^{n-2k}/\i \times \{\varepsilon\}$. A quotient $\bar
{\tilde Q}_{J_1}(\varepsilon)/T_{\bar {\tilde Q}}(\varepsilon)$ is
a polyhedron of self-intersection points of the mapping  $d_1$.
Consider the polyhedron of self-intersection of the mapping
 $d_1$ and  its subpolyhedron  $N_1$.
By the construction, if the positive parameter
$\varepsilon$ is small enough,  the structured mapping $\zeta: N_1
\to K(\H,1)$ admits a reduction to a mapping into the subspace
 $K(\Q,1) \cup
K(\H_b,1) \subset K(\H,1)$, the considered reduction is well
defined as the composition of the mapping $t_1: N_1 \to RK_1$ with
the mapping $\phi_1: RK_1 \to K(\Q,1) \cup K(\H_b,1)$ (see the diagram $(\ref{RK1})$). 

Let us prove that the mapping
 $t_1$ satisfies the boundary conditions from diagram  $(\ref{118.21})$ in Lemma $\ref{lemma291}$. 
For $\ell \ge 8$ the number  $r_1$ of the factors of the join $J_1$, which is calculated by the formula 
$(\ref{r_1})$, is greater then  $n_{\sigma}$. Because $\dim(N_1)=n_{\sigma-1}-1$, the boundary of the polyhedron
$N_1$ contains no strata of a deep greater then  $\frac{r_1-1}{2}$. Therefore the coordinate system 
in each component  $N_1$ of the type $\E_b$ is agree with boundary conditions. Lemma  $\ref{osnlemma2}$ is proved.

\section{Proof of Theorem  $\ref{th9}$}

Let us take a positive integer $k$ under the condition  $n-4k =
n_{\sigma}$, $k \ge 8$, this is possible if  $n \ge 127$. Let the
triple $[(g: N^{n-2k} \looparrowright \R^n,\eta,\Psi)]$ represent
the given element in the cobordism group $Imm^{\D}(n-2k,2k)$. Let
us denote by $L^{n-4k}_a$ the self-intersection manifold of the
immersion $g_a$, which is the restriction of $g$ on the marked component $N^{n-2k}_a \subset N^{n-2k}$. 
Let us consider a skew-framed immersion
$(f,\kappa,\Xi)$, such that  $\delta^{sf}_k([(f,\kappa,\Xi)]) =
[(g: N^{n-2k} \looparrowright \R^n,\eta,\Psi)]$. By Proposition
$\ref{prop23}$ we may assume that the triple  $[(g,\eta,\Psi)]$
admits a quaternionic structure in the sense of Definition
$\ref{def16}$.

In the first step let us assume that the classifying map $\eta$ of
the $\D$--framed immersion is cyclic in the sense of Definition
 $\ref{defcycl}$. This means that for the marked component the following equation is satisfied:
$$ \eta = i_{a} \circ \mu_a, $$
where $\mu_a: N^{n-2k}_a \to K(\I_a,1)$, $N^{n-2k}_a = N^{n-2k}$ and $i_{a}: K(\I_a,1) \to
K(\D,1)$ is the natural map  induced by the inclusion of the
subgroup.

Let us also assume that the classifying map $\zeta$ of the
$\Z/2^{[3]}$--framed immersion $\delta^{\D}_k(g, \eta,
\Psi)=(h,\zeta,\Lambda)$ is quaternionic in the sense of Definition
 $\ref{defQ}$. This means that $L^{n-4k}=L^{n-4k}_{\Q}$ and the following equation is satisfied:
\begin{eqnarray}\label{1.504}
\zeta = i_{\Q,\Z/2^{[3]}} \circ \lambda_{\Q},
\end{eqnarray}
where $\lambda_{\Q}: L^{n-4k}_{\Q} \to K(\Q,1)$, $L^{n-4k}_{\Q}=L^{n-4k}$, and $i_{\Q,\Z/2^{[3]}}:
K(\Q,1) \to K(\Z/2^{[3]},1)$ is the natural map, induced by the
inclusion of the subgroup (see Example $\ref{def17}$). Let us
prove the theorem in this case.

Let us consider the classifying mapping  $\eta: N^{n-2k} \to
K(\D,1)$. Let us denote by $\tilde N^{n-2k-2} \subset N^{n-2k}$
the submanifold, representing the Euler class of the vector bundle
$\eta^{\ast}(\psi_{\D})$, where by $\psi_{\D}$ is denoted the universal
2-dimensional vector bundle over the classifying space
 $K(\D,1)$. Because the classifying map $\eta$ is cyclic,
the submanifold $\tilde N^{n-2k-2} \subset N^{n-2k}$ is
co-oriented, moreover
 we have
$$\eta^{\ast}(\psi_{\D}) = \mu_a^{\ast}(\psi_+),$$
where by $\psi_+$ we denote the 2-dimensional universal
$SO(2)$--bundle over $K(\I_a,1)$.

 Let us denote by
 $\tilde g: \tilde N^{n-2k-2} \looparrowright \R^n$ the
 restriction of the immersion $g$ on the submanifold  $\tilde N^{n-2k-2} \subset N^{n-2k}$,
assuming that the immersion $\tilde g$ is generic. The immersion
 $\tilde g$ is a $\D$--framed immersion by
$\tilde \Psi$, the classifying map  $\tilde \eta$ of this
$\D$--framed immersion is the restriction of $\eta$ to the
submanifold, this map  is cyclic. The triple $(\tilde g, \tilde
\eta, \tilde \Psi)$ is constructed from the triple  $(g, \eta,
\Psi)$ by means of the transfer homomorphism $J^{\D}$ in the
bottom row of the diagram $(\ref{1.5})$ (in this diagram $k_1$ is
changed to $k$, $k$ is changed to $k+1$).

Let us denote by  $\tilde L^{n-4k-4}$ the self-intersection
manifold of the immersion  $\tilde g$. The manifold  $\tilde
L^{n-4k-4}$ is a submanifold of the manifold  $L^{n-4k}$, $\tilde
L^{n-4k-4} \subset L^{n-4k}$. The parameterized immersion
 $\tilde h: \tilde L^{n-4k-4}
\looparrowright \R^n$ is well defined, this immersion is a
$\Z/2^{[3]}$--framed immersion by means of  $\tilde \Lambda$, the
classifying map $\tilde \zeta$ of this $\Z/2^{[3]}$--framed
immersion is quaternionic. The triple $(\tilde h, \tilde \zeta,
\tilde \Lambda)$ is defined from the triple  $(h, \zeta, \Lambda)$
by means of the homomorphism $J^{\Z/2^{[3]}}$ in the bottom row of
the diagram $(\ref{1.5.1})$ (in this diagram  $k_1$ is changed to
$k$,  $k$ is changed to $k+1$).

By Lemma  $\ref{lemma11}$ the submanifold $\tilde L^{n-4k-4}
\subset L^{n-4k}$ represents the Euler class of the bundle
$\zeta^{\ast}(\psi_{[3]})$. This submanifold is the source manifold of a $\Z/2^{[3]}$-immersion,
representing the image of the left bottom horizontal homomorphism in the diagram $(\ref{1.5.1})$
(in the diagram $k_1=k$, $k=k+1$).

Let us consider the canonical 2-sheeted covering
$\tilde p: \bar{\tilde L}^{4k-4} \to \tilde L^{n-4k}$. The
submanifold $\bar{\tilde L}^{4k-4} \subset \bar L^{n-4k}$
represents the Euler class of the bundle $\tilde
p^{\ast}(\zeta^{\ast}(\psi_{[3]}))$. This vector bundle is
naturally isomorphic to the vector bundle $\bar
\zeta^{\ast}(\psi^!_{[3]})$, where $\bar \zeta: \bar L^{n-4k} \to
K(\H_c,1)$ is the canonical 2-sheeted covering over the
classifying map $\zeta$ ($\H_c \cong \D \times \D$),
$\psi^!_{[3]}$ is the pull-back of the universal vector bundle
$\psi_{[3]}$ over $K(\Z/2^{[3]},1)$ by means of the covering
$K(\H_c,1) \to K(\Z/2^{[3]},1)$.

 Because the classifying map $\zeta$ is quaternionic,
 the submanifold $\tilde L^{n-4k-4} \subset
L^{n-4k}$ is co-oriented and represents the homological Euler
class of the $SO(4)$--bundle $\lambda^{\ast}(\psi_{\Q})$, and moreover for the corresponding $O(4)$--bundles 
$\bar
\zeta^{\ast}(\psi_{\H_c}) = \bar \lambda^{\ast}(\psi_{\Q}^!)$,
where:

--$\psi_{\Q}$ is the universal $SO(4)$--vector bundle over
the classifying space $K(\Q,1)$. This bundle is given by the quaternionic-conjugated representation 
with respect to the representation $(\ref{Q a1})$ -$(\ref{Q
a3})$. The bundle $\psi_{\Q}$, as a $O(4)$--bundle, is defined by the
formula: $\psi_{\Q} = i_a^{\ast}(\psi_{[3]})$, $\psi^!_{\Q} =
i_{\I_a,\Q}^{\ast}(\psi_{\Q})$.

-- $\psi_{\H_c}$ is the universal $O(4)$--bundle over $K(\H_c,1)$ ($\H_c
\cong \D \times \D$).

--$\bar \lambda: \bar L^{n-4k} \to K(\I_a,1)$ is the 2-sheeted
covering over the classifying mapping $\lambda: L^{n-4k} \to
K(\Q,1)$, induced by the 2-sheeted covering $K(\I_a,1) \to
K(\Q,1)$ over the target space of the map $\lambda$.

For the universal $SO(4)$--bundle   $ \psi_{\Q}^!$ the following formula is
satisfied:
\begin{eqnarray}\label{!}
\psi_{\Q}^! = \psi_+ \oplus  \psi_-, 
\end{eqnarray}
where the bundle $\psi_+$ admits a lift $\psi_+^{U}$ to a complex
$U(1)$--bundle, the bundle $\psi_-$ is a $SO(2)$--bundle, obtained
from $\psi_+^{U}$ by means of the complex conjugation and
forgetting the complex structure. 

The proof of $(\ref{!})$ follows from the formulas $(\ref{Q a1})$-$(\ref{Q a3})$. 
This formulas correspond to Lemma  $\ref{TQ}$: $\psi_- = T_{\Q}^{\ast}(\psi_+)$.

The bundles $\psi_+$, $\psi_-$
satisfy the equation: $e(\psi_+)=-e(\psi_-)$, and the Euler class
$e(\psi_+)$ of the bundle  $\psi_+$ is equal to the generator $t
\in H^2(K(\I_a,1);\Z)$ in the standard basis, the Euler class
$e(\psi_-)$ of the bundle $\psi_-$ is equal to $-t$ and is
opposite to the generator $t$ of the standard basis.

Let us denote by
 $m \in H^{4k}( N^{n-2k};\Z)$ the cohomology class, dual to
 the fundamental class of the oriented submanifold
 $\bar L^{n-4k} \subset N^{n-2k}$ in the oriented manifold
 $N^{n-2k}$. Let us denote by  $e_g \in
H^{4k}(N^{n-2k};\Z)$ the Euler class
of the immersion $g$ (this is the top class of the normal bundle $\nu_g$).
By the Herbert theorem for the immersion
$g: N^{n-2k} \looparrowright \R^n$ with the self-intersection
manifold   $L^{n-4k}$ (see [E-G], Theorem 1.1 the case  $r=1$,
the coefficients is $\Z$) the following formula are satisfied:
\begin{eqnarray}\label{1.501}
e_g+m=0.
\end{eqnarray}

Let us denote by $\tilde m \in H^{4k-4}( N^{n-2k};\Z)$
 the cohomology class, dual to
 the fundamental class of the oriented submanifold
 $\overline{\tilde L}^{n-4k-4} \subset
\tilde N^{n-2k-2} \subset N^{n-2k}$ in the oriented manifold
 $N^{n-2k}$. Let us denote by  $e_{\tilde g} \in H^{4k-4}(N^{n-2k};\Z)$ the cohomology class, dual to the
image of the homology Euler class of the
immersion $\tilde g$ by the inclusion $\tilde N^{n-2k-2} \subset
N^{n-2k}$.  By the Herbert theorem for the immersion $\tilde g:
N^{n-2k} \looparrowright \R^n$ with the self-intersection manifold
$\tilde L^{n-4k}$ (see. [E-G], Theorem 1.1 the case $r=1$, the
coefficients is $\Z$) the following formula are satisfied:
\begin{eqnarray}\label{1.502}
e_{\tilde g}+ \tilde m=0.
\end{eqnarray}

 Because
$\bar \lambda = \mu_a$, we may use the equation:
 $\bar \lambda^{\ast}(\psi_{\Q}^!) =
\mu_a^{\ast}(\psi_+) \oplus \mu_a^{\ast}(\psi_-)$. The following
equation are satisfied: $\tilde m = m e(\mu_a^{\ast}(\psi_+))
e(\mu_a^{\ast}(\psi_-))$, where the right side is the product of
the three cohomology classes: $m$ and the two Euler classes of the
corresponding bundles. The following equation are satisfied:
$e_{\tilde g} = e_g e^2(\mu_a^{\ast}(\psi_+))$. The equation
$(\ref{1.502})$ can be rewritten in the following form:
\begin{eqnarray}\label{1.503}
e_g e^2(\mu_a^{\ast}(\psi_+)) + m e(\mu_a^{\ast}(\psi_+))
e(\mu_a^{\ast}(\psi_-)) =0.
\end{eqnarray}

Then we may take into account $(\ref{1.501})$ and the equation
$e(\mu_a^{\ast}(\psi_-))=-e(\mu_a^{\ast}(\psi_+))$. Let us rewrite
the previous formula as follows:
\begin{eqnarray}\label{1.504}
2e_g e^2(\mu_a^{\ast}(\psi_+))=0.
\end{eqnarray}
Because of the equation $e_g = e(\mu_a^{\ast}(\psi_+))^k$, we
obtain:
\begin{eqnarray}\label{1.505}
2e^{k+2}(\mu_a^{\ast}(\psi_+))  =0.
\end{eqnarray}

Let us recall that  $\dim(L)=n-4k = n_{\sigma} \ge 7$ and
$\dim(\tilde L) = n_{\sigma}-4 \ge 3$. The formula for the Hopf
invariant for $\D$-framed immersion $(g,\eta,\Psi)$ using
$(\ref{2909})$ is the following:
\begin{eqnarray}\label{1.506}
 h^{\D}_k((g,\eta,\Psi)) =  \langle e^{k+2}(\mu_a^{\ast}(\psi_+))
\mu_a^{\ast}(\tau_{n-4k-4});[N^{n-2k}] \rangle \pmod{2},
\end{eqnarray}
where $\tau_{n-4k-4} \in H^{n-4k-4}(K(\I_a,1);\Z/4)$ is the
generic class  modulo 4, the cohomology class
$e(\mu_a^{\ast}(\psi_+))$ is modulo 4, and the fundamental class
$[N^{n-2k}]$ of the oriented manifold $N^{n-2k}$ is modulo 4. The
condition $ h^{\D}_k((g,\eta,\Psi))=1$ implies the following
condition: the cohomology class $e^{k+2}(\mu_a^{\ast}(\psi_+))$ is
of order 4. This contradicts the formula $(\ref{1.505})$.
Therefore, $h^{sf}_k(f,\kappa,\Xi)= h^{\D}_k((g,\eta,\Psi))=0$ and
the theorem in the particular case is proved.

Let us prove the theorem in the general case. Let us consider the
pair of mappings  $(\mu_a,
 \lambda)$, where
$\mu_a: N^{n-2k}_a \to K(\I_a,1)$, $N^{n-2k}_a \subset N^{n-2k}$, $\lambda = \lambda_{\Q} \cup
\lambda_{\H_{b \times \bb}} : L^{n-4k}_{\Q} \cup L^{n-4k}_{\H_{b \times \bb}} \to K(\Q,1)
\cup K(\H_{b \times \bb},1)$, where $L^{n-4k}_a = L^{n-4k}_{\Q} \cup
L^{n-4k}_{\H_{b \times \bb}}$, $L^{n-4k}_a \subset L^{n-4k}$, these two mappings determine the quaternionic
structure of the $\D$--framed immersion $(g,\eta,\Psi)$ in the sense of Definition $\ref{def16}$.

 Let us consider the manifold $\bar L^{n-4k}_a = \bar L^{n-4k}_{\Q} \cup \bar
L^{n-4k}_{\H_{b \times \bb}}$, defined by the formula $(\ref{barA})$. 
The manifold  $\bar L^{n-4k}_a$ is the
canonical 2-sheeted covering over the manifold $L^{n-4k}_a$.

The formula  $(\ref{1.501})$ is valid, and additionally the cohomology class
$m$ (this class is dual to the fundamental class  $[\bar L_a]$ of
the submanifold  $\bar L^{n-4k}_a \subset N^{n-2k}_a$)  decomposes
into the following sum:
\begin{eqnarray}\label{mm}
m=m_{\Q} + m_{\H_{b \times \bb}},
\end{eqnarray}
corresponding to the type of the components  $L^{n-4k}_{\Q}$,
$L^{n-4k}_{\H_{b \times \bb}}$ of the self-intersection manifold (see the formula $(\ref{hQH})$).

Let us consider the submanifold  $\tilde N^{n-2k-2}_a \subset
N^{n-2k}_a$, representing the Euler class of the bundle
 $\mu_a^{\ast}(\psi^+)$. The following immersion
 $\tilde g_a: \tilde N^{n-2k-2}_a \looparrowright
\R^n$ is well defined by the restriction of the immersion  $g_a$
to the submanifold $\tilde N^{n-2k-2}_a \subset N^{n-2k}_a$. Let us
denote by  $\tilde L^{n-4k-4}_a$ the self-intersection manifold of
the immersion $\tilde g_a$ (compare with the corresponding definition of the previous step).

The inclusion $\tilde L^{n-2k-4}_a \subset L^{n-2k}_a$ is
well defined. In particular, the manifold $\tilde L^{n-4k-4}_a$ is
represented by the union of the following two components: $\tilde
L^{n-4k-4}_a = \tilde L^{n-4k-4}_{\Q} \cup \tilde
L^{n-4k-4}_{\H_{b \times \bb}}$.

\begin{lemma}\label{600}

The co-oriented submanifold $\tilde L^{n-2k-4}_{\Q} \subset L^{n-2k}_{\Q}$
represents the Euler class of the $SO(4)$--bundle
$\lambda_{\Q}^{\ast}(\psi_{\Q})$.

The submanifold $\tilde L^{n-2k-4}_{\H_{b \times \bb}} \subset L^{n-2k}_{\H_{b \times \bb}}$
represents the Euler class of the $SO(4)$--bundle
$\lambda_{\H_{b \times \bb}}^{\ast}(\psi_{\H_{b \times \bb}})$, where $\psi_{\H_{b \times \bb}}$ is the
universal $SO(4)$--bundle over the space $K(\H_{b \times \bb},1)$. The corresponding $O(4)$--bundle is
standardly defined as the inverse image of the bundle
$\psi_{\Z/2^{[3]}}$ over $K(\Z/2^{[3]},1)$ by means of the
inclusion $K(\H_{b \times \bb},1) \subset K(\Z/2^{[3]},1)$.
\end{lemma}

\subsubsection*{Proof of Lemma $\ref{600}$}
The proof follows from the arguments above in the proof of
commutativity of the left squares of the diagrams $(\ref{1.5})$
and $(\ref{1.5.1})$. 
\[  \]

The bundle $\psi_{\H_{b \times \bb}}$ is isomorphic to the Whitney sum of the two
2-bundles: $\psi_{\H_{b \times \bb}} = p_{\H_{b \times \bb},\I_a}^{\ast}(\psi_{\I_a}) \oplus
p_{\H_{b \times \bb},\I_a}^{\ast}(\psi_{\I_a}) \otimes l_{\Z/2}$, where
$p_{\H_{b \times \bb},\I_a}^{\ast}(\psi_{\I_a})$ is the 2-dimensional bundle,
defined as the pull-back of the canonical 2-dimensional bundle
$\psi^+$ over $K(\I_a,1)$ by means of the natural mapping
$p_{\H_{b \times \bb},\I_a}: K(\H_{b \times \bb},1) \to K(\I_a,1)$, induced by the
homomorphism  $p_{\H_{b \times \bb},\I_a} :\H_{b \times \bb} \to \I_a$, $l_{\Z/2}$ is a line
bundle, defined as the inverse image of the canonical line bundle
over $K(\Z/2,1)$ by means of the projection  $K(\H_{b \times \bb},1) \to
K(\Z/2,1)$, this projection corresponds to the epimorphism  $\H_{b \times \bb}
\to \Z/2$ with the kernel  $\I_{a} \subset H_{b \times \bb}$.

By analogous arguments the class
 $\tilde{m}$ is well defined as in the formula $(\ref{1.502})$, moreover, the following
formula is satisfied:
\begin{eqnarray}\label{tildemm}
\tilde m=\tilde m_{\Q} +   \tilde m_{\H_{b \times \bb}},
\end{eqnarray}
where the terms in the right side of the formula are defined as
the cohomology classes, dual to the fundamental classes
$[\bar{\tilde L}_{\Q}]$, $[\bar{\tilde L}_{\H_{b \times \bb}}]$ of the
canonical coverings over the corresponding component.

The formula relating $m_{\Q}$ and $\tilde m_{\Q}$ is the following:
$\tilde m_{\Q} = m_{\Q} e(\mu_a^{\ast}(\psi_+))
e(\mu_a^{\ast}(\psi_-))$. The formula relating $m_{\H_{b \times \bb}}$ and
$\tilde m_{\H_{b \times \bb}}$ is the following: $\tilde m_{\H_{b \times \bb}} = 
m_{\H_{b \times \bb}}
e^2(\mu_a^{\ast}(\psi_+))$. To prove the last equation we use the
following fact: the bundle $i_{\I_a,\H_{b \times \bb}}^{\ast}(\psi_{\H_{b \times \bb}})$,
where the mapping $i_{\I_a,\H_{b \times \bb}}: K(\I_a,1) \to K(\H_{b \times \bb},1)$ corresponds to the index 2 subgroup  $i_{\I_a,\H_{b \times \bb}}:
\I_a \subset \H_{b \times \bb}$, is isomorphic to the bundle $\psi^+ \oplus
\psi^+$.

The analog of the formula $(\ref{1.503})$ is the following:
\begin{eqnarray}\label{1.5033}
e_g e^2(\mu_a^{\ast}(\psi_+)) - m_{\Q} e^2(\mu_a^{\ast}(\psi_+)) +
m_{\H_{b \times \bb}}e^2(\mu_a^{\ast}(\psi_+)) =0.
\end{eqnarray}

Let us multiply both sides of the formula $(\ref{mm})$ by
the cohomology class $e^2(\mu_a^{\ast}(\psi_+))$ and take the sum
with the opposite sign with $(\ref{1.5033})$, we get:
\begin{eqnarray}\label{1.5044}
2m_{\Q}e^2(\mu_a^{\ast}(\psi_+))=0.
\end{eqnarray}
This is an analog of the formula $(\ref{1.504})$.

Let us prove that the Hopf invariant of  the $\D$--framed
immersion
 $(g,\eta,\Psi)$ is trivial. By  Corollary $\ref{cor18}$
the Hopf invariant is given by the formula $(\ref{hQH})$. Let us
prove that the each term in this formula is equal to zero. The
first term $h_{\lambda}(L_{\Q})$, according to $(\ref{3024})$,
is calculated as the reduction modulo 2 of the following
characteristic number modulo 4:
$$ h_{\lambda}(L_{\Q}) = \langle m_{\Q} \mu_a^{\ast}(x);[N^{n-2k}_a] \rangle, $$
where $x \in H^{n-4k}(\I_a;\Z/4)$ is the generator. Analogously,
the second term $h_{\lambda}(L_{\H_{b \times \bb}})$ is the reduction modulo
2 of the following number modulo 4:
$$ h_{\lambda}(L_{\H_{b \times \bb}}) = \langle m_{\H_{b \times \bb}} \mu_a^{\ast}(x);[N^{n-2k}_a] \rangle. $$

Note that  $x=\tau^2y$, where $\tau \in H^2(K(\I_a,1);\Z/4)$, $y
\in H^{n-4k-4}(K(\I_a,1);\Z/4)$ are the generators. We have
$\mu_a(\tau)=e(\mu_a(\psi^+))$, because  $\tau$ is the Euler class
of the bundle $\psi_+$. Therefore, from $(\ref{1.5044})$ we get
$$ h_{\lambda}(L_{\Q}) = 0, $$
because
$m_{\Q}e^2(\mu_a^{\ast}(\psi_+))=m_{\Q}(\mu_a^{\ast}(\tau))^2 =
m_{\Q} \mu_a^{\ast}(x)$. 

To calculate the second term
$h_{\lambda}(L_{\H_{b \times \bb}})$ it is sufficient to note that
$$ \langle m_{\H_{b \times \bb}} \mu_a^{\ast}(x);[N^{n-2k}_{a}] \rangle =
\langle \mu_a^{\ast}(x);[\bar L^{n-4k}_{\H_{b \times \bb}}] \rangle = \langle
p^{\ast}((\lambda_{\H_{b \times \bb}})^{\ast}(x'));[\bar L^{n-4k}_{\H_{b \times \bb}}] \rangle =
0,$$ where $p: \bar L^{n-4k}_{\H_{b \times \bb}} \to L^{n-4k}_{\H_{b \times \bb}}$ -- is the
2-sheeted covering, corresponding to the subgroup
$i_{\I_a,\H_{b \times \bb}}$, $x' \in H^{n-4k}(K(\H_{b \times \bb},1);\Z/4)$ is a cohomology
class, such that $i^{\ast}_{\I_a,\H_{b \times \bb}}(x')=x$,
$[\bar L^{n-4k}_{\H_{b \times \bb}}]$ is the fundamental class of the total manifold
of the canonical 2-sheeted covering $p$.

Theorem $\ref{th9}$ is proved.

\begin{remark}
A straightforward generalization of Theorem
 $\ref{th9}$ for mappings with singularities  $c: S^{n-2k}/\i \to
\R^n$, which admits a (relative) quaternionic structure in the sense of Definition $\ref{quaterniond}$ is not possible. 
\end{remark}

 \[  \]
Moscow Region, Troitsk, IZMIRAN, 142190,
 pmakhmet@izmiran.ru
 \[  \]

\end{document}